\documentclass[10pt,reqno]{amsart} 

\usepackage{amssymb,latexsym,amsmath,amsthm,amsfonts,amscd}
\usepackage{multirow,comment,ulem}
\usepackage{mathtools}
\usepackage{color}
\usepackage{parskip}
\usepackage{marginnote}

\begin{document}

\title{Rankin-Selberg $L$-functions for $\GSpin \times \GL$ Groups}

\author[Mahdi Asgari]{Mahdi Asgari} 
\address{Department of Mathematics \\ 
Oklahoma State University \\ 
Stillwater, OK 74078--1058 \\
USA} 
\email{asgari@math.okstate.edu} 

\author[James Cogdell]{James W. Cogdell} 
\address{Department of Mathematics \\ 
Ohio State University \\ 
Columbus, OH 43210 \\
USA} 
\email{cogdell@math.osu.edu} 
	 
\author[Freydoon Shahidi]{Freydoon Shahidi} 
\address{Mathematics Department \\  
Purdue University \\ 
West Lafayette, IN 47907 \\ 
USA} 
\email{shahidi@math.purdue.edu} 

\numberwithin{equation}{section}
\newtheorem{thm}{Theorem}[section]
\newtheorem{cor}[thm]{Corollary}
\newtheorem{lem}[thm]{Lemma}
\newtheorem{prop}[thm]{Proposition}
\newtheorem{con}[thm]{Conjecture}
\newtheorem{ass}[thm]{Assumption}
\newtheorem{defi}[thm]{Definition}
\newtheorem{exer}[thm]{Exercise}

\theoremstyle{remark}
\newtheorem{rem}[equation]{Remark}
\newtheorem{notation}[equation]{Notation}
\newtheorem{exam}[equation]{Example}

\newcommand{\A}{\mathbb A}
\newcommand{\GL}{\mathrm{GL}}
\newcommand{\GO}{\mathrm{GO}}
\newcommand{\GSO}{\mathrm{GSO}}
\newcommand{\GSp}{\mathrm{GSp}}
\newcommand{\GSpin}{\mathrm{GSpin}}
\newcommand{\SO}{\mathrm{SO}}
\newcommand{\Sp}{\mathrm{Sp}}
\newcommand{\Spin}{\mathrm{Spin}}

\newcommand{\C}{\mathbb C} 
\newcommand{\Q}{\mathbb Q} 
\newcommand{\Z}{\mathbb Z} 
\newcommand{\R}{\mathbb R} 

\newcommand{\cW}{\mathcal W} 

\renewcommand{\L}{\mathcal L} 
\newcommand{\W}{{\mathcal W}} 
\newcommand{\X}{\mathcal X} 
  
\newcommand{\bs}{\backslash}
\newcommand{\im}{{\rm im}}
\newcommand{\pr}{{\rm pr}} 
\newcommand{\sym}{{\rm Sym}} 
\newcommand{\tr}{\operatorname{tr}} 

\newcommand{\bpm}{\begin{pmatrix}}
\newcommand{\epm}{\end{pmatrix}}
\renewcommand{\i}{\mathfrak i} 
\renewcommand{\Re}{{\rm Re}}

 \newcommand{\adots}{\mathinner{\mkern1mu\raise1pt\hbox{.} 
\mkern2mu\raise4pt\hbox{.} 
\mkern2mu\raise7pt\hbox{.}\mkern1mu}} 

\newcommand{\e}{\ensuremath{\epsilon}}

\newcommand{\w}[1]{\ensuremath{\widetilde{#1}}} 
\newcommand{\h}[1]{\ensuremath{\widehat{#1}}} 

\begin{abstract}
We construct an integral representation for the global Rankin-Selberg (partial) $L$-function $L(s, \pi \times \tau)$ 
where $\pi$ is an irreducible globally generic cuspidal automorphic representation of a general spin group 
(over an arbitrary number field) and $\tau$ is one of a general linear group, generalizing the works 
of Gelbart, Piatetski-Shapiro, Rallis, Ginzburg, Soudry and Kaplan among others. We consider all ranks 
and both even and odd general spin groups including the quasi-split forms. The resulting facts about the location 
of poles of $L(s, \pi \times \tau)$ have, in particular, important consequences in describing the image of 
the Langlands funtorial transfer from the general spin groups to general linear groups. 
\end{abstract}

\maketitle

%%%%%%%%%%%%%% NEW SECTION %%%%%%%%%%%%%%%%%%%% 
%%%%%%%%%%%%%%%%%%%%%%%%%%%%%%%%%%%%%%%%%%%%%%% 
\section{Introduction}\label{sec:Intro}
The purpose of this article to develop {\sl integral representations} for the Rankin-Selberg convolution $L$-functions 
for the globally generic automorphic representations of $\GSpin \times \GL$ groups. We consider all semi-simple ranks of both 
the general spin groups and the general linear groups, including the case of the quasi-split, non-split, even, general spin groups. 
Our work in this article corresponds to steps (1), part of (2), and (4) in the ``$L$-function machine'' as described in \cite[\S 1.7]{gs}: 
establish a ``Basic Identity'' for the global zeta integrals and expand them as an Euler product; analyze the meromorphic 
behavior while we do not deal with the functional equation in this article; and complete the ``unramified computation'' which 
relates the local zeta integrals to local Langlands $L$-functions in the unramified case.

For us the main motivation for this work was its application to classifying the image of the generic functorial transfer 
from the general spin groups to the general linear groups, even though our results are otherwise of interest as well, just as it has 
been the case with all the cases of developing integral representations for automorphic $L$-functions. Indeed in \cite{manuscripta} 
the first and third authors already used some analytic properties of the partial Rankin-Selberg $L$-function $L^S(s, \pi \times \tau)$ 
for the description of the image of the generic functorial transfer from $\GSpin_{2n+1}$ or $\GSpin_{2n}$ to $\GL_{2n}$. 
Here, $\pi$ is a globally generic, unitary, cuspidal automprhic representation of the general spin group and $\tau$ is one of 
the general linear group. As usual, $S$ denotes a finite set of places, including the Archimedean places, outside of which all the 
data is unramified. See Proposition \ref{L-pole} for the precise statement.

Our results here are restricted to the globally generic representations, and their application is to the functorial transfer of the globally generic 
automorphic representations. However, in recent years there have been great strides in establishing functorial transfer of 
arbitrary (not necessarily generic) automorphic representations from the classical groups and their similitude versions, as well as unitary 
groups, to the general linear groups. In particular we mention J. Arthur's endoscopic classification of representations of 
special orthogonal and symplectic groups \cite{arthur} and the works of Y. Cai, S. Friedberg, D. Ginzburg, D. Gourevitch 
and E. Kaplan on global functoriality for non-generic representations \cite{cfk-ppt, cfk-gafa, cfgk-invent, cfgk-pams, gk23}.  
Arthur's book does not consider the case of general spin groups, even though Arthur's methods, with appropriate modifications, 
would be applicable to the general spin groups. On the other hand, the works of Y. Cai et al mentioned above for the non-generic 
representations do cover the general spin groups.  The generic transfer for the general spin groups was established earlier 
by the first and third named authors \cite{duke, manuscripta}.

S. Gelbart, I. I. Piatetski-Shapiro and S. Rallis gave the original methods of constructing integral representations that 
produce the Rankin-Selberg $L$-functions for $G \times \GL_n$ in \cite{gpsr}, where $G = \SO_{2n+1}$, $\SO_{2n}$, 
or $\Sp_{2n}$.  The integrals in each case look different and they needed substantially different methods to deal with 
each, calling them {\sl Method A}, {\sl Method B}, and {\sl Method C}, respectively. As we only deal with groups of 
Dynkin types $B$ and $D$, we will focus on the first two methods, which we briefly recall below.

In Method A, one takes a globally generic cuspidal automorphic representation $\pi$ of $\SO_{2n+1}(\A)$ and a cuspidal representation 
$\tau$ of $\GL_n(\A)$.  We consider $\GL_n$ as the Levi factor of the Siegel parabolic subgroup in $\SO_{2n}$ (i.e., ``doubling 
the number of variables'') and construct an Eisenstein series on it. One then embeds $\SO_{2n}$ in $\SO_{2n+1}$ and integrates 
a Whittaker function of $\pi$ against a Fourier coefficient of the Eisenstein series. Here, the integral is over the adelic points 
of $\SO_{2n}$ modulo its rational points. Gelbart, Piatetski-Shapiro and Rallis then study this integral by ``unfolding'' it and writing 
it as an Euler product. They then compute the local integral at an unramified finite place $v$, which will turn out to be expressed 
as a quotient of the local $L$-function $L(s, \pi_v \times \tau_v)$ and the exterior square $L$-function $L(s, \tau, \wedge^2)$. 
They prove this by considering the decomposition of a certain symmetric algebra and use some results of Ton-That \cite{ton-that1,ton-that2} 
along with the Caselman-Shalika formula. A similar construction can be done for $\SO_{2n} \times \GL_{n-1}$ as well, with an embedding 
of $\SO_{2(n-1)+1}$ inside $\SO_{2n}$, with the symmetric square $L$-function $L(s, \tau, \sym^2)$ replacing the exterior square $L$-function.

In Method B, the roles of cuspidal representation $\pi$ and the Eisenstein series constructed 
from $\tau$ are switched, so the cuspidal representation is on the smaller group $\SO_{2n}$ while the Eisenstein series is on the 
larger group $\SO_{2n+1}$, coming from an induced representation from the Siegel parabolic of $\SO_{2n+1}$. While the analysis 
in Method B is somewhat different, a similar unramified computation can be done. The result will be again the local Rankin-Selberg 
$L$-function $L(s, \pi \times \tau)$, divided by the symmetric square $L$-function $L(s, \tau, \sym^2)$.  Again, one can also 
consider the case of $\SO_{2n} \times \GL_n$ using Method B and again the exterior square $L$-function $L(s, \tau, \wedge^2)$ 
appears. While \cite{gpsr} 
mostly focuses on the split groups, they point out that the methods work for quasi-split groups as well, and even double covers of 
special orthogonal groups, i.e., the spin groups. (They also cover the Rankin-Selberg construction for $\Sp_{2n} \times \GL_n$ 
in their Method C as we mentioned above.)

D. Ginzburg \cite{ginzburg} generalized Method A from $\SO_{2n+1} \times \GL_n$, resp. $\SO_{2n} \times \GL_{n-1}$, to the case 
of $\SO_{2n+1} \times \GL_m$ with $m \le n$, resp. $\SO_{2n} \times \GL_m$ with $m \le n-1$.  The idea here is that one proceeds 
similarly as above, using $\tau$ on $\GL_m$ to construct an Eisenstein series on $\SO_{2m}$ and embeds $\SO_{2m}$ inside $\SO_{2m+1}$. 
However, one then  ``pads" the  integral with some unipotent integrals in order to produce the zeta integral that gives the Rankin-Selberg 
$L$-function for $\SO_{2n+1} \times \GL_m$. It is clear that $m \le n$ is necessary for this process to work. Again, there is also a similar 
procedure for $\SO_{2n} \times \GL_m$ with $m \le n-1$.  In order to do the unramified computations, Ginzburg uses an inductive 
argument that reduces the proof to that of the case of $m=n$, resp. $m=n-1$.  Fortunately, the results of Ton-That mentioned above 
for decomposing the symmetric algebra is still available for any rank and one could replace the inductive argument with the direct  
decomposition the symmetric algebra in this case, as we do in our work for the general spin groups (see below).

However, when one considers generalizing Method B for $\SO_{2n+1} \times \GL_m$ with $m > n$, a serious obstacle 
appears because the decomposition of the symmetric algebra appears to be much more complicated and not so easily tractable. 
D. Soudry \cite{soudry-mem} showed how to resolve this issue by finding a certain ``duality'' between the Rankin-Selberg $L$-function  
for $\SO_{2n+1} \times \GL_m$ with $m > n$ (in Method B) and that for $\SO_{2m} \times \GL_n$ with $m > n$ (in Method A). 
E. Kaplan \cite{kaplan-thesis} then extended this work to the case of $\SO_{2n} \times \GL_m$ with $m \ge n$. Kaplan also 
considers the quasi-split, non-split, forms that exist in the case of even special orthogonal groups.

In this article we give a construction of an integral representations for Rankin-Selberg $L$-functions $L(s, \pi \times \tau)$ where 
$\pi$ is an irreducible globally generic (i.e., having a Whttaker model with respect to a generic character) cuspidal automorphic representation 
of a general spin groups and $\tau$ is an irreducible cuspidal automorphic representation of a general linear group. 
See Thoerem \ref{thm-BasicIdentityI} and Theorem \ref{thm-BasicIdentityII}.  
We consider 
both the odd and the even cases, including the quasi-split non-split forms in the even case, and any rank of the general linear group. As such, 
we are generalizing all the works above, both in Method A and Method B. See the table in Section \ref{sec-euler-exp} for the details 
of the various cases.

Naturally, we follow similar constructions as in the works of Ginzburg, Soudry and Kaplan (and the original works of Gelbart, Piatetski-Shapiro and Rallis). 
The main difference in the similitude case, in addition to a careful analysis of the unfolding arguments, is the appearance of the 
{\it twisted} symmetric/exterior square $L$-functions of $\tau$. See Theorems \ref{unram-A} and \ref{unram-B} for the details.  
We also point out that in the quasi-split non-split case the expressions obtained from the unramified computations in the two theorems  
resemble the ones in the ``opposite  parity'' case.  This phenomenon is to be expected considering the Galois action on the Dynkin diagram 
in the even case, for example, and it is already present in the work of E. Kaplan for quasi-split $\SO_{2n}$.

As we mentioned above, for the symmetric algebra decomposition we use the results of Ton-That \cite{ton-that1,ton-that2}, which are 
available for the special orthogonal and symplectic groups.  We carefully study the effect of the presence of the center in the similitude case, 
cf. Section \ref{symm-alg-subsubsec}. 
We then combine this with a suitable version of the Casselman-Shalika formula that we review in Section \ref{cs-subsec} in order to 
complete the unramified computations.  As we already mentioned we then use ideas similar to \cite[Appendix]{gpsr} to show directly 
that the two sides of the equations we are claiming as the result of our unramifed computation given the same power series in $q^{-s}$. 
It should be possible to accomplish the same goal by arguing inductively as in the works of Ginzburg or Kaplan. 
As pointed above, there is  substantially difference analysis in Methods A and B and therefore the unramifed computations also look  
significantly different. That is why we do the two methods in two separate theorems in Section \ref{loc-id-sec} even though the final 
expressions for the local unramified integrals look similar in Theorems \ref{unram-A} and \ref{unram-B}.

In order to effect the duality argument of Soudry, mentioned above, for our cases we need to use some results about the $\gamma$-factors 
for the groups involved. The constructions and analysis of the $\gamma$-factors is part of the ``$L$-function machine" that we have not 
studied in this article. However, fortunately E. Kaplan, J. F. Lau and B. Liu have studied them for exactly the cases we need in \cite{kll}.  
As such we have simply used their result in the only place where we need to invoke the $\gamma$-factors in this article, namely the 
generalization of Soudry's duality argument.  Finally, we point out that in the quasi-split case in Method B, we also need to invoke a 
certain uniqueness result which is fortunately also provided by \cite{kll} as part of their work on local descent from general spin groups.

We note that $\GSpin_5 \cong \GSp_4$ and therefore our results, in particular, cover the Rankin-Selberg product 
$L$-functions for $\GSp_4 \times \GL_1$, $\GSp_4 \times \GL_2$ and $\GSp_4 \times \GL_3$ for generic representations. 
These cases have been studied extensively over the years. 
The twist by $\GL_1$ essentially amounts to the construction of the standard $L$-function of $\GSp_4$ while the 
twist by $\GL_2$ was studied by Novodvorsky \cite[\S 3]{novod} and Soudry \cite{soudry-duke}. Bump \cite[\S 3.3--3.5]{bump} 
surveys these two cases as well as the twist by $\GL_3$, where he gives a particular embedding of $\GL_3$ in $\GL_4$ and assumes 
that the representations have trivial central character.   We point out that the split group $\GSpin_6$ is isomorphic to 
$\left\{ (A,b) \in \GL_4 \times \GL_1 :  \det(A) = b^2 \right\}$ (cf. \cite[\S 2.2]{asgari-choiy-forum}).  
As such, the subgroup of $\GL_4$ that Bump works with for this case 
corresponds to one of the two (isomorphic, non-conjugate) Siegel Levi subgroups in $\GSpin_6$ and our construction would then 
agree with Bump's description in \cite[\S 3.5]{bump}.  
(There are also many more works for $\GSp_4 \times \GL_2$ that consider non-generic representations, which our results 
here do not cover.)

Recently P. Yan and Q. Zhang have considered a Rankin-Selberg integral for a general linear group and a product of two general linear 
groups in \cite{yan-zhang}.  Their study, in particular, gives another proof of Jacquet's local converse theorem. Since the small rank case 
of $\GL_2 \times \GL_2$ is very close to the group $\GSpin_4$, their integral for $(\GL_2 \times \GL_2) \times \GL_n$ and our case of 
$\GSpin_4 \times \GL_n$ are closely related.

{\it Acknowledgements.} 
We would like to thank D. Ginzburg, J. Hundley, E. Kaplan, J. F. Lau, E. Sayag, and D. Soudry for helpful discussions and/or their 
interest in this work. In particular, we thank Eyal Kaplan for very helpful conversations, explaining his work to us, and for a careful 
reading of an earlier draft of this paper. Also, the first author 
thanks Pan Yan for some very helpful discussions.

%%%%%%%%%%%%%% NEW SECTION %%%%%%%%%%%%%%%%%%%% 
%%%%%%%%%%%%%%%%%%%%%%%%%%%%%%%%%%%%%%%%%%%%%%% 
\tableofcontents

%%%%%%%%%%%%%% NEW SECTION %%%%%%%%%%%%%%%%%%%% 
%%%%%%%%%%%%%%%%%%%%%%%%%%%%%%%%%%%%%%%%%%%%%%% 
\section{The Preliminaries} \label{sec-pre}
%%%%%%%%%%%%%%%%%%%%%%%%%%%%%%%%%%
\subsection{Notation} 
Let $k$ be a number field and let $\A=\A_k$ be the ring of ad\`eles of $k.$ 
Also, let $F$ be a local field of characteristic zero. Often we have $F=k_v$ for 
some place $v$ of the number field $k$.

We consider both odd and even 
(split and quasi-split) 
general spin groups defined over $k$ 
or $F$. 
We will 
use $G = G_{n'} = \GSpin_{n'}$ and $H = H_{m'} = \GSpin_{m'}$ with either 
\[ \begin{array}{lll}
\mbox{(case A)} & n' = 2n+1,  m' = 2m, & \mbox{ or} \\  
\mbox{(case B)} & n' = 2n,  m' = 2m+1 
\end{array} \]
with $m$ and $n$ positive integers. 

Later on we will assume that $n' < m'$ and introduce an embedding 
$G \hookrightarrow H$ with $n'$ and $m'$ of opposite parity.

As in \cite[Theorem 4.3.1]{hs} there exist surjections $G_{n'} \longrightarrow \SO_{n'}$ and we  
fix one such surjection and denote it by $\pr,$ so that we have 
\begin{equation} \label{pr}
 \pr : G_{n'} \longrightarrow \SO_{n'}. 
\end{equation}
The projection map also gives maps at the level of $k$-points, $F$-points, 
and the ad\`elic points, all of which will also be 
denoted $\pr$.

%%%%%%%%%%%%%%%%%%%%%%%%%%%%%%%%%%
\subsection{Structure of $\GSpin$ Groups}\label{G-odd-even} 
There are several constructions one could give for the general spin groups. One construction is via the 
introduction of a based root datum for each group as in \cite[\S 7.4.1]{spr}, 
which we do below. 
More detailed descriptions can also 
be found in \cite[\S 2]{duke} and \cite[\S 4]{hs}.

%%%%%%%%%%%%%%%%%%%%%%%%%%%%%%%%%%
\subsubsection{The root data of  $\GSpin$ groups} 
Let $n' \ge 3$. (See Remark \ref{low-rank} below.) 
The based root datum of the split 
$\GSpin_{n'}$ is given by $(X,R,\Delta, X^{\vee},R^{\vee}, \Delta^\vee)$, where 
$X$ and $X^{\vee}$ are $\mathbb Z$-modules generated by generators $e_{0},e_{1},\cdots,e_{n}$ 
and $e^{*}_{0},e^{*}_{1},\cdots,e^{*}_{n}$, respectively. The pairing 
\begin{equation}\label{pairing} 
\langle \ , \ \rangle : X \times X^{\vee} \longrightarrow \mathbb Z
\end{equation}
is the standard pairing. 

When $n' = 2n+1$ the roots and coroots are given by 
\begin{eqnarray}
\label{R-odd}
R = R_{2n+1} &=& \left\{ \pm (e_{i} \pm e_{j}) \, : \, 1 \le i<j \le n \right\} \cup 
\left\{\pm e_{i} \, : \, 1\le i \le n \right\} \\
\label{R-odd-co}
R^{\vee} = R^{\vee}_{2n+1} &=& \left\{ \pm(e^{*}_{i} - e^{*}_{j}) \, : \, 1 \le i<j \le n \right\} \cup \\ 
&&\nonumber 
\left\{ \pm(e^{*}_{i} + e^{*}_{j} - e_0^*) \, : \, 1 \le i<j \le n \right\} \cup 
\left\{\pm (2e^{*}_{i} - e^{*}_{0}) \, : \, 1\le i \le n \right\} 
\end{eqnarray}
along with the bijection $R \longrightarrow R^{\vee}$ given by 
\begin{eqnarray}
\label{v-odd-1}
(\pm(e_{i} - e_{j}))^{\vee} &=&  \pm (e^{*}_{i} - e^{*}_{j}) \\ 
\label{v-odd-2}
(\pm(e_{i} + e_{j}))^{\vee} &=&  \pm (e^{*}_{i} + e^{*}_{j} - e^{*}_{0}) \\ 
\label{v-odd-3}
(\pm e_{i})^{\vee} &=& \pm (2e^{*}_{i} - e^{*}_{0}).  
\end{eqnarray} 
Moreover, we fix the following choice of simple roots and coroots: 
\begin{eqnarray}
\Delta &=& \left\{ e_1-e_2, e_2-e_3,\cdots,e_{n-1}-e_n, e_n \right\}, \\ 
\Delta^\vee &=& \left\{ e^*_1-e^*_2, e^*_2-e^*_3,\cdots,e^*_{n-1}-e^*_n, 2e^*_n-e^*_0 \right\}.
\end{eqnarray} 
This based root datum determines the group $\GSpin_{2n+1}$ uniquely, equipped with a Borel subgroup 
$B$ containing a maximal torus $T.$

When $n' = 2n$ we have 
\begin{eqnarray}
\label{R-even}
R = R_{2n} &=& \left\{ \pm (e_{i} \pm e_{j}) \, : \, 1 \le i < j \le n \right\}  \\
\label{R-even-co}
R^{\vee} = R^{\vee}_{2n} &=& \left\{ \pm (e^{*}_{i} - e^{*}_{j}) \, : \, 1 \le i<j \le n \right\} \cup \\ 
\nonumber
&& \left\{ \pm (e^{*}_{i} + e^{*}_{j} - e^{*}_{0}) \, : \, 1 \le i<j \le n \right\} 
\end{eqnarray}
along with the bijection $R \longrightarrow R^{\vee}$ given by 
\begin{eqnarray}
\label{v-even-1}
(\pm(e_{i} - e_{j}))^{\vee} &=&  \pm (e^{*}_{i}-e^{*}_{j}) \\ 
\label{v-even-2}
(\pm (e_{i} + e_{j}))^{\vee} &=& \pm (e^{*}_{i} + e^{*}_{j} - e^{*}_{0}).  
\end{eqnarray} 
and 
\begin{eqnarray}
\Delta &=& \left\{ e_1-e_2, e_2-e_3,\cdots,e_{n-1}-e_n, e_{n-1}+e_n \right\}, \\ 
\Delta^\vee &=& \left\{ e^*_1-e^*_2, e^*_2-e^*_3,\cdots,e^*_{n-1}-e^*_n, e^{*}_{n-1}+e^{*}_{n}-e^{*}_{0} \right\}.
\end{eqnarray} 
Again, this based root datum determines the split group $\GSpin_{2n}$ uniquely, equipped with a Borel subgroup 
$B$ containing a maximal torus $T.$

\begin{rem} \label{low-rank} 
We should also mention that $\GSpin_0 \cong \GL_1$,  $\GSpin_1 \cong \GL_1$ and $\GSpin_2 \cong \GL_1 \times \GL_1$. 
While some of the notation above make sense for these small rank cases as well, $\Delta$ and $\Delta^\vee$ are 
empty.  Finally, for the quasi-split non-split even general spin group $\GSpin_2^*$, 
associated with a quadratic extension $K/k$ (see below), we have 
$\GSpin_2^* \cong \operatorname{Res}_{K/k} \GL_1$. 
\end{rem}

%%%%%%%%%%%%%%%%%%%%%%%%%%%%%%%%%%
\subsubsection{Abstract group structure of $\GSpin$ groups}
We proved in \cite{asgari} that the above is equivalent to a second construction of the split  $\GSpin_{n'}$ given 
as a suitable quotient of $\GL_{1} \times \Spin_{n'},$ where $\Spin_{n'}$ is the split, simply-connected, simple, 
connected group of type $B_n$ if $n' = 2n+1,$ or of type $D_n$ if $n' = 2n.$ If $n' \ge 3,$ we have 
\begin{equation} 
\GSpin_{n'} \cong \left( \GL_{1} \times \Spin_{n'} \right) \slash \left\{(1,1), (-1,c) \right\}, 
\end{equation} 
where $c$ is a particular element of the center of $\Spin_{n'}$ as follows: 
\begin{itemize}
\item[(A)] If $n' = 2n+1,$ then $Z(\Spin_{n'}) = \{1, c\} \cong \Z \slash 2 \Z,$   
\item[(B1)] If $n' = 2n$ with $n$ even, then $Z(\Spin_{n'}) = \{1, c, z, cz\} \cong \Z \slash 2 \Z \times \Z \slash 2 \Z,$ and 
\item[(B2)] If $n' = 2n$ with $n$ odd, then $Z(\Spin_{n'}) = \{1, z, c=z^2, z^3\} \cong \Z \slash 4 \Z,$  
\end{itemize}
where $c = \alpha^\vee_{n}(-1)$ when $n' = 2n+1$ and $c = \alpha^\vee_{n-1}(-1) \alpha^\vee_{n}(-1)$ when 
$n' = 2n,$ an element of order 2 in all cases.

%%%%%%%%%%%%%%%%%%%%%%%%%%%%%%%%%%
\subsubsection{Quasi-split forms of even $\GSpin$} \label{sec:qsplit}
When $n' = 2n,$ we also have the quasi-split $\GSpin_{n'}$ groups, which we describe below.  
They are of type $^{2}D_{n}.$ We refer to \cite{hs} and \cite{cpss-clay} for more details about them. 

Let $G$ be a quasi-split form of $\GSpin_{2n}.$  We have a fixed Borel subgroup $B$ and a 
Cartan subgroup $T \subset B$. We fix a {\it pinning} (or splitting)  
$(B, T, \{x_{\alpha}\}_{\alpha\in\Delta})$, where $\{x_{\alpha}\}$ is a collection of root vectors, 
one for each simple root of $T$ in $B$. 
We also denote the maximal $k$-split subtorus of $T$ by $T_s.$ 
Following \cite{hs} we give the following two parametrizations of the quasi-split forms of $\GSpin_{2n}.$

\  
%%%%%%
\paragraph{\it First Parametrization.}  

By \cite[\S16.2]{spr}, the quasi-split forms of $\GSpin_{2n}$ over $k$ are determined by the indexed root data 
$(X,\Delta, X^\vee, \Delta^\vee, \Delta_0,\nu).$  Here, $\Delta_0$ is empty since the 
group is quasi-split. Also, $\nu$ denotes a Galois action on $X$ and $X^\vee.$  The Galois action is 
either trivial or switches the simple roots $e_{n-1}-e_n$ and $e_{n-1}+e_n$ while keeping all other simple roots 
fixed. In fact, the nontrivial Galois element acts on $X$ and on $X^\vee$ via 
\begin{equation} 
\nu(e_i) = 
\begin{cases} 
e_0 + e_n, & i=0, \\
e_i, & 1 \le i \le n-1, \\
-e_n, & i=n,
\end{cases}
\quad \mbox{ and } \quad
\nu(e^*_i) = 
\begin{cases} 
e^*_0, & i=0, \\
e^*_i, & 1 \le i \le n-1, \\
-e^*_n + e^*_0, & i=n.
\end{cases}
\end{equation}
Moreover, the $k$-rational character lattice $_kX$ is spanned by $e_1,\dots,e_{n-1}, e_n+2 e_0$ and the $k$-rational 
cocharacter lattice $_kX^\vee$ is spanned by $e^*_0,e^*_1,\dots,e^*_{n-1}$ (cf. \cite[\S 4.3]{hs}.)  In particular, the 
root system of $G$ relative to $T_s$ is of type $B_{n-1}$ with $k$-rational simple roots and $k$-rational simple coroots given by 
$_k\Delta = \left\{e_1 - e_2, \dots, e_{n-2} - e_{n-1}, e_{n-1}\right\}$ and 
$_k\Delta^\vee = \left\{e^*_1 - e^*_2, \dots, e^*_{n-2} - e^*_{n-1}, 2 e^*_{n-1} - e^*_0\right\}.$

There is a one-to-one correspondence between 
\begin{itemize}
\item[(i)] the quasi-split $k$-groups $G$ with connected component of $L$-group 
${}^LG^0 \cong \GSO_{2n}(\C)$ 
\item[] and 
\item[(ii)] the characters $\mu : \operatorname{Gal}(\overline{k} / k) \longrightarrow S,$ where 
\begin{equation*} 
S = \left\{ \sigma \in \operatorname{Aut}(X(T)) : \sigma \mbox{ permutes $\Delta$ via an 
automorphism of the Dynkin diagram} \right\}. 
\end{equation*}
\end{itemize} 
We have $S \cong \Z \slash 2 \Z$ and by class field theory the characters of order two of 
$\operatorname{Gal}(\overline{k} / k)$ are in bijection with 
\begin{itemize}
\item[(iii)] the quadratic characters $\mu : k^\times \backslash \A^\times \longrightarrow \left\{\pm 1 \right\}.$ 
\end{itemize}
Therefore, the quasi-split forms of $\GSpin_{2n}$ are parametrized by the quadratic idele class characters of $k.$ 
When $\mu$ is nontrivial, we denote the associated quasi-split non-split group by $\GSpin^\mu_{2n}$ 
or simply $\GSpin^*_{2n}$ when the particular $\mu$ is unimportant. We will also denote the 
quadratic extension of $k$ associated with $\mu$ by $K^\mu/k$ or simply $K/k$.

\ 
%%%%%%
\paragraph{\it Second Parametrization.}  For $a \in k^\times,$ denote its square class 
in $k^\times \slash (k^\times)^2$ by $\underline{a} = a (k^\times)^2$. 
Let $k(\sqrt{a})$ be the smallest extension of $k$ in which the elements of 
$\underline{a}$ are squares, so that $k(\sqrt{\underline{a}}) = k(\sqrt{a}).$ We then let 
$\GSpin^{\underline{a}}_{2n}=\GSpin^{a}_{2n}$  denote the quasi-split form of $\GSpin_{2n}$ such that the 
associated map $\operatorname{Gal}(\overline{k} / k) \longrightarrow \{\pm 1\}$ factors through 
$\operatorname{Gal}\left(k(\sqrt{\underline{a}})\right / k) = \operatorname{Gal}\left(k(\sqrt{a})/k\right)$.

%%%%%%%%%%%%%%%%%%%%%%%%%%%%%%%%%%
\subsubsection{Dual of $\GSpin$ Groups} \label{duals}
Yet another construction for $\GSpin_{n'}$ is via their dual groups as follows \cite[\S 4.3]{hs}: 
\begin{itemize}
\item[(A)] If $n' = 2n+1,$ then the split 
$G=\GSpin_{n'}$ is the $k$-split, connected, reductive group having 
the based root datum dual to $\GSp_{2n}$ (so of type $B_n$).  
(See \cite[\S 2.3]{duke} for the precise description.)    
Hence, ${}^LG = \GSp_{2n}(\C) \times  \operatorname{Gal}(\overline{k} / k)$, with the 
Galois group acting trivially.

\item[(B)] If $n' = 2n,$ then the split 
$G=\GSpin_{n'}$ is the $k$-split connected, reductive group having 
based root datum dual to that of $\GSO_{2n}$ (so of typle $D_n$). 
Here, $\GSO_{2n}$ is the connected component of the group $\GO_{2n}$ with all groups defined over $k.$ 
(Again, see \cite[\S 2.3]{duke} for the precise description.)    
Hence, ${}^LG = \GSO_{2n}(\C) \times  \operatorname{Gal}(\overline{k} / k)$, with the 
Galois group acting trivially.

\item[(C)] 
If $n' = 2n$,  $G = \GSpin^a_{2n}$, the quasi-split group associated with $K = k(\sqrt{a})$ as above, 
and its $L$-group can be given by 
${}^LG = \GSO_{2n}(\C) \rtimes  \operatorname{Gal}(\overline{k} / k)$, a semi-direct 
product where the Galois action on $\GSO_{2n}(\C)$ is given as follows.  
If $\gamma \in \operatorname{Gal}(\overline{k} / k)$ with $\gamma\vert_K = 1$, then  
the action is trivial. If $\gamma\vert_K \not= 1$, then the action of $\gamma$ on $g \in \GSO_{2n}(\C)$ 
is given by $j g j^{-1}$, where $j = \operatorname{diag}(I_{n-2}, \operatorname{diag}(w,w), I_{n-2})$ and 
$w = \begin{pmatrix}0 & 1 \\ 1 &0 \end{pmatrix}$.  

\end{itemize}

%%%%%%%%%%%%%%%%%%%%%%%%%%%%%%%%%%
\subsection{Weyl Groups} 
By \cite[Lemma 6.2.1]{hs}, we know that the Weyl group of $G_{n'} = \GSpin_{n'}$ is isomorphic to 
the Weyl group of $\SO_{n'}.$

When $n' = 2n,$ for the split $G_{2n}$, as for $\SO_{2n},$ we have that the Weyl group 
$W_{n'} \cong \mathfrak S_n \rtimes \{ \pm 1\}^{n-1}.$  Choose representatives  
$(p,\underline{\e}) \in \mathfrak S_n \rtimes \{ \pm 1\}^{n-1}$ as in \cite[\S 6]{hs}. 
Similarly, when $n' = 2n+1,$ for the split $G_{2n+1},$ as for $\SO_{2n+1},$ we have that the Weyl group 
$W_{n'} \cong \mathfrak S_n \rtimes \{ \pm 1\}^n.$ Again, take representatives 
$(p, \underline{\e}) \in \mathfrak S_n \rtimes \{ \pm 1\}^n$. 
In both cases, we have the following Weyl actions on the root and coroot lattices of $G_{n'}$: 
\begin{equation*}
(p,\underline{\e}) \cdot e_i = 
\begin{cases} 
e_{p(i)} & i>0, \e_{p(i)} = 1, \\ 
-e_{p(i)} & i>0, \e_{p(i)} = -1, \\ 
e_0 + \sum\limits_{\e_{p(i)=-1}} e_{p(i)} & i=0, 
\end{cases}
\end{equation*}
and 
\begin{equation*}
(p,\underline{\e}) \cdot e^*_i = 
\begin{cases} 
e^*_{p(i)} & i>0, \e_{p(i)} = 1, \\ 
e^*_0 - e^*_{p(i)} & i>0, \e_{p(i)} = -1, \\ 
e^*_0  & i=0.  
\end{cases}
\end{equation*}

The above assertions for $n'$ even are proved in \cite[Lemma 6.2.3]{hs}. The assertions in 
the odd case were intended to be the content of \cite[Lemma 13.2.2]{hs} although 
it appears that Lemma 13.2.2 and its proof were copies of Lemma 6.2.3 and its proof without modification. 
However, a similar calculation gives the above.

%%%%%%%%%%%%%% NEW SECTION %%%%%%%%%%%%%%%%%%%% 
%%%%%%%%%%%%%%%%%%%%%%%%%%%%%%%%%%%%%%%%%%%%%%% 
\section{Unipotent Periods, Parabolic Subgroups, and Embeddings}  \label{sec-embed} 

%%%%%%%%%%%%%%%%%%%%%%%%%%%%%%%%%%
\subsection{Unipotent Periods} 
We recall the following facts from \cite[\S 4,11]{hs}: 
\begin{itemize}
\item The kernel of the projection map $\pr$ in (\ref{pr}) lies in the center $Z(G_{n'}).$ In fact, 
if we write 
\begin{equation*} 
1 \longrightarrow \left\{(1,1), (-1,c) \right\} \longrightarrow \GL_{1} \times \Spin_{n'} \longrightarrow G_{n'} \longrightarrow 1, 
\end{equation*}
then $\ker(\pr)$ is the image of the $\GL_{1}$ factor and so it is central. From this it follows that the action of $G_{n'}$ on itself 
by conjugation factors through $\pr.$

\item If $u$ is a unipotent element of $G_{n'}(\A)$ and $g \in G_{n'}(\A),$ then $\pr(gug^{-1})$ is unipotent in 
$\SO_{n'}(\A)$ and $gug^{-1}$ is the unique unipotent element of its preimage in $G_{n'}(\A).$

\item $\pr : G_{n'} \longrightarrow \SO_{n'}$ induces an isomorphism of unipotent varieties. We may specify 
unipotent elements of subgroups by their images under $\pr.$ This defines coordinates for any unipotent element 
or subgroup. Hence, we may write $u_{i,j}$ for the $(i,j)$-entry of $\pr(u).$ In particular, unipotent periods in $G_{n'}$ 
and $\SO_{n'}$ can be identified. Therefore, any identity or relationship between unipotent periods in $\SO_{n'}$ which 
is proved only by conjugation or ``swapping'' (root exchange) extends to $G_{n'}.$ 
\end{itemize}

Following the notation of \cite{hs}, if $G$ is any reductive algebraic group defined over $k,$ $U$ is a unipotent 
subgroup of $G,$ and $\psi_U$ is a character of $U(k) \backslash U(\A),$ we define 
\begin{equation} \label{u-period} 
\varphi^{(U,\psi_U)} (g) = \int\limits_{U(k) \backslash U(\A)} \varphi(ug) \psi^{-1}_U(u) \, du. 
\end{equation} 
This unipotent integral is a Fourier coefficient (cf. \cite{grs, soudry-mem}). Here, one can take $\varphi$ to be an automorphic form on $G.$

%%%%%%%%%%%%%%%%%%%%%%%%%%%%%%%%%%
\subsection{Parabolic Subgroups} 

There are two types of parabolic subgroups that play a role.

%%%%%%%%%%%%%%%%%%%%%%%%%%%%%%%%%%
\subsubsection{The Maximal Parabolic Subgroups $P_\ell$} 
For 
$1 \le \ell < n$ if $n'= 2n+1$  or  
$1 \le \ell < n-1$ if $n' = 2n,$ 
let $P_\ell$ be the standard maximal parabolic subgroup of $G_{n'}$ with Levi isomorphic to 
$\GL_{\ell} \times \GSpin_{n'-2\ell}$.  

If $\gamma\in \GL_\ell(k)$ with $k$ a field, we let $\gamma^\wedge$ denote 
the image of $\gamma$ under the isomorphism of $\GL_{\ell} \times \GSpin_{n'-2\ell}$ 
with the Levi component of $P_\ell$.  (Note that the character and cocharacter lattices of 
$\GSpin_{n'-2\ell}$ are sublattices of $X$ and $X^\vee$, spanned by the 
generators $e_{0},e_{\ell+1},\cdots,e_{n}$ and $e^{*}_{0},e^{*}_{\ell+1},\cdots,e^{*}_{n}$, 
respectively.) 
Since there is a corresponding parabolic subgroup of $SO_{n'}$ we can write
\[
\gamma^\wedge=\begin{pmatrix} \gamma & & \\ & I_{n'-2\ell} \\ & & \gamma^*\end{pmatrix}\in SO_{n'}(k)
\]
which we can identify with its lift to $G_{n'}(k)$ if $\gamma$ is a unipotent or a Weyl group element.

%%%%%%%%%%%%%%%%%%%%%%%%%%%%%%%%%%
\subsubsection{The Siegel Parabolic Subgroup $P_n$} \label{subsec:siegel} 
This is a standard maximal parabolic subgroup with Levi isomorphic to $\GL_{n} \times \GL_{1}.$ 
When $G_{n'}$ is split with $n' = 2n,$ there are two parabolic subgroups $P = M U$ with Levi subgroup 
$M \cong \GL_{n} \times \GL_{1}.$  Following \cite[\S 6]{hs}, we denote by $P_n$ the one in which one 
deletes the root $e_{n-1}+e_n$ and the coroot $e^*_{n-1}+e^*_{n}-e^*_0.$ 

If we consider $G_{n'}$ with $n' = 2n+1,$ then there is a unique standard parabolic subgroup $P = M U$ with 
$M \cong \GL_{n} \times \GL_{1}.$ We denote this parabolic by $P_n$ as well. 

In either case, the subgroup $\mathfrak S_n \subset W_{n'}$ is isomorphic to $W_M,$ the Weyl group 
of $M.$ Also, $\pr\left(\GL_{n} \times \GL_{1}\right) = \GL_{n}$ and $\ker(\pr) = \im (e^*_0).$

%%%%%%%%%%%%%%%%%%%%%%%%%%%%%%%%%%
\subsubsection{The Parabolic Subgroups $Q_\ell$} \label{subsec:Q}
For $G_{n'}$ split with $n' = 2n,$ and for $1 \le \ell < 
n-1,
$ we let $Q_\ell = L_\ell N_\ell$ denote 
the standard parabolic subgroup with Levi  
\begin{equation*} 
L_\ell \cong \GL_{1}^\ell \times G_{2n - 2\ell} 
\end{equation*} 
and unipotent radical given by 
\begin{equation*} 
N_\ell = \left\{ u \, \vert \, u_{i,j} = 0 \mbox{ for } i > \ell \mbox{ with } i < j 
\le 2n 
-i \right\}. 
\end{equation*} 

For $n' = 2n+1$ and $1 \le \ell < n,$ we let $Q_\ell = L_\ell N_\ell$ be the parabolic subgroup with Levi 
\begin{equation*} 
L_\ell \cong \GL_{1}^\ell \times G_{2n+1 - 2\ell} 
\end{equation*} 
and unipotent radical given by 
\begin{equation*} 
N_\ell = \left\{ u \, \vert \, u_{i,j} = 0 \mbox{ for } i > \ell \mbox{ with } i < j 
\le 2n+1 
-i \right\}. 
\end{equation*} 
%

%%%%%%%%%%%%%%%%%%%%%%%%%%%%%%%%%%
\subsubsection{Stabilizers} 
First, consider the parabolic $Q_\ell \subset G_{n'}$ with $n' = 2n$ even. In this case, a general character 
of $N_\ell$ is of the form  (cf. \cite[Remark 9.1.2]{hs}) 
\begin{equation*} 
\psi_0 \left(c_1 u_{1,2} + \cdots + c_{\ell-1} u_{\ell-1,\ell} + d_1 u_{\ell,\ell+1} + \cdots + d_{2n-2\ell} u_{\ell,2n-\ell} \right),  
\end{equation*} 
where $\psi_0$ is a fixed non-trivial additive character (of $k$, $\A$, or $F$, depending on the context).

The Levi $L_\ell$ acts on the space of characters and over an algebraically closed field there is an open orbit, 
consisting of all those elements with $c_i \not= 0$ for all $i$ and 
$\underline{d} \, J \, {}^t\!\underline{d} \not= 0,$ where 
$\underline{d} = \left(d_1, d_2, \dots, d_{2n-2\ell} \right)$ 
and $J$ is the matrix with $1$'s on the skew diagonal and 
zeros elsewhere. Over a general field $k$ two such elements are in the same $k$-orbit if and only if the values of 
$\underline{d} \, J \, {}^t\!\underline{d}$ 
are in the same square class.

Let $\Psi_\ell$ be the character of $N_\ell$ defined by 
\begin{equation*} 
\Psi_\ell(u) = \psi_0\left(u_{1,2} + \cdots + u_{\ell-1,\ell} + u_{\ell,n} - u_{\ell,n+1} \right). 
\end{equation*} 
Then one can see (\cite[Remark 9.1.2]{hs}) that 
\begin{itemize}
\item the stabilizer $L_\ell^{\Psi_\ell}$ has two connected components, 
\item the connected component of the identity $\left( L_\ell^{\Psi_\ell} \right)^0 \cong G_{2n-2\ell-1},$ 
\item there is an ``obvious'' choice of isomorphism $\iota : G_{2n-2\ell-1} \longrightarrow \left( L_\ell^{\Psi_\ell} \right)^0$ 
having the following property: If $\left\{e_i^* : i=0,1,\dots, n\right\}$ is the basis for the cocharacter lattice of $G_{2n},$ 
and $\left\{\bar{e}_i^* : i=0,1,\dots, n-\ell-1\right\}$ is the basis for $G_{n-\ell-1},$ then 
\begin{equation*} 
\iota \circ \bar{e}_i^* = \begin{cases} 
e_0^* & i=0, \\ 
e_{\ell+1}^* & i=1,\dots, n-\ell-1. \\
\end{cases}
\end{equation*} 
\end{itemize}  
Note that it follows from $\iota \circ \bar{e}_0^* = \bar{e}_0^*$ that the induced map between 
$\bar{e}_0^* \left(\GL_1\right)$ in $G$ and $\bar{e}_0^* \left(\GL_1\right)$ in $H$ is the identity (and not, for example, 
the inversion map).

Next, consider $Q_\ell \subset G_{n'}$ with $n' = 2n+1$ odd. In this case, the general character of $N_\ell$ is of the form 
\begin{equation*} 
\theta(u) = 
\psi_0 \left(c_1 u_{1,2} + \cdots + c_{\ell-1} u_{\ell-1,\ell} + d_1 u_{\ell,\ell+1} + \cdots + d_{2n+1-2\ell} u_{\ell,2n+1-\ell} \right). 
\end{equation*} 
The Levi $L_\ell$ acts on the space of characters and over an algebraically closed field there is an open orbit, consisting of 
all those elements with $c_i \not= 0$ for all $i$ and 
$\underline{d} \, J \, {}^t\!\underline{d} \not= 0,$ where 
$\underline{d} = \left(d_1, d_2, \dots, d_{2n+1-2\ell} \right)$ 
and $J$ is the matrix with $1$'s on the skew diagonal and 
zeros elsewhere. Over a general field $k$ two such elements are in the same $k$-orbit if and only if the two values of 
$\underline{d} \, J \, {}^t\!\underline{d}$ 
are in the same square class. 
With $\theta$ as above, define $\operatorname{Inv}(\theta)$ to be the square class of 
$\underline{d} \, J \, {}^t\!\underline{d}.$ 
The character $\theta$ is said to be in general position if $c_i \not= 0$ for $i=1,\dots, \ell-1$ and $\operatorname{Inv}(\theta) \not= 0.$

By \cite[Lemma 16.1.7]{hs} we have that if $\theta$ is in general position, then its stabilizer in $L_\ell,$ namely 
$L_\ell^{\theta},$ has two components and $\left( L_\ell^{\theta} \right)^0 \cong G_{2n-2\ell}^{\operatorname{Inv}(\theta)}.$

For $a \in k^\times$ let $\Psi_\ell^a$ be the character of $N_\ell$ defined by 
\begin{equation*} 
\Psi_\ell^a(u) = 
\psi_0 \left(u_{1,2} + \cdots + u_{\ell-1,\ell} + u_{\ell,n} + \frac{a}{2} u_{\ell,n+2} \right). 
\end{equation*} 
Then the orbit of $\Psi_\ell^a$ is determined by the square class of $a.$  The character 
$\theta$ is in the same 
orbit as $\Psi_\ell^1.$  For each $a \in k^\times,$ we have $\left(L_\ell^{\Psi_\ell^a}\right)^0 \cong G_{2n-2\ell}^{\underline{a}},$ 
where $\underline{a}$ is the square class of $a.$ 

We note here that \cite{hs} states the above assertions only in the case of even $n$ because only this case is needed there. However, 
the assertions are valid for all $n.$

%%%%%%%%%%%%%%%%%%%
\begin{rem}[{\it Notaion for the parabolic subgroups of \cite{grs}}] \label{par-grs-acs} 
We note that \cite{grs} uses a different notation for its parabolic subgroups.  However, while we have switched the notation for the 
parabolics to match both \cite{hs} and what seems to us more standard notation, the notation for the unipotent radicals do agree with 
those of \cite{grs}.  As this caused confusion at one point, let us make this explicit. 

In \cite{grs}, the parabolic $P_\ell$ is that preserving a maximal isotropic flag of length $\ell$ \cite[p. 42]{grs}.  
Its Levi decomposition is $P_\ell = M_\ell N_\ell$ with $M_\ell \cong (\GL_1)^\ell \times h(W_{m,\ell})$ in their notation.  
They then define a character $\psi_\ell$ or $\psi_{\ell,\alpha}$ which agrees with our $\Psi_\ell^a$ and then have 
$\operatorname{Stab}_{h(W_{m,\ell})} (\psi_{\ell,\alpha}) \cong G.$  Therefore,  
\begin{equation}
P_\ell^{GRS} = Q_\ell^{ACS} \quad M_\ell^{GRS} = L_\ell^{ACS} \quad N_\ell^{GRS} = N_\ell^{ACS}.
\end{equation}

In \cite{grs} the parabolic $Q_\ell$ is that preserving a maximal isotropic subgroup of dimension $\ell$ \cite[pp. 65--66]{grs}.  
Its Levi decomposition is $Q_\ell = D_\ell U_\ell$ with $D_\ell \cong \GL_\ell \times h(W_{m,\ell})$ in their notation \cite[p. 81]{grs}.  
Therefore, 
\begin{equation}
Q_\ell^{GRS} = P_\ell^{ACS} \quad D_\ell^{GRS} = M_\ell^{ACS} \quad U_\ell^{GRS} = U_\ell^{ACS}.
\end{equation}
Here, the labels GRS refer to the subgroups in \cite{grs} and the labels ACS refer to the corresponding ones in this paper. 
\end{rem}

%%%%%%%%%%%%%% NEW SECTION %%%%%%%%%%%%%%%%%%%% 
%%%%%%%%%%%%%%%%%%%%%%%%%%%%%%%%%%%%%%%%%%%%%%% 
\section{Global Integrals I (Case B)} \label{sec-integralB} 

In this section we translate \cite[\S 10.3]{grs} into the context of $\GSpin$ groups.  
This corresponds to 
Method B  in the original work of Gelbart and Piatetski-Shapiro \cite[Part B]{gpsr}, which dealt with the special 
orthogonal and general linear groups with equal (or nearly equal) ranks.  As such, we refer 
to the integrals of this section as Case B. (See the table in Section \ref{sec-euler-exp} for a summary of various  
cases.) Also, recall Remark \ref{par-grs-acs} on our notation for the parabolic subgroups.

We begin with a number field $k$ and a $k$-vector space $V$ of dimension $\dim V = m' \ge 3,$ where 
$m'$ can be even or odd.  We take a non-degenerate quadratic form on $V$  and let $h(V) = \SO(V) = \SO_{m'}$ 
denote the special orthogonal group of $V.$  We let $H = H_{m'} = \GSpin_{m'}$ be the associated 
$\GSpin$ group that covers $h(V),$ so associated to a quadratic space of dimension $m'.$  We assume that $h(V)$ and $H$ are split so that 
$m = [\frac{1}{2} \dim V]$ is the Witt index of $V.$  (This is Assumption 2.1 of \cite{grs}.) Therefore, either $m' = 2m$ or $m' = 2m+1.$  Recall that 
we have fixed a projection 
\begin{equation*}
\pr : H_{m'} \longrightarrow h(V) 
\end{equation*}
that induces an isomorphism of unipotent varieties and allows us to identify unipotent periods on $H$ and $h(V)$. (See \cite{hs}.)

Following \cite{hs}, we take an integer $\ell$ such that $1 \le \ell < m$ and let $Q_\ell \subset H$ be the standard 
parabolic subgroup of $H,$ where $Q_\ell = L_\ell N_\ell$ with $L_\ell \cong \GL_{1}^\ell \times H_{m' - 2\ell}$ and 
the unipotent radical $N_\ell$ is as in Section \ref{subsec:Q}.  We consider the following characters of $N_\ell:$ 
\begin{equation*}
\Psi_\ell (u) = \psi_0 \left(u_{1,2} + \cdots + u_{\ell-1,\ell} + u_{\ell,2m} - u_{\ell,2m+1}\right)  \mbox{ if } m' = 2m, 
\end{equation*}
or 
\begin{equation*}
\Psi_\ell^a = \psi_0 \left(u_{1,2} + \cdots + u_{\ell-1,\ell} + u_{\ell,2m} + \frac{a}{2} u_{\ell,2m+2}\right)  \mbox{ if } m' = 2m+1, 
\end{equation*}
where $a \in k^\times.$  We let 
\begin{equation*}
G = \operatorname{Stab}_{L_\ell} \left( \Psi_\ell \right)^0 \cong \GSpin_{2m-2\ell-1} \mbox{ if } m' = 2m,  
\end{equation*}
by \cite[Section 9.1]{hs} or 
\begin{equation*}
G = \operatorname{Stab}_{L_\ell} \left( \Psi_\ell^a \right)^0 \cong \GSpin^a_{2m-2\ell} \mbox{ if } m' = 2m+1  
\end{equation*}
by \cite[Lemma 16.1.7]{hs}, which allows for the case of $G$ being quasi-split.  In what follows, we will simply use $\Psi_\ell$ in either case 
and suppress the $a \in k^\times.$ 

In \cite{hs}, the authors do not consider the case $\ell = 0$ which would correspond to the construction of 
Gelbart and Piatetski-Shapiro \cite[Part B]{gpsr}.
When $\ell = 0$ there is no $\Psi_\ell$ and we just ``restrict''.  
We are able to ``pull back'' 
the embedding in \cite[p. 43]{grs}. The embedding proceeds as follows. When $\ell=0$ we take 
\[
w_0=y_a=e_{\tilde{m}}+(-1)^{m'+1}\frac{a}{2}e_{-\tilde{m}}\in SO_{m'} 
\]
in the notation of \cite{grs}. Then $y_a\in W_{m, \tilde{m}-1}$ and hence in $W_{m, \ell}$ for all $\ell$ including $\ell=0$. 
We have $(y_a,y_a)=(-1)^{m'+1}a$ with $a\in k^\times$. Then we take $G$ to be $pr^{-1}(y_a^\perp)\subset H$. In this case, 
\[
 G\cong \GSpin_{2m-1}\mbox{ if }m'=2m
\]
 or 
\[
G\cong\GSpin^a_{2m}\mbox{ if }m'=2m+1.
\]
(Note in \cite{grs} in the case $B_{m'}$  they take a vector $\ell_0$ such that $V=X\oplus\langle \ell_0\rangle\oplus X^\vee$ with $(\ell_0,\ell_0)=1$. 
The vector $y_a$ should play the same role as their $\ell_0$. Note that in either case, the larger group $H$ is of type $B_{m'}$.) The statements below  now hold in the $\ell=0$ case.

Let $P_m = M_m U \subset H$ denote the Siegel parabolic subgroup of $H.$  We made a choice of $P_m$ 
in Section \ref{sec-embed} and we have $M_m \cong \GL_{m} \times \GL_{1}$.  We fix the isomorphism between $\GL_{m} \times \GL_{1}$ 
and $M_m$ as in \cite[\S 1]{cfk-ppt}, where the isomorphism is denoted by $i_M$.  (The choice of this isomorphism matters later, such as 
in \eqref{fs-I}, where we use $\omega^{-1}$, not $\omega$.)  With this isomorphism fixed, we have 
$\pr : M_m \longrightarrow \GL_{m}$ with 
$\ker(\pr) = \im(e_0^*) = Z(H)^0.$  Let $\tau$ be a cuspidal automorphic representation of $\GL_{m}(\A)$ 
and let $\eta$ be an idele class character of $\GL_{1}(\A) = \A^\times.$  We form the following 
induced representation (with $s$ being replaced by $s-1/2$ later on below): 
\begin{equation*} 
\rho = \rho_{\tau,\eta,s} = \operatorname{Ind}_{P_m(\A)}^{H(\A)} \left( \tau |\det|^{s} \otimes \eta \right). 
\end{equation*}  
For $\Re(s) \gg 0$ we can choose a section $f_s \in V_\rho,$ the space of $\rho,$ form the Eisenstein series 
\begin{equation*}
E(h, f_s) = \sum\limits_{\delta \in P_m(k) \bs H(k)} f_s(\delta h),
\end{equation*}  
continue as 
a function of $s,$ and form, as in (\ref{u-period}), the unipotent period 
\begin{equation*} 
E^{(N_\ell,\Psi_\ell)} (g, f_s), 
\end{equation*}  
 which is naturally an automorphic form on 
 $G(\A) = \left( \operatorname{Stab}_{L_\ell} \Psi_\ell \right)^0(\A).$ Note that in the $\ell=0$ case we interpret as above, i.e., 
 $G\hookrightarrow H$ and simple restriction from $H$ to $G$.

Take $G \cong \GSpin_{m'-2\ell-1}$ and $H \cong \GSpin_{m'}$ as above with the embedding $G \hookrightarrow H.$ 
Let $(\pi , V_\pi)$ be a cuspidal automorphic representation of $G(\A)$ with central character $\omega_\pi$ and denote by $\omega$ 
the idele class character such that for $a\in\A^\times$, we have 
\begin{equation} \label{omega} 
\pi(e^*_0(a)) = \omega(a) \operatorname{Id}_{V_\pi}. 
\end{equation} 
Let $(\tau , V_\tau)$ be a cuspidal automorphic representation of $\GL_{m}(\A).$ Let $f_s$ be a $K$-finite 
section (with $K$ denoting the maximal compact at the Archimedean place) in the space of 
\begin{equation}\label{fs-I} 
\rho_{\tau,\omega^{-1},s-1/2} = 
\operatorname{Ind}_{P_m(\A)}^{H(\A)} \left( \tau |\det|^{s - 1/2} \otimes \omega^{-1} \right). 
\end{equation} 
For $\varphi \in V_\pi,$ a cusp form on $G(\A),$ consider 
\begin{equation*} 
\L(\varphi,f_s) = \int\limits_{Z(\A) G(k) \bs G(\A)} \varphi(g) E^{(N_\ell,\Psi_\ell)} (g,f_s) \, dg, 
\end{equation*} 
where $Z = Z(G)^0$ is the identity component of the center of $G$. 
(When $\ell = 0$ we would just restrict and there is no unipotent period.)

%%%%%%%%%%%
\begin{lem} \label{Siegel-estimates}
The integral $\L(\varphi,f_s)$ converges absolutely and uniformly in vertical strips in $\C$ away from poles of the Eisenstein series. 
Therefore, it defines a meromorphic function on $\C.$ 
\end{lem}
%%%%%%%%%%%

\begin{proof}
This is essentially basic estimates on Siegel sets in $G.$  It is verified in the $\SO$ case in \cite{grs} and 
the same proof works for $\GSpin.$ 
\end{proof}

%%%%%%%%%%%
\begin{thm} \label{thm-BasicIdentityI}
Let the notation be as above. 
\begin{itemize}
\item[(i)] Assume that $\L(\varphi,f_s)$ is not identically zero.  Then $\pi$ is globally generic (for a suitable $\psi$-Whittaker model).  

\item[(ii)] If $\Re(s) \gg 0,$ then we have an identity 
\begin{equation*} 
\L(\varphi,f_s) = \int\limits_{N_G(\A)Z(\A) \bs G(\A)} W_\varphi^{\psi}(g) 
\int\limits_{N_\ell(\A) \cap \beta^{-1} P_m(\A) \beta \bs N_\ell(\A)} f_s^{(Z_m,\psi)} \left( \beta u g \right) \Psi_\ell(u)^{-1}
\, du \, dg, 
\end{equation*} 
where 
\begin{itemize} 
\item[$\bullet$] $N_G \subset G$ is a standard maximal unipotent subgroup of $G,$ 
\item[$\bullet$] $\psi = \psi_{N_G}$ is the standard Whittaker character on $N_G$ obtained from $\psi_0$ as in \cite[p. 289]{grs},  
\item[$\bullet$] $W_\varphi^\psi(g)$ is the corresponding $\psi$-Whittaker function of $\varphi,$
\item[$\bullet$] $\beta = \beta_{\ell,a}$ is a product of a certain Weyl group element and a (rational) ``diagonal'' element in 
$H,$ or more precisely, in a unipotent subgroup $J \subset N_\ell$ described in Lemma \ref{group-J} below,   
\item[$\bullet$] 
$Z_m$ is the 
maximal unipotent subgroup of $\GL_{m} \subset P_m \subset H$ and $\psi$ is the standard Whittaker character of $N_m,$ and 
\item[$\bullet$] we have 
\begin{equation*} 
f_s^{(Z_m,\psi)}(h) = \int\limits_{Z_m(k) \bs Z_m(\A)} f_s(\hat{z}h) \psi^{-1}(z) \, dz,
\end{equation*}
with $\hat{z}$ the image of $z \in \GL_{m}$ in $H.$ 
\end{itemize}  
\end{itemize} 
\end{thm}
%%%%%%%%%%%

%%%%%%%%%%%
\begin{proof}
We prove this in the $\ell\neq 0$ case. The $\ell=0$ case follows as in \cite{gpsr}. The proof involves several steps as we detail below, establishing several intermediate results along the way.

\noindent{\bf Step 1.} 
The integral converges for all $s,$ because 
\begin{itemize}
\item $\varphi$ is rapidly decreasing (mod $Z$), 
\item $E$ is of moderate growth (mod $Z$), 
\item $E^{(N_\ell,\Psi_\ell)}$ is a compact integration. 
\end{itemize}

\noindent{\bf Step 2.}   
For $\Re(s) \gg 0,$ i.e., in the realm of absolute convergence of the Eisenstein series, we replace $E$ by its absolutely convergent series. Then, 
\begin{equation*} 
E^{(N_\ell,\Psi_\ell)} (h,f_s) = 
\int\limits_{N_\ell(k) \bs N_\ell(\A)} 
\sum\limits_{\delta \in P_m(k) \bs H(k)}
f_s(\delta u h) \Psi_\ell^{-1}(u) \, du. 
\end{equation*} 

\noindent{\bf Step 3.}   
We factor the sum through the double cosets of $P_m \bs H / P_\ell,$ where $P_\ell$ is the maximal parabolic subgroup 
of $H$ with Levi isomorphic to $\GL_{\ell} \times \GSpin_{m' - 2\ell}.$

%%%%%%%%%%%
\begin{lem}
The coset representatives for $P_m(k) \bs H(k) / P_\ell(k)$ are as in \cite[pp. 70-71]{grs}, or \cite[p. 285]{grs}, i.e.,  
$\epsilon_r \in W_H = W_{\SO_{m'}}$ for $0 \le r \le \ell,$ with  
\begin{equation*} 
\pr(\epsilon_r) = 
\begin{pmatrix} I_r &&\\ & 0 & I_{m-\ell} \\ & I_{\ell-r} & 0 \end{pmatrix}^\wedge
\begin{pmatrix} I_r &&&&\\ &&& I_{\ell-r} & \\ && I_{m'-2\ell} && \\ & I_{\ell-r} & &&& \\ &&&& I_r \end{pmatrix}
w_b^{\ell-r}, 
\end{equation*} 
where the $\wedge$ notation is defined in 3.2.1, taken for the group $H$,  and where $w_b$ is as in \cite[pp. 70--71]{grs},  so an auxiliary Weyl group element. 
\end{lem}

%%%%%%%%%%%
\begin{proof} 
We know that $W_H \cong W_{\SO_{m'}}$ and that $\ker(\pr) \subset Z(H) \subset Z(T),$ where $T$ is the maximal 
torus in $H.$  From \cite{grs} we know that 
\begin{equation*} 
\SO_{m'} = \bigsqcup\limits_{r=0}^\ell \pr\left( P_m(k) \right) \pr(\epsilon_r) \pr\left( P_\ell(k) \right). 
\end{equation*} 
Since $\ker(\pr) \subset Z(H) \subset P_m,$ we have
\begin{equation*} 
H(k) = \bigsqcup\limits_{r=0}^\ell P_m(k) \epsilon_r P_\ell(k). 
\end{equation*} 
\end{proof}

Using this decomposition, we can partially unfold the Eisenstein series: 
\begin{equation*} 
E^{(N_\ell,\Psi_\ell)} (h, f_s) = 
\sum\limits_{r=0}^\ell 
\int\limits_{N_\ell(k) \bs N_\ell(\A)} 
\sum\limits_{\delta \in P_\ell^{(r)}} 
f_s\left( \epsilon_r \delta u h \right) \Psi_\ell^{-1}(u) \, du, 
\end{equation*} 
where $P_\ell^{(r)} = P_\ell \cap \epsilon_r^{-1} P_m \epsilon_r.$

The authors of \cite{grs} make $P_\ell^{(r)}$ explicit in their equation (4.19)--(4.20) or (4.22)--(4.23).  
If we write $P_\ell^{(r)} = M_\ell^{(r)} \ltimes U_\ell^{(r)},$ then the unipotent part $U_\ell^{(r)}$ agrees 
with \cite{grs} and the $\GSpin_{m'} / \SO_{m'}$ difference is in $M_\ell^{(r)}.$  
We will need the Levi part of this decomposition when $r = 0.$

\noindent{\bf Step 4.}  
Now, we factor the innermost sum through the double cosets $P_\ell^{(r)}\backslash P_\ell / R_\ell,$ where 
$R_\ell = R_{\ell,a} = \operatorname{Stab}_{L_\ell} \left( \Psi_\ell^a \right) \ltimes N_\ell = G \ltimes N_\ell.$ 
The authors of \cite{grs} compute the representatives for these double cosets for $\SO_{m'}.$  
In the case of $\SO_{m'}$ they are of the form   
\begin{equation*} 
\begin{pmatrix} \epsilon && \\ & \gamma & \\ &&\epsilon^* \end{pmatrix}
\end{equation*} 
with $\epsilon$ running through a set of representatives for $W_{\GL_{r}\times\GL_{\ell-r}} \bs W_{\GL_{\ell}}$ and 
\begin{equation*} 
\gamma = \begin{cases} 
I_{m'-2\ell} & \mbox{ if } m' = 2m, \mbox{ or } m' = 2m+1 \mbox{ and } a \not= t^2, \\ \\
\begin{pmatrix} 
I_{m'-\ell-1} &&&& \\ 
& 1 &&& \\ 
& \pm t & 1 && \\ 
& - \frac{t^2}{2} & \mp t & 1 & \\ 
&&&& I_{m' - \ell - 1}
\end{pmatrix}
& \mbox{ if } m' = 2m+1 \mbox{ and }  a = t^2. \\ 
\end{cases}
\end{equation*} 
Note that these are either Weyl group representatives, which we choose representatives for in $H$, 
or Weyl group representatives times unipotents, which have unique lifts to $H.$  So we obtain the following.

%%%%%%%%%%%
\begin{lem} 
The representatives $\{ \eta \}$ for $P_\ell^{(r)}(k) \bs P_\ell(k) / R_\ell(k)$ are either 
Weyl group representatives or Weyl group representatives times unipotents and so are uniquely 
determined by 
\begin{equation*} 
\pr(\eta) = 
\begin{pmatrix} 
\epsilon && \\ & \gamma & \\ && \epsilon^*
\end{pmatrix} 
\in \SO_{m'}(k)
\end{equation*} 
as above.
\end{lem}

Then we can further unfold the Eisenstein series 
\begin{eqnarray*} 
E^{(N_\ell,\Psi_\ell)}(h,f_s) =&  \\
\sum\limits_{r=0}^{\ell} \, 
\sum\limits_{\eta \in P_\ell^{(r)}(k) \bs P_\ell(k) / R_\ell(k)}^{}  \, \, 
\int\limits_{N_\ell(k) \bs N_\ell(\A)} \, 
\sum\limits_{\delta \in R_\ell(k)\cap \eta^{-1} P_\ell^{(r)}(k) \eta \bs P_\ell(k)}^{} \,
&
f_s(\epsilon_r \eta \delta u h) \Psi_\ell^{-1}(u) \, du.
\end{eqnarray*} 

\noindent{\bf Step 5.}   
We next describe some decompositions.  Since $R_\ell = G \cdot N_\ell,$ we have 
\begin{equation*} 
R_\ell \cap \eta^{-1} P_\ell^{(r)} \eta = \left( G \cap \eta^{-1} M_\ell^{(r)} \eta \right) \cdot \left( N_\ell \cap \eta^{-1} P_\ell^{(r)}\eta\right) 
\end{equation*} 
and 
\begin{equation*} 
(R_\ell \cap \eta^{-1} P_\ell^{(r)} \eta) \bs R_\ell = \left( (G \cap \eta^{-1} M_\ell^{(r)} \eta) \bs G \right) 
\cdot \left( (N_\ell \cap \eta^{-1} P_\ell^{(r)} \eta) \bs N_\ell \right). 
\end{equation*} 
We utilize this decomposition inside the $N_\ell$ integration. For fixed $r$ and $\eta$ the inner integration becomes 
\begin{equation*} 
\int\limits_{N_\ell(k) \bs N_\ell(\A)} \, 
\sum\limits_{\delta_2 \in (G(k) \cap \eta^{-1} M_\ell^{(r)}(k) \eta) \bs G(k)}^{} \,
\sum\limits_{\delta_1 \in (N_\ell(k) \cap \eta^{-1} P_\ell^{(r)}(k) \eta) \bs N_\ell(k)}^{} \,
f_s(\epsilon_r \eta \delta_1 \delta_2 u h) \Psi_\ell^{-1}(u) \, du.
\end{equation*} 
We next interchange the $u$ integral and the $\delta_2$ sum.  We can do this, as in \cite{grs}, by the absolute 
convergence for $\Re(s) \gg 0.$  Note that any modulus character will be $1$ on $\delta_2$ since it is 
rational, and since $\delta_2 \in G(k)$ it stabilizes the character $\Psi_\ell.$  After interchanging, we can collapse 
the $N_\ell$ integration and the $\delta_1$ summation to obtain 
\begin{equation*} 
\sum\limits_{\delta \in (G(k) \cap \eta^{-1} M_\ell^{(r)}(k) \eta) \bs G(k)}^{} \, 
\int\limits_{(N_\ell(k) \cap \eta^{-1} P_\ell^{(r)}(k) \eta) \bs N_\ell(\A)} \, 
f_s(\epsilon_r \eta u \delta h) \Psi_\ell^{-1}(u) \, du.
\end{equation*} 

\noindent{\bf Step 6.}   
We next fix $\delta$ and factor the $N_\ell$ integration as 
\begin{eqnarray*} 
\int\limits_{(N_\ell(k) \cap \eta^{-1} P_\ell^{(r)}(k) \eta) \bs N_\ell(\A)} \, 
f_s(\epsilon_r \eta u \delta h) \Psi_\ell^{-1}(u) \, du = & 
\\ 
\int\limits_{(N_\ell(\A) \cap \eta^{-1} P_\ell^{(r)}(\A) \eta) \bs N_\ell(\A)} \, 
\int\limits_{(N_\ell(k) \cap \eta^{-1} P_\ell (k) \eta) \bs (N_\ell(\A) \cap \eta^{-1} P_\ell^{(r)}(\A) \eta)} 
&
f_s(\epsilon_r \eta u' u \delta h) \Psi_\ell^{-1}(u' u) \, du' \, du. 
\end{eqnarray*} 
%

%%%%%%%%%%%
\begin{lem} \label{group-J}
If $r>0,$ then there exists a unipotent $k$-group $J \subset N_\ell \cap \eta^{-1} P_\ell^{(r)} \eta$ such that 
\begin{itemize}
\item $\Psi_\ell$ is non-trivial on $J(\A)$ and  
\item $\epsilon_r \eta J \eta^{-1} \epsilon_r^{-1} \subset U_m,$ the unipotent radical of $P_m.$ 
\end{itemize}
\end{lem}

\begin{proof} 
This exists in $\SO_{m'}$ by \cite[page 287 or the proof of Proposition 5.1]{grs}.  However, this is a unipotent element and the unipotent varieties are the same 
for $\SO_{m'}$ and $H.$  Hence, the group theoretic statements remain true in $H.$ 

The characters $\psi_{\ell, \alpha}$ of \cite[(3.10)]{grs} and $\Psi_\ell^{a}$ of \cite[Definition 16.1.9]{hs} are equal, so we get the first statement. 
\end{proof}

So if $r>0,$ since $\epsilon_r \eta J \eta^{-1} \epsilon_r^{-1} \subset U_m$ and the Eisenstein series is induced 
from $P_m,$ for $j \in J(\A)$, we see that 
\begin{eqnarray*} 
f_s\left(\epsilon_r \eta j u'' u \delta h\right) 
&=& 
f_s\left(\epsilon_r \eta j \eta^{-1} \epsilon_r^{-1} \epsilon_r \eta u'' u \delta h\right) \\
&=& f_s\left(\epsilon_r \eta u'' u \delta h\right)
\end{eqnarray*} 
and 
\begin{equation*} 
\int\limits_{J(k) \bs J(\A)} \Psi_\ell^{-1}(j) \, dj = 0. 
\end{equation*} 
Therefore, the unipotent $(N_\ell,\Psi_\ell)$ integration is zero except in the case of $r=0.$  Hence, we have the following.

%%%%%%%%%%%
\begin{prop}
\begin{eqnarray*}
E^{(N_\ell,\Psi_\ell)} (h,f_s) 
&=& 
\sum\limits_{\eta \in P_\ell^{(0)} \bs P_\ell(k) / R_\ell(k)} \, 
\sum\limits_{\delta \in (G(k) \cap \eta^{-1} M_\ell^{(0)} \eta) \bs G(k)} \\ 
&& 
\int\limits_{(N_\ell(k) \cap \eta^{-1} P_\ell^{(0)}(k) \eta) \bs N_\ell(\A)} \, 
f_s\left( \epsilon_0 \eta u \delta h \right) \Psi_\ell^{-1}(u) \, du.   
\end{eqnarray*}
\end{prop}

\noindent{\bf Step 7.}   
At this point we consider two cases:

{\bf Case I:} $m'$ is even, or $m'$ is odd and $a \not= t^2$ is not a square.  
In this case there is only one $\eta,$ namely $\eta = I_{m'}.$

{\bf Case II:} $m'$ is odd and $a=t^2$ is a square.  In this case there are 
three choices for $\eta,$ namely, $\eta = I_{m'},$ $\eta_+,$ and $\eta_-,$ where 
\begin{equation*} 
\eta_{\pm} = 
\begin{pmatrix} 
I_\ell && \\ & \gamma_{\pm} & \\ && I_\ell
\end{pmatrix}
\end{equation*} 
and 
\begin{equation*} 
\gamma_{\pm} = 
\begin{pmatrix} 
I_{m'-\ell-1} &&&& \\ 
& 1 &&& \\ 
& \pm t & 1 && \\ 
& - \frac{t^2}{2} & \mp t & 1 & \\ 
&&&& I_{m' - \ell - 1}
\end{pmatrix}. 
\end{equation*} 
%

%%%%%%%%%%%
\begin{lem}
The above are still the coset representatives for $P_\ell^{(0)}(k) \bs P_\ell(k) / R_\ell(k).$ 
\end{lem}

\begin{proof} 
This lemma follows from the agreement of the unipotent varieties for $\GSpin$ and $\SO$ along 
with the fact that the differences between $\GSpin$ and $\SO$ will lie in $P_\ell(k)$ and $R_\ell(k).$ 
\end{proof}

We can now eliminate Case II and $\eta = \eta_{\pm}.$

%%%%%%%%%%%
\begin{prop}
In the case of $m' = 2m+1$ and $\eta = \eta_{\pm},$ the contributions of these terms 
to the unipotent period of the Eisenstein series vanish. 
\end{prop}

\begin{proof}
We let $\eta = \eta_{\pm} \not= I_{m'}.$  We have 
$N_\ell(k) \cap \eta^{-1} P_\ell^{(0)}(k) \eta = \eta^{-1} \left( N_\ell(k) \cap P_\ell^{(0)} \right) \eta.$ 
Therefore, 
\begin{eqnarray*} 
\int\limits_{(N_\ell(k) \cap \eta^{-1} P_\ell^{(0)}(k) \eta) \bs N_\ell(\A)} \, 
f_s\left(\epsilon_0 \eta u \delta h\right) \Psi_\ell^a(u)^{-1} \, du = & 
\\ 
\int\limits_{(N_\ell(\A) \cap \eta^{-1} P_\ell^{(0)}(\A) \eta) \bs N_\ell(\A)} \, 
\int\limits_{(N_\ell(k) \cap P_\ell^{(0)}(k)) \bs (N_\ell(\A) \cap P_\ell^{(0)}(\A))} 
&
f_s(\epsilon_0 v \eta u \delta h) \Psi_\ell^0(\eta^{-1} v \eta u)^{-1} \, dv \, du. 
\end{eqnarray*} 

In terms of matrices, as in \cite{hs}, using the agreement of unipotent varieties, 
\begin{equation*} 
N_\ell \cap P_\ell^{(0)} = 
\left\{ 
v = \begin{pmatrix} 
z & 0 & 0 & y & 0 \\ 
& I_{m-\ell} & 0 & 0 & y' \\ 
&& 1 & 0 & 0 \\ 
&&& I_{m-\ell} & 0 \\ 
&&&&z^*
\end{pmatrix} : z \in Z_\ell
\right\}. 
\end{equation*} 
Note that 
\begin{equation*} 
\Psi_\ell^a(\eta^{-1}v\eta) = \psi(z_{1,2} + \cdots + z_{\ell-1,\ell}) = \psi_\ell(z) 
\end{equation*} 
and this is independent of $y.$  If we conjugate this past $\epsilon_0,$ we obtain 
\begin{equation*} 
\epsilon_0 \left( N_\ell \cap P_\ell^{(0)} \right) = 
\left\{ 
\epsilon_0 v \epsilon_0^{-1} = \hat{z}'_\ell = 
\begin{pmatrix} 
I_{m-\ell} & x \\ 
& \zeta
\end{pmatrix}^\wedge
\right\}  = \hat{Z}'_\ell,
\end{equation*} 
where for $g \in \GL_{t},$ we set $\hat{g}$ to be the lift of $g$ into the Levi $\GL_{t} \times \GL_{m'-2t}$ 
of $P_t.$  Note that even though we have not specified $\zeta$, as in \cite[p. 287]{grs} we know that 
$\zeta \in Z_\ell$.  In these coordinates 
\begin{equation*} 
\Psi_\ell^a(\eta^{-1} v \eta) = \psi(\zeta_{1,2} + \cdots + \zeta_{\ell-1,\ell}) 
= \psi(z'_{m-\ell+1,m-\ell+2} + \cdots + z'_{m-1,m}). 
\end{equation*} 

If we denote this character of $Z'_\ell$ by $\psi_{Z'_\ell}^0$ and 
\begin{equation*} 
f_s^{(Z'_\ell,\psi_{Z'_\ell}^0)} (h) = 
\int\limits_{Z'_\ell(k) \bs Z'_\ell(\A)} f_s(\hat{z}'_\ell h) \psi_{Z'_\ell}^0 (z'_\ell)^{-1} \, dz'_\ell, 
\end{equation*} 
then we have 
\begin{equation*} 
\int\limits_{N_\ell(k) \cap P_\ell^{(0)}(k) \bs N_\ell(\A) \cap P_\ell^{(0)}(\A)} 
f_s\left( \epsilon_0 v \eta u \delta h \right) \Psi_\ell^0 \left( \eta^{-1} v \eta \right)^{-1} \, dv 
= f_s^{ (Z'_\ell,\psi^0_{Z'_\ell})} 
\left( \epsilon_0 \eta u \delta h \right) 
\end{equation*} 
and finally 
\begin{equation*} 
\int\limits_{\left. \left(N_\ell(k) \cap \eta^{-1} P_\ell^{(0)}(k) \eta \right) \middle\backslash N_\ell(\A) \right.}
f_s\left( \epsilon_0 \eta u \delta h \right) \Psi_\ell^a (u)^{-1} \, du 
=  
\int\limits_{ \left. \left(N_\ell(\A) \cap \eta^{-1} P_\ell^{(0)}(\A) \eta\right) \middle\backslash N_\ell(\A) \right.} 
f_s^{ (Z'_\ell,\psi^0_{Z'_\ell})}
\left( \epsilon_0 \eta u \delta h \right) \Psi_\ell^a (u)^{-1} \, du. 
\end{equation*} 
However, in the unipotent period $f_s^{(Z'_\ell,\psi^0_{Z'_\ell})}(\epsilon_0 \eta u h),$ 
as an inner integral we have the constant term of $f_s$ along the unipotent radical of 
the standard parabolic subgroup of $\GL_{m}$ which corresponds to the 
partition $(m-\ell,\ell)$ of $m.$  Since we are inducing from a cuspidal representation $\tau$ of 
$\GL_{m}(\A)$ and $1 \le \ell \le m,$ we have that this constant term is $\equiv 0.$  Hence, 
\begin{equation*} 
\int\limits_{(N_\ell(k) \cap \eta^{-1} P_\ell^{(0)}(k) \eta) \bs N_\ell(\A)}
f_s\left( \epsilon_0 \eta u \delta h \right) \Psi_\ell^a (u)^{-1} \, du 
\equiv  0
\end{equation*} 
as desired. 
\end{proof}
%%%%%%%%%%%

Now, in either case I or Case II we are reduced to $\eta = I_{m'}.$  As a 
consequence of the previous proposition, we have the following corollary.

%%%%%%%%%%%
\begin{cor} 
\begin{equation*} 
E^{(N_\ell,\Psi_\ell)} \left( h,f_s \right) = 
\sum\limits_{\delta \in (G(k) \cap M_\ell^{(0)}(k)) \bs G(k)} \, 
\int\limits_{(N_\ell(k) \cap P_\ell^{(0)}(k)) \bs N_\ell(\A)} 
f_s\left( \epsilon_0 u \delta h \right) \Psi_\ell(u)^{-1} \, du.
\end{equation*} 
\end{cor}

\noindent{\bf Step 8.}   
If we take this last expression for the  
Fourier coefficient of the Eisenstein series and insert it into 
our global integral we arrive at the expression
\begin{eqnarray*} 
\L(\varphi,f_s) 
&&= 
\int\limits_{Z(\A) G(k) \bs G(\A)} 
\varphi(g) E^{(N_\ell,\Psi_\ell)}(g,f_s) \, dg  
\\ 
=&& 
\int\limits_{Z(\A) G(k) \bs G(\A)} 
\varphi(g) 
\left( 
\sum\limits_{\delta \in (G(k) \cap M_\ell^{(0)}(k)) \bs G(k)} \, 
\int\limits_{ (N_\ell(k) \cap P_\ell^{(0)}(k)) \bs N_\ell(\A) } 
f_s\left( \epsilon_0 u \delta g \right) \Psi_\ell(u)^{-1} \, du 
\right) \, dg.  
\end{eqnarray*} 
Since $\delta \in G(k)$ and $\varphi$ is automorphic, we can bring $\varphi$ into the sum over 
$\delta$ and replace $\varphi(g)$ by $\varphi(\delta g).$  Then we can absorb 
the $\delta$ sum into the integral over $G(\A).$

%%%%%%%%%%%
\begin{prop} 
\begin{equation*} 
\L(\varphi,f_s) = 
\int\limits_{Z(\A) \left( G(k) \cap M_\ell^{(0)}(k) \right) \bs G(\A)} 
\varphi(g) 
\int\limits_{ (N_\ell(k) \cap P_\ell^{(0)}(k)) \bs N_\ell(\A) } 
f_s\left( \epsilon_0 u g \right) \Psi_\ell(u)^{-1} \, du \, dg.  
\end{equation*} 
\end{prop}

As noted above, if we conjugate the unipotent past $\epsilon_0$ we obtain 
\begin{equation*} 
\epsilon_0 \left( N_\ell \cap P_\ell^{(0)} \right) \epsilon_0^{-1}
= 
\left\{ 
\epsilon_0 u \epsilon_0^{-1} = \hat{z}'_\ell = 
\begin{pmatrix} 
I_{m-\ell} & x \\
& \zeta
\end{pmatrix}^\wedge
\right\}  
= \hat{Z}'_\ell, 
\end{equation*} 
where for $g \in \GL_{t},$ we denote by $\hat{g}$ the lift of $g$ into the Levi 
$\GL_{t} \times \GSpin_{m' - 2t}$ of $P_t.$  If we transfer the character 
$\Psi_\ell = \Psi_\ell^{a}$ to $Z'_\ell,$ we find 
\begin{equation*} 
\Psi_\ell^a \left( \epsilon_0^{-1} \hat{z}'_\ell \epsilon \right) 
= \psi\left( \zeta_{1,2} + \cdots + \zeta_{\ell-1,\ell} + \frac{a}{2} x_{m-\ell,1} \right)
\end{equation*} 
and so on the group $Z'_\ell$ we set 
\begin{equation*} 
\Psi_{Z'_\ell} \left( z'_\ell  \right) 
= 
\Psi_{Z'_\ell}^a \left( z'_\ell  \right) 
= 
\psi\left( \frac{a}{2} z'_{m-\ell,m-\ell+1} + z'_{m-\ell+1,m-\ell+2} + \cdots + z'_{m-1,m} \right). 
\end{equation*} 
Then if we set 
\begin{eqnarray*} 
f_s^{ (Z'_\ell,\Psi_{Z'_\ell}) }(h) 
&=& 
\int\limits_{Z'_\ell(k) \bs Z'_\ell(\A)} 
f_z(\hat{z}'_\ell) \Psi_{Z'_\ell}(z'_\ell)^{-1} \, dz'_\ell 
\\
&=& 
\int\limits_{Z'_\ell(k) \bs Z'_\ell(\A)} 
f_z(\hat{z}'_\ell) \Psi_\ell (\epsilon_0^{-1} z'_\ell \epsilon_0)^{-1} \, dz'_\ell,  
\end{eqnarray*} 
then part of the integral that appears in the global integral is the unipotent period 
\begin{equation*} 
\int\limits_{ (N_\ell(k) \cap P_\ell^{(0)}(k)) \bs (N_\ell(\A) \cap P_\ell^{(0)}(\A)) } 
f_s\left( \epsilon_0 u' u g \right) \Psi_\ell(u')^{-1} \, du'  
= 
f_s^{(Z'_\ell,\Psi_{Z'_\ell})}\left( \epsilon_0 u g \right) 
\end{equation*} 
and the integral representation itself becomes 
\begin{equation*} 
\L(\varphi,f_s) 
= 
\int\limits_{Z(\A) \left( G(k) \cap M_\ell^{(0)}(k) \right) \bs G(\A)} 
\varphi(g) 
\int\limits_{ (N_\ell(\A) \cap P_\ell^{(0)}(\A)) \bs N_\ell(\A) } 
f_s^{(Z'_\ell,\Psi_{Z'_\ell})}\left( \epsilon_0 u g \right) \Psi_\ell(u)^{-1} \, du \, dg.  
\end{equation*} 

Next we turn to the Whittaker coefficient. Recall that in Step 1 we commented on convergence issues 
and steps 2 through 8 only deal with the Eisenstein series and have nothing to do with the $dg$ integration. 
The purpose in these steps was to unfold and simplify the Eisenstein series and show that  certain 
terms vanish. In step 9 below we show that integral factors through a Whittaker-Fourier coefficient of $\varphi$. 

\noindent{\bf Step 9.}   
We now have to analyze $G(k) \cap M_\ell^{(0)}(k).$ Recall that 
\begin{equation*} 
G \cong 
\begin{cases} 
\GSpin_{2(m-\ell)-1} & \mbox{ if } m' = 2m, \\
\GSpin^a_{2(m-\ell)} & \mbox{ if } m' = 2m+1. 
\end{cases}
\end{equation*} 

As in \cite{grs}, we have $G \cap M_\ell^{(0)} \subset P_{G,m-\ell},$ the ``Siegel'' parabolic subgroup of $G$.   
(We placed ``Siegel'' in quotes since in the case of $a \neq t^2.$ we did not define it as such.) 
If we consider the cases, then we have the following:  

Case (i): $m' = 2m$ even. In this case $G \cap M^{(0)} \subset P_{G,m-\ell-1},$ the Siegel parabolic 
of $G.$  The Levi of the Siegel is then $\GL_{m-\ell-1} \times \GL_{1}.$  We may have 
$G \cap M^{(0)} = \GL_{m-\ell-1} \ltimes U_{G,m-\ell-1}$ or $G \cap M^{(0)} = P_{G,m-\ell-1}.$ 

Case (ii): $m' = 2m+1$ and $a \neq t^2.$  In this case we have that $G \cap M^{(0)} \subset P_{G,m-\ell-1}.$  
This parabolic subgroup of $G$ has Levi subgroup given by $\GL_{m-\ell-1} \times \GSpin^a_{2}.$  Again,  
$G \cap M^{(0)} = \GL_{m-\ell-1} \ltimes U_{G,m-\ell-1}.$ 

Case (iii): $m' = 2m+1$ and $a = t^2$ a square.  We still 
have $G \cap M^{(0)} \subset P_{G,m-\ell}$. (This case does not appear in \cite{grs} since it was not needed there.)

On the other hand, the actual form of the Levi part of $G \cap M^{(0)}$ does not play a role.  What one needs is the following.  
Write $G \cap M^{(0)} = P' = M' U'.$  Then $U'$ is the same as in \cite{grs} since it is unipotent.

Define 
\begin{equation*} 
C = C_{G,m-\ell}= \left\{ u' \in U' : u' e_m = e_m \right\}. 
\end{equation*} 
%

%%%%%%%%%%%
\begin{lem}
We have 
\begin{itemize}
\item $ \epsilon_0 C_{G,m-\ell} \epsilon_0^{-1} \subset U_{H,m},$ the unipotent radical of the Siegel parabolic $H.$ 
\item $C_{G,m-\ell}$ normalizes $N_\ell$ and commutes with $N_\ell \cap P_\ell^{(0)}.$ 
\item $C_{G,m-\ell}$ preserves $\Psi_\ell.$ 
\end{itemize} 
\end{lem}

\begin{proof}
Note that $C_{G,m-\ell}$ has these properties in the $\SO$ case, so by the equality of unipotent varieties 
and lifting of Weyl elements remains true in the $\GSpin$ case. 
\end{proof}

Therefore, this unipotent subgroup leaves $f_s^{(Z'_\ell,\Psi_{Z'_\ell})}\left( \epsilon_0 u g \right) \Psi_\ell(u)^{-1}$ invariant.  
So we can factor it through and take this unipotent period of $\varphi.$  Note that since this is not the unipotent radical 
of a parabolic subgroup, there is no reason for this period of $\varphi$ to vanish.  Then 
\begin{equation*} 
\L(\varphi,f_s) 
=
\int\limits_{P'(k)C(\A)Z(\A) \bs G(\A)} \varphi^{(C,{\mathbf 1})}(g)
\int\limits_{( N_\ell(\A) \cap P_\ell^{(0)}(\A) ) \bs N_\ell(\A)} 
f_s^{(Z'_\ell,\Psi_{Z'_\ell})}\left( \epsilon_0 u g \right) \Psi_\ell(u)^{-1} \, du dg.
\end{equation*} 

\noindent{\bf Step 10.}   
We next need the following Lemma.

%%%%%%%%%%%
\begin{lem} 
$C \bs P' \cong P^1_{m-\ell} \subset \GL_{m-\ell},$ the mirabolic subgroup of $\GL_{m-\ell}.$ 
\end{lem}

\begin{proof}
We can see this in a similar way as in \cite[p. 288--289]{grs}.  There are obvious similitude analogs of \cite[(10.16) and (10.17)]{grs} 
which give the elements of $P^1_{m-\ell}$. The similitude analogs would only have the similitude character in the term $d^*$ in the 
notation of \cite[(10.16) and (10.17)]{grs} (so that the element is indeed in the $\GSpin$ group) and otherwise the same formulas apply. 
We could then give the isomorphism by sending the cosets of $C$ in $C \bs P'$ 
to their corresponding elements in $P^1_{m-\ell}$ as in \cite[p. 289]{grs}. 
\end{proof}

Since $\varphi^{(C,{\mathbf 1})}(g)$ is left invariant under both $C(\A),$ by taking period, and $P'(k),$ 
since $\varphi$ is cuspidal, we can use the process of Piatetski-Shapiro and Shalika to Fourier expand 
along $C \bs P' \cong P^1_{m-\ell}.$  Note that 
$P^1_{m-\ell} \cong \GL_{m-\ell-1} \ltimes k^{m-\ell-1}$ and the $\GL_{m-\ell-1}$ lies in the 
Levi of the parabolic subgroup of $G$ with Levi $\GL_{m-\ell-1} \times \GL_{1}$ if $m' = 2m$ is even and 
$\GL_{m-\ell-1} \times \GSpin^a_{2}$ if $m' = 2m+1$ is odd.  The $F^{m-\ell-1}$ represents a unipotent subgroup.  
So this construction should lift from $\SO_{2(m-\ell)-1}$ (resp. $\SO^a_{2(m-\ell)}$) to $G.$  
Then, according to \cite{grs} we have 
\begin{equation*} 
\varphi^{(C,{\mathbf 1})}(g) 
=
\sum\limits_{d \in Z_{m-\ell-1}(k) \bs \GL_{m-\ell-1}(k)} 
W^\psi_\varphi\left( \begin{pmatrix} I_\ell && \\ & d & \\ && 1\end{pmatrix}^\wedge g \right)
= 
\sum\limits_{q \in N_G(k) \bs P'(k) } 
W^\psi_\varphi(qg),
\end{equation*} 
where 
\begin{equation*} 
W^\psi_\varphi(g) 
= 
\int\limits_{N_G(k) \bs N_G(\A)} 
\varphi(vg) \psi(v)^{-1} \, dv 
\end{equation*} 
and 
$\psi = \psi_{N_G}$ is the standard Whittaker character on $N_G.$

If we replace $W^\psi_\varphi(g)$ with its Fourier expansion and unfold the sum over 
$N_G(k) \bs P'(k)$ against the integral over $P'(k) C(\A) \bs G(\A),$ we obtain 
\begin{equation*} 
\L(\varphi,f_s) 
= 
\int\limits_{N_G(k) C(\A) Z(\A) \bs G(\A)} 
W^\psi_\varphi(g) 
\int\limits_{N_\ell(\A) \cap P_\ell^{(0)}(\A) \bs N_\ell(\A) } 
f_s^{(Z'_\ell,\Psi_{Z'_\ell})}\left( \epsilon_0 u g \right) \Psi_\ell(u)^{-1} \, du dg.  
\end{equation*} 

\noindent{\bf Step 11.}   
We next factor the $dg$ integration through $N_G(k) \bs N_G(\A).$  
The Whittaker function will return $\psi_{N_G}(u').$  So we first obtain 
\begin{eqnarray*} 
\L(\varphi,f_s) 
= 
\int\limits_{N_G(\A) Z(\A) \bs G(\A)} 
W^\psi_\varphi(g) 
&&
\int\limits_{N_G(k) C(\A) \bs N_G(\A)} 
\psi_{N_G}(u') 
\\
&&
\int\limits_{N_\ell(\A) \cap P_\ell^{(0)}(\A) \bs N_\ell(\A) } 
f_s^{(Z'_\ell,\Psi_{Z'_\ell})}\left( \epsilon_0 u u' g \right) \Psi_\ell(u)^{-1} \, du du' dg.  
\end{eqnarray*} 

Now $u \in N_\ell(\A)$ and $u' \in N_G(\A).$  Since $G = L_{\ell,a}$ lies in the Levi 
subgroup $L_\ell,$ we know that $u'$ normalizes $u$ and does not change $\Psi_\ell.$  
Therefore, we can interchange the (compact) $u$ and $u'$ integrations and write the 
argument of $f_s$ as 
\begin{equation*} 
f_s^{(Z'_\ell,\Psi_{Z'_\ell})}\left( \epsilon_0 u' u g \right). 
\end{equation*} 
Following \cite{grs}, we decompose $\epsilon_0 u' \epsilon_0^{-1} = z' u''$ 
with $z' \in Z_m,$ the maximal unipotent of $\GL_{m} \in H$ and $u'' \in U_m,$  
the unipotent radical of the Siegel parabolic $P_m \subset H.$

%%%%%%%%%%%
\begin{lem}
If we write $\epsilon_0 u' \epsilon_0^{-1} = z' u''$ as above, then 
\begin{itemize}
\item $Z_m = Z'_\ell \cdot \left\{ z' :  \epsilon_0^{-1} z' \epsilon_0 \in N_G  \right\}$ 
\item $\psi_{N_G} \left( \epsilon_0^{-1} z' \epsilon_0 \right) = 
\Psi_m^a \vert_{\{ z' : \epsilon_0^{-1} z' \epsilon_0 \in N_G \}}$ 
for $\Psi_m^a$ some non-degenerate character of $Z_m.$ 
\end{itemize}
\end{lem}

\begin{proof} 
This is true in $\SO$ and, by the agreement of unipotent varieties and lifts of Weyl elements $(\epsilon_0),$ 
it also holds in $\GSpin.$
\end{proof} 

Therefore 
\begin{equation*} 
\int\limits_{N_G(k) C(\A) \bs N_G(\A)} 
\psi_{N_G}(u') f_s^{(Z'_\ell,\Psi_{Z'_\ell})}\left( \epsilon_0 u' u g \right) \, du' 
= 
f_s^{(Z_m,\Psi_m^a)}\left( \epsilon_0 u g \right), 
\end{equation*} 
where $\Psi_m^a$ is a non-degenerate character of $Z_m$ and 
\begin{equation*} 
\L(\varphi,f_s) 
= 
\int\limits_{N_G(\A) Z(\A) \bs G(\A)} 
W^\psi_\varphi(g) 
\int\limits_{N_\ell(\A) \cap P_\ell^{(0)}(\A) \bs N_\ell(\A) } 
f_s^{(Z_m,\Psi^a_m)}\left( \epsilon_0 u g \right) \Psi_\ell(u)^{-1} \, du dg.  
\end{equation*} 

When $G = \GSpin_{2m}$ it may be occasionally useful, particularly in the local context, to consider an arbitrary 
Whittaker character $\psi$ of $N_G$ for the Whittaker model of $\pi$.  One can do this by making $\Psi^a_m$ or 
$\beta$ explicit (each determining the other).

\noindent{\bf Step 12.}   
Choose an element $d_a \in T_m \subset \GL_{m}$ which conjugates $\Psi_m^a$ to the standard 
character $\psi_m$ of $Z_m,$ i.e., 
\begin{equation*} 
\Psi^a_m\left( d_a x d_a^{-1} \right) = \psi_m(z), \quad z \in Z_m. 
\end{equation*} 
This all takes place in the $\GL_{m}$ Levi of the Siegel parabolic in $\SO_{m'}$ and so the same 
is true in $\GSpin_{m'} = H.$   Then 
\begin{equation*} 
f_s^{(Z_m,\Psi^a_m)}(h) = 
f_s^{(Z_m,\psi_m)}\left( \hat{d}_a h \right).  
\end{equation*} 
So if we let  
\begin{equation*} 
\beta = \beta_{\ell,a} = \hat{d}_a \epsilon_0, 
\end{equation*} 
then 
\begin{equation*} 
\L(\varphi,f_s) 
= 
\int\limits_{N_G(\A) Z(\A) \bs G(\A)} 
W^\psi_\varphi(g) 
\int\limits_{N_\ell(\A) \cap P_\ell^{(0)}(\A) \bs N_\ell(\A) } 
f_s^{(Z_m,\psi_m)}\left( \beta u g \right) \Psi_\ell(u)^{-1} \, du dg.  
\end{equation*} 
We finally note that 
\begin{equation*} 
N_\ell \cap P^{(0)} = N_\ell \cap P_\ell \cap \epsilon_0^{-1} P_m \epsilon_0 = N_\ell \cap \beta^{-1} P_m \beta 
\end{equation*} 
since $\beta = \hat{d}_a \epsilon_0$ and $\hat{d}_a \in T_m \subset \GL_{m} \subset P_m.$ Then 
\begin{equation*} 
\L(\varphi,f_s) 
= 
\int\limits_{N_G(\A) Z(\A) \bs G(\A)} 
W^\psi_\varphi(g) 
\int\limits_{N_\ell(\A) \cap \beta^{-1} P_m(\A) \beta \bs N_\ell(\A) } 
f_s^{(Z_m,\psi_m)}\left( \beta u g \right) \Psi_\ell(u)^{-1} \, du dg.  
\end{equation*} 

This finishes the proof of Theorem \ref{thm-BasicIdentityI}
\end{proof}

%%%%%%%%%%%%%% NEW SECTION %%%%%%%%%%%%%%%%%%%% 
%%%%%%%%%%%%%%%%%%%%%%%%%%%%%%%%%%%%%%%%%%%%%%% 
 \section{Global Integrals II (Case A)} \label{sec-integralA} 

In this section we translate \cite[\S 10.4]{grs} into the $\GSpin$ context. This corresponds to 
Method A  in the original work of Gelbart and Piatetski-Shapiro \cite[Part B]{gpsr}, which dealt with the special 
orthogonal and general linear groups with equal (or nearly equal) ranks.  As such, we refer 
to the integrals of this section as Case A. (See the table in Section \ref{sec-euler-exp} for a summary of various  
cases.) We again note that we follow 
\cite{hs} in the labeling of various parabolic subgroups (cf. Remark \ref{par-grs-acs}). 

We still take $H = \GSpin_{m'}$ to be the larger group.  Now, $(\pi, V_\pi)$ will be a cuspidal 
representation of $H(\A)$ and $H$ will be either split or quasi-split.  The group $G = \GSpin_{n'}$ will be 
the smaller group, split, and having the opposite parity to $H.$

We have that $H$ is the $\GSpin$ cover of an orthogonal group $\SO_{m'}(V)$ with $V$ a quadratic 
space.  Let $\tilde{m}$ be the Witt index of $V.$  Hence, 
\begin{equation*}
\tilde{m} = 
\begin{cases}
m & \mbox{ if } m' = 2m+1, \\
m & \mbox{ if } m' = 2m \mbox{ with } H \mbox{ split, } \\
m-1 & \mbox{ if } m' = 2m \mbox{ with } H \mbox{ quasi-split. } 
\end{cases}
\end{equation*}
Let $ 0 \le \ell$ be such that $ \ell < \tilde{m}$ if $m'$ is even and $\ell < \tilde{m} -1$ if $m'$ is odd.  
Let $Q_\ell = L_\ell \ltimes N_\ell$ be the parabolic subgroup of $H$ with Levi of the form 
$(\GL_{1})^\ell \times \GSpin_{m' - 2\ell}.$  Let $N_\ell$ be its unipotent radical. Let $\Psi_\ell$ be 
a character of $N_\ell$ so that $\operatorname{Stab}_{L_\ell}(\Psi_\ell)$ is a split $\GSpin_{m'-2\ell-1}.$  
If $H$ is split, we can take the character $\Psi_\ell$ from the previous section.

When $H$ is quasi-split, we find a character $\Psi_\ell$ of $N_\ell$ such that 
$\operatorname{Stab}_{L_\ell}(\Psi_\ell)$ is split and the $\GSpin$ cover of the $\SO(W_{m,\ell}\cap w_0^\perp)$ of 
\cite{grs}.  This is done in the same as the $\Psi_\ell$ from the previous section for $a=1$.   
In \cite{hs} the authors do not consider this case since for descent they can always take the larger group $H$ to be split.  
Note that if $H$ is quasi-split, so an even $\GSpin$, then $\SO(W_{m,\ell}\cap w_0^\perp)$ is an odd $\SO$ 
and so automatically split. Hence, it does not matter which anisotropic vector we take.

Let $G = \operatorname{Stab}_{L_\ell}(\Psi_\ell).$  Then $G = \GSpin_{m' - 2\ell - 1}.$ 
We let $n' = m' - 2\ell -1.$  If $m' = 2m+1$ is odd, then $n' = 2n$ with $n = m-\ell$ and 
if $m' = 2m$ is even, then $n' = 2m - 2\ell - 1 = 2n+1$ with $n = m-\ell-1.$ In either case $G$ is split. 

Again, when $\ell = 0$ there is no $\Psi_\ell$ and we just ``restrict'' and  the statements below 
hold in the $\ell=0$ case as well.

Let $P_G = M_G \ltimes U_G$ be the Siegel parabolic subgroup of $G.$  Then $M_G = \GL_{n} \times \GL_{1}$ 
with the $\GL_{1}$ factor being the connected component of the center of $G$.  Let $\tau$ be a cuspidal automorphic 
representation of $\GL_{n}(\A)$ and let 
\begin{equation} \label{fs-II} 
\rho = \rho_{\tau, \omega_\pi^{-1}, s} = \operatorname{Ind}_{P_G}^G \left( \tau \vert\det\vert^{s} \otimes \omega_{\pi}^{-1} \right)
\end{equation} 
so that $\rho_s$ transforms by $\omega_\pi^{-1}$ under the connected component of the center of $G.$  
(Again, we will replace $s$ by $s-1/2$ below.) 
Let $f_s$ be a $K$-finite holomorphic section of $\rho$ and form the Eisenstein series 
\begin{equation*} 
E(g,f_s) = \sum\limits_{\delta \in P_G(k) \bs G(k)} f_s(\delta g)
\end{equation*} 
which is absolutely convergent, uniformly on compact subsets, for $\Re(s) \gg 0.$

For $\varphi \in V_\pi,$ a cusp form on $H(\A),$ define 
\begin{equation*} 
\L(\varphi, f_s) = \int\limits_{G(k) Z(\A) \bs G(\A)} \varphi^{(N_\ell,\Psi_\ell)}(g) E(g,f_s) \, dg,  
\end{equation*} 
where $Z = Z(G)^0$ is the identity component of the center of $G$ as in the earlier case. 

Similarly to Lemma \ref{Siegel-estimates}, the integral $\L(\varphi,f_s)$ converges absolutely 
and uniformly in vertical strips in $\C$ away from poles of the Eisenstein series and hence 
it defines a meromorphic function on $\C.$

%%%%%%%%%%%
\begin{thm}\label{thm-BasicIdentityII}
Let the notaion be as above.  
\begin{itemize}
\item[(i)] If $\L(\varphi, f_s)$ is not identically zero, then $\pi$ is globally generic with respect to a certain 
``standard'' Whittaker character $\psi.$ (For a precise description of the character $\psi$ of $N_G$ is all cases 
we refer to \cite{soudry-mem} and \cite{kaplan-thesis}.) 
\item[(ii)] For $\Re s \gg 0$ we have an identity 
\begin{equation*} 
\L(\varphi, f_s) = \int\limits_{N_G(\A) Z(\A) \bs G(\A)} \int\limits_{\X(\A)} 
W_\varphi^\psi (\lambda \delta_\ell g) f_s^{(Z_n,\psi_n)}(g) \, d\lambda \, dg, 
\end{equation*} 
where $W_\varphi^\psi$ is the appropriate Whittaker function of $\varphi,$ $\X$ is isomorphic to a unipotent subgroup 
of $\GL_{m}$ if $m'=2m+1,$ resp. $\GL_{m-1}$ if $m'=2m,$ and is of the form 
\begin{equation} \label{X}
\X = \left\{ \begin{pmatrix} I_n & \\ \lambda & I_\ell \end{pmatrix}^\wedge \right\},
\end{equation} 
the element 
\begin{equation} \label{delta-ell} 
\delta_\ell = \begin{pmatrix} 0 & I_n \\ I_\ell & 0 \end{pmatrix}^\wedge
\end{equation} 
is a Weyl group element of $\GL_{m},$ resp. $\GL_{m-1},$ as above, and 
$f_s^{(Z_n,\Psi_n)}$ is the same unipotent period that appears in 
Theorem \ref{thm-BasicIdentityI} 
(so $Z_n$ is the maximal unipotent subgroup of $\GL_{n}.$) 
\end{itemize}
\end{thm}

%%%%%%%%%%%
\begin{proof}
For $\Re(s) \gg 0$ the Eisenstein series $E(g,f_s)$ is an absolutely convergent series.  We replace it by its definition: 
\begin{equation*} 
\L(\varphi,f_s) = \int\limits_{G(k) Z(\A) \bs G(\A)}  \varphi^{(N_\ell,\Psi_\ell)}(g)  
\sum\limits_{\delta \in P_G(k) \bs G(k)} f_s(\delta g) \, dg.
\end{equation*} 
Since $\varphi$ is automorphic and $G = \operatorname{Stab}(\Psi_\ell),$ we can move the summation outside 
$\varphi^{(N_\ell,\Psi_\ell)}$ 
\begin{equation*} 
\L(\varphi,f_s) = \int\limits_{G(k) Z(\A) \bs G(\A)} 
\sum\limits_{\delta \in P_G(k) \bs G(k)} 
\varphi^{(N_\ell,\Psi_\ell)}(\delta g)  
f_s(\delta g) \, dg 
\end{equation*} 
and collapse the sum against the integral to obtain 
\begin{equation*} 
\L(\varphi,f_s) = \int\limits_{P_G(k) Z(\A) \bs G(\A)} 
\left(\varphi^{(N_\ell,\Psi_\ell)}\right)^{(U_G,{\mathbf 1})} (g)  
f_s(g) \, dg.  
\end{equation*} 

We now come to a crucial result.

%%%%%%%%%%%
\begin{prop} \label{sum-int-exp}
We have 
\begin{equation*} 
\left(\varphi^{(N_\ell,\Psi_\ell)}\right)^{(U_G,{\mathbf 1})} (g) = 
\sum\limits_{\gamma \in Z_n(k) \bs \GL_{n}(k)} 
\int\limits_{\X(\A)} 
W_\varphi^\psi\left(\hat{\gamma} \lambda \delta_\ell g\right) \, d\lambda. 
\end{equation*} 
\end{prop}

%%%%%%%%%%%
\begin{proof}
We will prove this Proposition in Section \ref{root-exchange-sec} after we review the process of ``root exchange''.  
The analogous result for the special orthogonal groups is \cite[Theorem 7.3]{grs}.  
\end{proof}

With this, we have 
\begin{equation*} 
\L(\varphi,f_s) = \int\limits_{M_G(k) U_G(\A) Z(\A) \bs G(\A)} 
\sum\limits_{\gamma \in Z_n(k) \bs \GL_{n}(k)} 
\int\limits_{\X(\A)} 
W_\varphi^\psi\left(\hat{\gamma} \lambda \delta_\ell g\right) f_s(g) \, d\lambda \, dg,  
\end{equation*} 
where $U_G$ is the unipotent radical of the Siegel parabolic $P_G$ as in above \eqref{fs-II}. (In particular, we have 
$Z_n U_G = N_G$, the maximal unipotent subgroup of $G$.)

Now, for $\gamma \in Z_n(k) \bs \GL_{n}(k),$ we have that $\hat{\gamma}$ normalizes the group $\X.$ 
This is true whether we are in the $\SO$ context or the $\GSpin$ context because it is taking place in the $\GL_{m}$ 
or $\GL_{m-1}$ Levi subgroups. 

We claim that $\delta_\ell^{-1} \hat{\gamma} \delta_\ell$ is a general element of $M_G(k) / Z(k) = \GL_{n}(k).$  
Again this is true because it takes palce in the $\GL_{m},$ resp. $\GL_{m-1},$ Levi subgroup.

So we now move $\hat{\gamma}$ to the right and then collapse the sum over $Z_n(k) \bs \GL_{n}(k)$ with the 
integral to obtain 
\begin{equation*} 
\L(\varphi,f_s) = 
\int\limits_{Z_n(k) U_G(\A) Z(\A) \bs G(\A)} 
\int\limits_{\X(\A)} 
W_\varphi^\psi\left(\lambda \delta_\ell g\right) f_s(g) \, d\lambda \, dg. 
\end{equation*} 

If we now integrate over $Z_n(k) \bs Z_n(\A),$ then $W_\varphi^\psi\left(\lambda \delta_\ell g\right)$ 
will produce a character $\psi_n$ of $Z_n.$  We then integrate this unipotent period for $f_s$ and obtain 
\begin{equation*} 
\L(\varphi,f_s) = 
\int\limits_{Z_n(k) U_G(\A) Z(\A) \bs G(\A)} 
\int\limits_{\X(\A)} 
W_\varphi^\psi\left(\lambda \delta_\ell g\right) f_s^{(Z_n,\psi_n)}(g) \, d\lambda \, dg. 
\end{equation*} 
Finally, it follows from $Z_n U_G = N_G$ that 
\begin{equation*} 
\L(\varphi,f_s) = 
\int\limits_{N_G(\A) Z(\A) \bs G(\A)} 
\int\limits_{\X(\A)} 
W_\varphi^\psi\left(\lambda \delta_\ell g\right) f_s^{(Z_n,\psi_n)}(g) \, d\lambda \, dg. 
\end{equation*} 
This finishes the proof of Theorem \ref{thm-BasicIdentityII}. 
\end{proof}

%%%%%%%%%%%%%%%
\section{Constant Terms} \label{root-exchange-sec}
We now translate \cite[Section 7]{grs} into the $\GSpin$ context. 

%%%%%%%%%%%%%%%
\subsection{The Process of Root Exchanges}\label{SRE} 
We need the process of ``root exchanges'' of \cite{grs}.  A version of this  occurs in \cite{hs} as well,  but the Euler product expansion of 
Theorem \ref{thm-BasicIdentityII} requires the more elaborate version of \cite{grs}. Let us work in the following context. 

Let $H$ be a  connected reductive algebraic $k$-group, such as one of our $\GSpin$ groups. For our purposes we can have it be quasi-split.  
Let $U<G$ be a maximal unipotent $k$-subgroup.  Suppose $C<U$ is a $k$-subgroup of $U$, and $\psi=\psi_C$ is a non-trivial character 
of $C(k)\bs C(\A)$.  Suppose we have two other $k$-subgroups $X$ and $Y$ of $U$ such that the following six axioms hold:
%%%%%%%%%%%
\begin{enumerate}
\item $X$ and $Y$ normalize $C$. 
\item $X\cap C$ is normal in $X$ and $(X\cap C)\bs X$ is abelian, and similarly $Y\cap C$ is normal in $Y$ and $(Y\cap C)\bs Y$ is abelian. 
\item When $X(\A)$ and $Y(\A)$ act on $C(\A)$ by conjugation, they preserve $\psi_C$. 
\item $\psi_C$ is trivial on $(X\cap C)(\A)$ and $(Y\cap C)(\A)$. 
\item The commutator $(X,Y)\subset C$. (Recall that $(x,y)=x^{-1}y^{-1}xy$ and $(X,Y)$ denotes the subgroup generated by all the commutators.)
\end{enumerate}
%%%%%%%%%%%
Note that these five conditions imply that, for a fixed $y\in Y(\A)$, the map $x\mapsto \psi_C((x,y))$ defines a character of 
$X(\A)$, trivial on $(X\cap C)(\A)$ and similarly, for a fixed $x\in X(\A)$, the map $y\mapsto \psi_C((x,y))$ defines a character of $Y(\A)$ 
which is trivial on $(Y\cap C)(\A)$. (This is checked in \cite[\S 7.1]{grs}.) Finally
%%%%%%%%%%%
\begin{enumerate}
\item[(6)] The pairing of $(X\cap C)(\A)\bs X(\A) \times (Y\cap C)(\A)\bs Y(\A)$ given by 
$(x,y)\mapsto \psi_C((x,y))$ is bilinear and non-degenerate and identifies 
\begin{equation*}
(Y\cap C)(k)\bs Y(k)\simeq \left[(X(k)(X\cap C)(\A))\bs X(\A)\right]^{\wedge}  
\end{equation*} 
and 
\begin{equation*}
(X\cap C)(k)\bs X(k)\simeq \left[(Y(k)(Y\cap C)(\A))\bs Y(\A)\right]^\wedge.
\end{equation*} 
\end{enumerate}
%%%%%%%%%%%
Now let $B=CY=YC$, $D=CX=XC$ and $A=CXY$. Extend $\psi_C$ to a character $\psi_B$ of $B(\A)$ trivial on $B(k)$ by making it trivial on 
$Y(\A)$ and to a character $\psi_D$ of $D(\A)$ trivial on $D(k)$ by making it trivial on $X(\A)$.

%%%%%%%%%%%
\begin{prop}[Root exchange] \label{PRE} 
Let $f$ be an automorphic function on $H(\A)$ which is smooth and of uniform moderate growth. Then
\[
\int\limits_{B(k)\bs B(\A)}f(v)\psi_B(v)^{-1}\ dv=\int\limits_{(Y\cap C)(\A)\bs Y(\A)}\int\limits_{D(k)\bs D(\A)} f(uy)\psi_D(u)^{-1}\ du dy,
\]
where the right hand side converges in the sense that
\[
\int\limits_{(Y\cap C)(\A)\bs Y(\A)} \left| \int\limits_{D(k)\bs D(\A)}f(uyh)\psi_D(u)^{-1}\ du \right| dy <\infty
\]
uniformly in any compact subset of $H(\A)$.
\end{prop}
%%%%%%%%%%%

%%%%%%%%%%%
\begin{proof} 
The proof of this proposition in \cite{grs} depends only on the six axioms above and the properties of smooth automorphic functions 
of uniform moderate growth. As $G$ is connected reductive, the results of \cite{MW}, as used by \cite{grs}, hold. The proposition follows as in \cite{grs}.
\end{proof}
%%%%%%%%%%%

The proof of the proposition has the following corollary that is used in the proof of Proposition \ref{sum-int-exp}.

%%%%%%%%%%%
\begin{cor} \label{PREC} 
Let $f$ be a smooth automorphic function of uniform moderate growth on $H(\A)$. Then there exist smooth, uniform moderate growth, 
automorphic functions  
$f_1,\dots, f_r$ and Schwartz functions $\phi_1,\dots,\phi_r\in \mathcal S((Y\cap C)(\A)) \bs Y(\A)$ such that for all 
$y\in (Y\cap C)(\A)\bs Y(\A)$ we have
\[
\int\limits_{D(k)\bs D(\A)} f(uy)\psi^{-1}_D(u)\ du=\sum_{i=1}^r \phi_i(y)\int\limits_{D(k)\bs D(\A)} f_i(uy)\psi^{-1}_D(u)\ du.
\]
\end{cor}
%%%%%%%%%%%

%%%%%%%%%%%
\begin{proof} 
The proof in \cite{grs} uses the theory of smooth automorphic functions of uniform moderate growth from \cite{MW} as well as the theorem of 
Dixmier and Malliavin from \cite{DM}, both of which hold in our context. The result follows as in \cite{grs}.
\end{proof}

%%%%%%%%%%%%%%%
\subsection{Proof of Proposition \ref{sum-int-exp}}  \label{pf-prop}
We keep the notation from Section \ref{sec-integralA}.  We will use repeatedly that a cusp form $\varphi \in V_\pi$ on $H(\A)$ 
is smooth and of uniform moderate growth, so that we may use Section \ref{SRE}. 
The proof will be an induction, and for the induction we will need the following subgroups.

Recall that if $m'=2m+1$ is odd, then $n+\ell=m-1$ (by assumption) and $\GL_{n+\ell} \subset H$ is the Levi subgroup of the parabolic subgroup  
of $H$ with Levi subgroup $\GL_{n-\ell} \times H_3$ where $H_3$ is the split $\GSpin_{3}$. 
If $m'=2m$ is even then $n+\ell=m-1$ and $\GL_{m}$ is part of the Levi of the parabolic subgroup with Levi $\GL_{n+\ell} \times H_2$, 
where $H_2$ denotes the (split or quasi-split) $\GSpin(2)$. 
In either case, for $\gamma\in \GL_{n+\ell}$ we let $\gamma^\wedge$ denote the lift of $\gamma$ to an element of $H$. 
Since there is a corresponding parabolic subgroup of $\SO_{m'}$ we can follow \cite{grs} and write
\[
\gamma^\wedge=\bpm \gamma \\ & I_{m'-2(n+\ell)} \\ & & \gamma^*\epm\in \SO_{m'} 
\]
which we can identify with its lift to $H_{m'}$ if $\gamma$ is unipotent or a Weyl group element. Note that $m'-2(n+\ell)$ is either equal to $3$ or $2$.

We have defined 
\[
\X=\left\{\bpm I_n \\ \lambda & I_\ell\epm^\wedge\right\}
\]
and we set
\[
\X^{(i)}=\left\{\bpm I_n \\ \lambda & I_\ell\epm^\wedge \in\X \, : \, \lambda=\bpm \lambda_1 \\ \vdots \\ \lambda_\ell\epm \text{ with } \lambda_j=0 \text{ for } j\neq \ell-i\right\}
\]
and 
\[
\X_i=\left\{\bpm I_n \\ \lambda & I_\ell\epm^\wedge\in \X \, : \, \lambda_{\ell-i}=\cdots=\lambda_\ell=0\right\}
\]
Note that for each $i$, $\X_{i-1}=\X_{i}\X^{(i)}$ where we set $\X_{-1}=\X$. In what follows we will also need the groups 
\[
U_{n+\ell}^i=\left\{\bpm I_{n+\ell-i} & * \\ & z\epm^\wedge \, : \, z\in Z_i\right\}\cdot U_{n+\ell}. \subset N_{n+\ell}
\]

By definition we have
\[
\varphi^{(N_\ell,\Psi_\ell)}(g)=\int\limits_{N_\ell(k)\bs N_\ell(\A)}\varphi(vg)\Psi_\ell^{-1}(v)\ dv
\]
and then
\[
(\varphi^{(N_\ell,\Psi_\ell)})^{(U_G,1)}(g)=\int\limits_{U_G(k)\bs U_G(\A)}  
\left( \int\limits_{N_\ell(k)\bs N_\ell(\A)}\varphi(vrg)\Psi_\ell^{-1}(v)\ dv \right) \ dr.
\]
As $\varphi$ is automorphic on $H$, it is left invariant by $\delta_\ell\in H(k)$. Thus
\[
\varphi(vrg)=\varphi(\delta_\ell vrg)=\varphi(\delta_\ell vr \delta_\ell^{-1}\delta_\ell g).
\]

The authors of \cite{grs} analyze the conjugates $\delta_\ell N_\ell \delta_\ell^{-1}$ and $\delta_\ell U_G \delta_\ell^{-1}$ 
in terms of matrices in $\SO_{m'}$. To this end, let
\[
s(z;u, a, d, e;x,y)=\bpm I_n && x & d & y\\u & z & a & e & d'\\& & I_{m'-2(n+\ell)} & a' & x'\\& & & z^* & \\ & & & u'& I_n\epm
\]
where $z\in Z_\ell$, the maximal unipotent subgroup of $\GL_{\ell}$. With this notation,
\[
\begin{aligned}
\delta_\ell U_G \delta_\ell^{-1}&=\{ s(I_\ell; 0,0,0,0,;x^\circ,y)\in H \, : \, \delta_\ell^{-1}s(I_\ell; 0,0,0,0;x^\circ,y)\delta_\ell\in U_G\}\\
& =\left\{\bpm I_n && x^\circ & 0 & y\\0 & I_\ell & 0 & 0 & 0\\& & I_{m'-2(n+\ell)} & 0 & x^{\circ \prime}\\& & & I_\ell & \\ & & & 0& I_n\epm\right\}
\end{aligned}
\]
where there is a condition on $x^\circ$ to make sure this comes from an element of $U_G$, and
\[
\begin{aligned}
\delta_\ell N_\ell \delta_\ell^{-1}&=\{s(z; u,a,d,e;0,0)\, : \, z\in Z_\ell\}\\
&=\left\{\bpm I_n && 0 & d & 0\\u & z & a & e & d'\\& & I_{m'-2(n+\ell)} & a' &0\\& & & z^* & \\ & & & u'& I_n\epm\right\}.
\end{aligned}
\]
Let 
\[
S=\delta_\ell U_G N_\ell \delta_\ell^{-1}
\]
and let the character $\Psi_S$ be the character of $S$ obtained as follows: we extend the character 
$\Psi_\ell$ of $N_\ell(k)\bs N_\ell(\A)$ to $U_G N_\ell$ by making it trivial on $U_G(\A)$ and then set $\Psi_S(s)=\Psi_\ell(\delta_\ell^{-1} s \delta_\ell)$.
Let
\[
\begin{aligned}
\tilde{S}&=\{ s\in S\, : \, u=0, z=I_\ell\}\\
&=\left\{\bpm I_n && x^\circ & d & y\\0 & I_\ell & a & e & d'\\& & I_{m'-2(n+\ell)} & a' & x^{\circ\prime}\\& & & I_\ell & \\ & & & 0& I_n\epm\right\}.
\end{aligned}
\]
Note that $\X\subset S$ consists of the lower triangular elements of $S$, that is, $\X=\{ s(0;u,0,0,0;0,0)\}$.

%%%%%%%%%%%
\begin{lem} 
$\Psi_S$ is trivial on $\X$.
\end{lem}
%%%%%%%%%%%

%%%%%%%%%%%
\begin{proof}  
By definition $\Psi_S(s)=\Psi_\ell(\delta_\ell^{-1}s\delta_\ell)$. 
We have $\delta_\ell=\bpm 0 & I_n \\ I_\ell & 0 \epm^\wedge$ and  $\X=\left\{\bpm I_n & 0 \\ \lambda & I_\ell\epm^\wedge \right\}$.  
Letting $\tilde\lambda=\bpm I_n & 0 \\ \lambda & I_\ell\epm$ we find that $\Psi_S(\tilde\lambda)=\Psi_\ell(\delta_\ell^{-1}\tilde\lambda\delta_\ell)$. 
However
\[
\delta_\ell^{-1}\tilde\lambda\delta_\ell = 
\left(\bpm 0 & I_\ell \\  I_n & 0 \epm\bpm I_n & 0 \\ \lambda & I_\ell  \epm\bpm 0 & I_n \\ I_\ell & 0\epm\right)^\wedge 
= \bpm I_\ell & \lambda \\ 0 & I_n\epm^\wedge
\]
and
\[
\Psi_\ell\left(\bpm I_\ell & \lambda \\ 0 & I_n\epm^\wedge\right) 
= 
\begin{cases} 
\psi_0(\lambda_{\ell,m} - \lambda_{\ell, m+1}) & \mbox{if } m'=2m, \\ 
\psi_0(\lambda_{\ell,m}+\frac{a}{2}\lambda_{\ell, m+2}) & \mbox{if } m'=2m+1. 
\end{cases}
\]
by the definitions in Section \ref{sec-integralB}.  If $m'=2m+1$,  then $n+\ell=m-1$ and so $\lambda_{\ell, m}=\lambda_{\ell, m+2}=0$.  
If $m'=2m$, then $n+\ell=m-1$ and again $\lambda_{\ell, m}=\lambda_{\ell, m+1}=0$.  
(For these conditions on the relation between $n, \ell, m$ see the second paragraph of Section \ref{pf-prop}.)
Hence $\Psi_S(\tilde\lambda)=1$.
\end{proof}
%%%%%%%%%%%

With this notation we can write
\begin{equation} \label{phi-exp}
\begin{aligned}
(\varphi^{(N_\ell,\Psi_\ell)})^{(U_G,1)}(g) 
&= \int\limits_{U_G(k) \bs U_G(\A)}\left[ \int\limits_{N_\ell(k)\bs N_\ell(\A)}\varphi(vrg)\Psi_\ell^{-1}(v)\ dv\right] \ dr \\
&=\int\limits_{S(k) \bs S(\A)} \varphi(s\delta_\ell g)\Psi_S^{-1}(s)\ ds \\
&=\int\limits_{Z_\ell(k) \bs Z_\ell(\A)} \int\limits_{\X(k) \bs \X(\A)} \int\limits_{\tilde{S}(k) \bs \tilde{S}(\A)} 
\varphi(z\lambda \tilde{s}\delta_\ell g)\Psi_S^{-1}(\tilde{s})\Psi_{n+\ell}^{-1}( z) \, d\tilde{s} \, d\lambda \, dz .
\end{aligned}
\end{equation}

Since $\X=\X_0\X^{(0)}$ we can decompose the integral over $\X$ as
\[
\int\limits_{\X(k)\bs \X(\A)}=\int\limits_{\X_0(k)\bs \X_0(\A)}\int\limits_{\X^{(0)}(k)\bs \X^{(0)}(\A)}.
\]
Then we obtain an inner integral of the form
\[
\int\limits_{\X^{(0)}(k)\bs \X^{(0)}(\A)}\int\limits_{\tilde{S}(k) \bs \tilde{S}(\A)} 
\varphi(z\lambda_0\lambda^{(0)} \tilde{s}\delta_\ell g)\Psi_S^{-1}(\tilde{s})\Psi_{n+\ell}^{-1}( z)\ d\tilde{s}d\lambda^{(0)}
\]
or
\[
\int\limits_{\X^{(0)}(k)\bs \X^{(0)}(\A)}\int\limits_{\tilde{S}(k)\bs\tilde{S}(\A)} 
\left[\varphi(z\lambda_0\lambda^{(0)} \tilde{s}\delta_\ell g)\Psi_{n+\ell}^{-1}( z)\right]\Psi_S^{-1}(\tilde{s})\ d\tilde{s}d\lambda^{(0)}.
\]

We would now like to perform a root exchange on this integral. To this end, let 
$J=\{s(I_\ell; 0,0,0,0;x,y)\in H\}\supset \delta_\ell U_G \delta_\ell^{-1}$ and 
$J_0=J\cap S=\delta_\ell U_G\delta_\ell^{-1}$.  Now in the setting of \S\ref{SRE}, 
set $C=\tilde{S}$, $Y=\X^{(0)}$, and  $X=J$ so that $B=\tilde{S}\X^{(0)}$,  
$D=\tilde{S}J=U_{n+\ell}$, and $A=\tilde{S}\X^{(0)}J$. 
The authors of \cite{grs} verify that these satisfy (1) -- (6) of \S\ref{SRE}.   
Applying Proposition \ref{PRE} to this integral then gives
\[
\begin{aligned}
\int\limits_{\X^{(0)}(k) \bs \X^{(0)}(\A)} 
& \int\limits_{\tilde{S}(k)\bs\tilde{S}(\A)} 
\left[\varphi(z\lambda_0\lambda^{(0)} \tilde{s}\delta_\ell g)\Psi_{n+\ell}^{-1}( z)\right] 
\Psi_S^{-1}(\tilde{s})\ d\tilde{s}d\lambda^{(0)} \\
&= \int\limits_{\X^{(0)}(\A)} \int\limits_{U_{n+\ell}(k) \bs U_{n+\ell}(\A)}  
\left[\varphi(uz\lambda_0\lambda^{(0)} \delta_\ell g)\Psi_{n+\ell}^{-1}( z)\right]\Psi^{-1}_{n+\ell}(u)\ dud\lambda^{(0)} \\
&=\int\limits_{\X^{(0)}(\A)} \int\limits_{U_{n+\ell}(k)\bs U_{n+\ell}(\A)}  
\varphi(uz\lambda_0\lambda^{(0)} \delta_\ell g)\Psi_{n+\ell}^{-1}( zu)\ dud\lambda^{(0)}
\end{aligned}
\] 
with convergence of 
\[
\int\limits_{\X^{(0)}(\A)} \left| \, \int\limits_{U_{n+\ell}(k) \bs U_{n+\ell}(\A)}  
\varphi(uz\lambda_0\lambda^{(0)} \delta_\ell g)\Psi_{n+\ell}^{-1}(zu)\ du \right| d\lambda^{(0)}
\]
uniformly on compact subsets.

Moreover, by Corollary \ref{PREC} there are smooth automorphic functions 
$\varphi_1,\dots,\varphi_r$ of uniform moderate growth on $H(\A)$ and Schwartz functions 
$\phi_1,\dots,\phi_r$ on $\X^{(0)}(\A)$ such that
\[
\begin{aligned}
\int\limits_{U_{n+\ell}(k)\bs U_{n+\ell}(\A)} 
& \varphi(uz\lambda_0\lambda^{(0)}\delta_\ell g)\Psi_{n+\ell}^{-1}( zu)\ du \\
&=\sum_{i=1}^r\phi_i(z\lambda_0\lambda^{(0)}\delta_\ell g)\int\limits_{U_{n+\ell}(k)\bs U_{n+\ell}(\A)} 
\varphi_i(uz\lambda_0\lambda^{(0)}\delta_\ell g)\Psi_{n+\ell}^{-1}( zu)\ du.\\
\end{aligned}
\]
If we insert this in the result of our root exchange we have
\[
\begin{aligned}
\int\limits_{\X^{(0)}(\A)} &\int\limits_{U_{n+\ell}(k)\bs U_{n+\ell}(\A)}  
\varphi(uz\lambda_0\lambda^{(0)} \delta_\ell g)\Psi_{n+\ell}^{-1}( zu)\ dud\lambda^{(0)} \\
&=\int\limits_{\X^{(0)}(\A)} \left[\sum_{i=1}^r\phi_i(z\lambda_0\lambda^{(0)}\delta_\ell g) 
\int\limits_{U_{n+\ell}(k)\bs U_{n+\ell}(\A)} 
\varphi_i(uz\lambda_0\lambda^{(0)}\delta_\ell g)\Psi_{n+\ell}^{-1}( zu)\ du\right] d\lambda^{(0)} \\
\end{aligned}
\]
\[
= \int\limits_{U_{n+\ell}(k)\bs U_{n+\ell}(\A)}\left[\sum_{i=1}^r \int\limits_{\X^{(0)}(\A)} 
\phi_i(z\lambda_0\lambda^{(0)}\delta_\ell g) 
\varphi_i(uz\lambda_0\lambda^{(0)}\delta_\ell g) \Psi_{n+\ell}^{-1}( zu)\ d\lambda^{(0)}\right]  du
\]
so if we set
\[
\varphi_0(u\lambda_0z\delta_\ell g) = \sum_{i=1}^r \int\limits_{\X^{(0)}(\A)} 
\phi_i(z\lambda_0\lambda^{(0)}\delta_\ell g)\varphi_i(uz\lambda^{(0)}\lambda_0z\delta_\ell g) \ d\lambda^{(0)},
\]
then we have 
\[
\int\limits_{U_{n+\ell}(k)\bs U_{n+\ell}(\A)} \varphi_0(uz\lambda_0\delta_\ell g)\Psi_{n+\ell}^{-1}( zu)\  du.
\]

If we return now to \eqref{phi-exp}, noting the decomposition of integration right after \eqref{phi-exp}, and insert this we have
\[ 
\left( \varphi^{(N_\ell,\Psi_\ell)} \right)^{(U_G,1)}(g) 
= \int\limits_{Z_\ell(k)\bs Z_\ell(\A)}\int\limits_{\X_0(k)\bs \X_0(\A)}\int\limits_{U_{n+\ell} (k)\bs U_{n+\ell}(\A)} 
\varphi_0(uz\lambda_0 \delta_\ell g)\Psi_{n+\ell}^{-1}( zu) \, du \, d\lambda_0 \, dz .
\]

We next use the decompositions $\X_0=\X_1\X^{(1)}$ and $Z_\ell=Z_{\ell-1}Z_\ell^{(1)}$, where
\[
Z_\ell^{(1)}=\left\{ \bpm I_{\ell-1} & * \\ 0 & 1\epm\right\}\subset Z_\ell.
\] 
To this end, we write $\lambda_0 = \lambda_1 \lambda^{(1)}$ and, with abuse of notation, 
$z=z_\ell=z_{\ell-1}z^{(1)}=zz^{(1)}$ and decompose our integrals as
\[
\int\limits_{\X_0(k)\bs \X_0(\A)}=\int\limits_{\X_1(k)\bs \X_1(\A)}\int\limits_{\X^{(1)}(k)\bs \X^{(1)}(\A)}
\]
and
\[
\int\limits_{Z_\ell(k)\bs Z_\ell(\A)}=\int\limits_{Z_{\ell-1}(k)\bs Z_{\ell-1}(\A)} \int\limits_{Z_\ell^{(1)}(k)\bs Z_\ell^{(1)}(\A)}.
\]
For $h\in H(\A)$ set
\[ 
\begin{aligned} 
& I_1(\varphi_0,\Psi_{n+\ell})(h) = \\ 
& 
\int\limits_{\X^{(1)}(k)\bs \X^{(1)}(\A)}\int\limits_{Z_\ell^{(1)}(k)\bs Z_\ell^{(1)}(\A)} 
\int\limits_{U_{n+\ell}(k)\bs U_{n+\ell}(\A)} 
\varphi_0(u z^{(1)}\lambda^{(1)} h)\Psi_{n+\ell}^{-1}(uz^{(1)}) \, du \, dz^{(1)} \, d\lambda^{(1)}
\end{aligned}
\] 
so that
\[
\left(\varphi^{(N_\ell,\Psi_\ell)}\right)^{(U_G,1)}(g) = 
\int\limits_{Z_{\ell-1}(k)\bs Z_{\ell-1}(\A)}\int\limits_{\X_1(k)\bs \X_1(\A)} 
I_1(\varphi_0,\Psi_{n+\ell})(z\lambda_1 \delta_\ell g)\Psi_{n+\ell}^{-1}( z) \, d\lambda_1 \, dz .
\]

At this point, we apply another root exchange (or swap in the languate of \cite{hs}) in $I_1$.  
In the language of Section \ref{SRE} we set
\[
\begin{aligned}
C &=U_{n+\ell}Z^{(1)} \quad\quad &Y=\X^{(1)} \quad\quad &X=\left\{ \bpm I_n & & x\\& I_{\ell-1} & \\ & & 1\epm^\wedge\right\} \\
\psi_C&=\Psi_{n+\ell}|_C \quad\quad & B=CY \quad\quad &D=CX=U_{n+\ell}^1.
\end{aligned}
\]
Applying this to $I_1(\varphi_0, \Psi_{n+\ell})$ gives
\[
I_1(\varphi_0,\Psi_{n+\ell})(h) = \int\limits_{\X^{(1)}(\A)}\int\limits_{U^1_{n+\ell}(k)\bs U^1_{n+\ell}(\A)} 
\varphi_0(u^1\lambda^{(1)}h)\Psi_{n+\ell}^{-1}(u^1)\ du^1\ d\lambda^{(1)}.
\]
We insert this into the above formula for $(\varphi^{(N_\ell,\Psi_\ell)})^{(U_G,1)}(g)$ and interchange the order of integration 
to obtain
\[
\begin{aligned}
\left(\varphi^{(N_\ell,\Psi_\ell)}\right)^{(U_G,1)}(g) 
& =\int\limits_{\X^{(1)}(\A)} \int\limits_{Z_{\ell-1}(k)\bs Z_{\ell-1}(\A)} 
\int\limits_{\X_1(k)\bs\X_1(\A)}\int\limits_{U^1_{n+\ell}(k)\bs U^1_{n+\ell}(\A)}  \\ 
&\quad \quad \varphi_0(u^1\lambda_1z\lambda^{(1)})\Psi^{-1}_{n+\ell}(u^1z) \, du^1 \, d\lambda_1 \, dz \, d\lambda^{(1)}.
\end{aligned}
\]
We next apply Corollary \ref{PREC} (essentially Dixmier--Malliavin \cite{DM}) to write
\[
\begin{aligned}
& \int\limits_{U^1_{n+\ell}(k)\bs U^1_{n+\ell}(\A)}  
\varphi_0 \left(u^1(\lambda_1z\lambda^{(1)}\delta_\ell g)\right) \Psi^{-1}(u^1z) \, du & \\ 
& =\sum_{I=1}^r \phi_{0,i}(\lambda_1z\lambda^{(1)}\delta_\ell g) 
\int\limits_{U^1_{n+\ell}(k)\bs U^1_{n+\ell}(\A)} 
\varphi_{0,i} \left(u^1\lambda_1z\lambda^{(1)}\delta_\ell g \right) \Psi^{-1}_{n+\ell}(u^1z) \, du^1.
\end{aligned}
\] 
Note that the $\varphi_{0,i}$ are smooth and of uniform moderate growth (in fact in $W(\varphi)$) 
and the $\phi_{0,i}$ are Schwartz functions. If we let
\[
\varphi_1 \left(u^1\lambda_1zh \right) 
= \sum_{I=1}^r \int\limits_{\X^{(1)}(\A)} 
\varphi_{0,i}\left( u^1\lambda_1z\lambda^{(1)}h \right) \phi_{0,i}\left(\lambda_1z\lambda^{(1)}h \right) \, d\lambda^{(1)}, 
\] 
then by Corollary \ref{PREC} we have
\[
\begin{aligned} 
&\int\limits_{\X^{(1)}(\A)} \int\limits_{U^1_{n+\ell}(k) \bs U^1_{n+\ell}(\A)}   
\varphi_0\left(u^1(\lambda_1z\lambda^{(1)}\delta_\ell g)\right) \Psi^{-1}(u^1z) \, du \, d\lambda^{(1)} \\ 
& = \int\limits_{U^1_{n+\ell}(k)\bs U^1_{n+\ell}(\A)} \varphi_1(u^1(\lambda_1z\delta_\ell g))\Psi^{-1}_{n+\ell}(u^1z)\ du^1
\end{aligned}
\]
and so 
\[
\begin{aligned}
(\varphi^{(N_\ell,\Psi_\ell)})^{(U_G,1)}(g) 
&= \int\limits_{Z_{\ell-1}(k)\bs Z_{\ell-1}(\A)}\int\limits_{\X_1(k)\bs \X_1(\A)}\int\limits_{U^1_{n+\ell}(k)\bs U^1_{n+\ell}(\A)}  \\ 
&\quad\quad 
\varphi_1(u^1\lambda_1z\delta_\ell g)\Psi^{-1}_{n+\ell}(u^1z)\ du^1\ d\lambda_1\ dz.
\end{aligned}
\]
%

%%%%%%%%%%%
\begin{rem}
In \cite{grs} they keep this integral, as well as the original integral involving $\varphi$. Then this would be
\[
\begin{aligned}
(\varphi^{(N_\ell,\Psi_\ell)})^{(U_G,1)}(g) = & 
\int\limits_{\X^{(1)}(\A)}\int\limits_{\X^{(0)}(\A)} 
\int\limits_{Z_{\ell-1}(k)\bs Z_{\ell-1}(\A)} 
\int\limits_{\X_1(k)\bs \X_1(\A)} 
\int\limits_{U^1_{n+\ell}(k)\bs U^1_{n+\ell}(\A)}  \\ 
& \varphi(u^1\lambda_1z\lambda^{(0)}\lambda^{(1)}\delta_\ell g) \Psi^{-1}(u^1z)
 \, du^1 \, d\lambda_1 \, dz \, d\lambda^{(0)} \, d\lambda^{(1)}.
\end{aligned}
\] 
They perform  $\varphi\mapsto \varphi_0\mapsto \varphi_1$ to obtain convergence and 
they want the original integral with $\varphi$ and the integrations over the $\X^{(i)}$ to obtain 
an Euler product for decomposable data (i.e., $\varphi$ decomposable) and  need the full 
adelic integrals over the $\X^{(i)}$.
\end{rem}
%%%%%%%%%%%

We next do an induction to either remove the various $\X_i$ integrations or replace them with full adelic integrations over the $\X^{(i)}$.

%%%%%%%%%%%
\subsubsection{The Induction} 
For $1\leq i <\ell$ set 
\[
Z_{\ell-i}=\left\{ z_{\ell-i}=\bpm z & \\ & I_i\epm\in Z_\ell\right\}\subset H.
\]
Assume, by induction, that for $1\leq i \leq \ell-2$ we have:
\[
\begin{aligned}
(i)\quad\quad \left(\varphi^{(N_\ell,\Psi_\ell)}\right)^{(U_G,1)}(g) = & 
\int\limits_{\X^{(i,\dots, 1,0)}(\A)} \int\limits_{Z_{\ell-i}(k)\bs Z_{\ell-i}(\A)} 
\int\limits_{\X_i(k)\bs \X_i(\A)} 
\int\limits_{U^i_{n+\ell}(k)\bs U^i_{n+\ell}(\A)}  \\ 
& 
\varphi\left(u^i\lambda_iz_{\ell-i}\lambda^{(\prime)}\delta_\ell g\right) 
\Psi^{-1}(u^iz_{\ell-i}) \ du^i \ d\lambda_i\ dz_{\ell-i}\ d\lambda^{(\prime)},
\end{aligned}
\]
where $\lambda^{(\prime)}\in \X^{(i,\dots,1,0)}=\X^{(i)} \X^{(i-1)}\cdots\X^{(1)}\X^{(0)}$ and 
$d\lambda^{(\prime)}=d\lambda^{(i)}\cdots d\lambda^{(0)}$, and
\[
\begin{aligned}
(ii)\quad\quad \left(\varphi^{(N_\ell,\Psi_\ell)}\right)^{(U_G,1)}(g) = & 
\int\limits_{Z_{\ell-i}(k)\bs Z_{\ell-i}(\A)} 
\int\limits_{\X_i(k)\bs \X_i(\A)} 
\int\limits_{U^i_{n+\ell}(k)\bs U^i_{n+\ell}(\A)}  \\ 
& 
\quad\quad \varphi_i\left(u^i\lambda_iz_{\ell-i}\delta_\ell g\right) 
\Psi^{-1}_{n+\ell}(u^iz_{\ell-i})\ du^i\ d\lambda_i\ dz_{\ell-i}.
\end{aligned}
\]
with $\varphi_i \in W(\varphi)$.

We then repeat the above process  to obtain the next step, 
formulas \cite[(7.29)--(7.30)]{grs}.  For $i=\ell-2$ these formulas give 
\[
(i)\quad\quad \left(\varphi^{(N_\ell,\Psi_\ell)}\right)^{(U_G,1)}(g) =  
\int\limits_{\X(\A)} \int\limits_{U^{\ell-1}_{n+\ell}(k)\bs U^{\ell-1}_{n+\ell}(\A)}  
\varphi\left(u^{\ell-1}\lambda\delta_\ell g\right) \Psi^{-1}(u^{\ell-1}) \, du^{\ell-1} \, d\lambda,
\] 
and
\[
\begin{aligned}
(ii) \quad\quad \left(\varphi^{(N_\ell,\Psi_\ell)}\right)^{(U_G,1)}(g) &= 
\int\limits_{U^{\ell-1}_{n+\ell}(k)\bs U^{\ell-1}_{n+\ell}(\A)} 
\varphi_{\ell-1}\left(u^{\ell-1}\delta_\ell g\right) \Psi^{-1}_{n+\ell}(u^{\ell-1}) \, du^{\ell-1}. \\
\end{aligned}
\]
%

%%%%%%%%%%%
\begin{lem} 
Let $p\in P^1_{n-1,1}(\A) \subset \GL_{n}(\A)$. (This is the notation of \cite{grs} for the mirabolic subgroup of $\GL_{n}$.) Then 
\[
\begin{aligned} 
& \int\limits_{\X(\A)} \int\limits_{U^{\ell-1}_{n+\ell}(k)\bs U^{\ell-1}_{n+\ell}(\A)}  
\varphi\left(u^{\ell-1}\hat{p}\lambda\delta_\ell g\right) \Psi^{-1}(u^{\ell-1}) \, du^{\ell-1} \, d\lambda \\
&= \int\limits_{U^{\ell-1}_{n+\ell}(k)\bs U^{\ell-1}_{n+\ell}(\A)} 
\varphi_{\ell-1}\left(u^{\ell-1}\hat{p} \delta_\ell g\right) \Psi^{-1}_{n+\ell}(u^{\ell-1})\ du^{\ell-1}.
\end{aligned} 
\]
\end{lem}
%%%%%%%%%%%

%%%%%%%%%%%
\begin{proof} 
The proof proceeds in the same way as in \cite{grs}. 
\end{proof}
%%%%%%%%%%%

%%%%%%%%%%%
\begin{rem} 
%\noindent{\bf Remark:} 
In their derivation, the authors of \cite{grs}  only keep track of 
$\left(\varphi^{(N_\ell,\Psi_\ell)}\right)^{(U_G,1)}(e)$ and so $g$ does not appear in their formulas.
\end{rem} 
%%%%%%%%%%%

For $x\in \A^n$ let $n_\ell(x)=\bpm I_n & x \\ & 1\epm^\wedge\in N_{n+\ell}(\A)$. Set
\[
\phi_g(x)=\int\limits_{U^{\ell-1}_{n+\ell}(k)\bs U^{\ell-1}_{n+\ell}(\A)} 
\varphi_{\ell-1}\left(u^{\ell-1}n_\ell(x)\delta_\ell g\right) \Psi^{-1}_{n+\ell}(u^{\ell-1}) \, du
\]
so that
\[
\phi_g(0) = \left(\varphi^{(N_\ell,\Psi_\ell)}\right)^{(U_G,1)}(g).
\]
This is a smooth function on $ \A^n$ which is left invariant under $k^n$.  
We can write its Fourier expansion along $k^n\bs \A^n$ and then evaluate at $x=0$.  
The general character of $\A^n$ which is trivial on $k^n$ is of the form 
$\psi({^t\eta}\cdot x)$ for $\eta\in k^n$. By abelian Fourier analysis
 \[
 \phi_g(x)=\sum_{\eta\in k^n} a_\eta(\phi_g)\psi({^t\eta}\cdot x)
 \]
 where
 \[
 a_\eta(\phi_g)=\int\limits_{k^n\bs\A^n}\phi_g(x)\psi^{-1}({^t\eta}\cdot x)\ dx
 \]

\noindent{\bf The zero Fourier coefficient.} 
The Fourier coefficient corresponding to the trivial character $\eta=0$ is:
\[
a_0(\phi_g) = \int\limits_{k^n\bs \A^n} \int\limits_{U^{\ell-1}_{n+\ell}(k)\bs U^{\ell-1}_{n+\ell}(\A)} 
\varphi_{\ell-1}\left(u^{\ell-1}n_\ell(x)\delta_\ell g\right) 
\Psi^{-1}_{n+\ell}(u^{\ell-1}) \, du \, dx.
\] 
Note that $\{n_\ell(x)\mid x\in \A^n\}\cdot U_{n+\ell}^{\ell-1}=U^\ell_{n+\ell}\simeq N_\ell \ltimes U_n$, 
where we view $N_\ell$ as being in the $\GSpin$ part of the Levi subgroup of 
$P_n\simeq (\GL_{n} \times \GSpin)\ltimes U_n$, and $\Psi_{n+\ell}$ on $U_{n+\ell}^{\ell-1}$ equals 
$\psi_\ell$ on $N_\ell$, so that this $0$-Fourier coefficient gives
\[
\int\limits_{k^n\bs \A^n} \int\limits_{U^{\ell-1}_{n+\ell}(k)\bs U^{\ell-1}_{n+\ell}(\A)} 
\varphi_{\ell-1}\left(u^{\ell-1}n_\ell(x)\delta_\ell g\right) \Psi^{-1}_{n+\ell}(u^{\ell-1}) \, du \, dx 
= \left(\varphi_{\ell-1}^{(U_n,1)}\right)^{(N_\ell, \psi_\ell)}(\delta_\ell g).
\]

\noindent{\bf The non-zero Fourier coefficients.}  
If $\eta\neq 0$, then we can write ${^t\eta}=(0, \dots, 0,1)\gamma$ where  
$\gamma\in P^1_{n-1}(k)\bs \GL_{n}(k)$ with $P^1_{n-1,1} = \operatorname{Stab}_{\GL_{n}}(0, \dots, 0,1)$ 
is the mirabolic subgroup of $\GL_{n}$.  We can manipulate this Fourier coefficient into
\[
\begin{aligned} 
a_\eta(\phi_g) 
&= \int\limits_{k^n\bs \A^n} \phi_g(n_\ell(x))\psi^{-1}(^t\eta\cdot x) \, dx \\ 
&= \int\limits_{k^n\bs\A^n} \int\limits_{U^{\ell-1}_{n+\ell}(k)\bs U^{\ell-1}_{n+\ell}(\A)} 
\varphi_{\ell-1}\left(u^{\ell-1}n_\ell(x)\delta_\ell g\right) \Psi_{n+\ell}^{-1}(u^{\ell-1}n_\ell(\gamma x)) \, du \, dx.
\end{aligned}
\]

Since $\varphi_{\ell-1}$ is still automorphic, it is left invariant under $\hat{\gamma}$ since $\gamma\in \GL_{n}(k)$. Then
\[
\varphi_{\ell-1}\left(un_\ell(x)\delta_\ell g\right) 
= \varphi_{\ell-1}\left(\hat{\gamma}un_\ell(x)\delta_\ell g\right) 
= \varphi_{\ell-1}\left(\hat{\gamma}u\hat{\gamma}^{-1}\hat{\gamma}n_\ell(x)\delta_\ell g\right) 
=\varphi_{\ell-1}\left(\hat{\gamma}u\hat{\gamma}^{-1}n_\ell(\gamma x)\hat{\gamma}\delta_\ell g\right).
\] 
So the $\eta$-Fourier coefficient of $\phi_g$ is then
\[
\begin{aligned} 
a_\eta(\phi_g) 
&= \int\limits_{k^n\bs \A^n} \phi_g(n_\ell(x))\psi^{-1}(^t\eta\cdot x) \, dx\\ 
&= \int\limits_{k^n\bs\A^n} \int\limits_{U^{\ell-1}_{n+\ell}(k)\bs U^{\ell-1}_{n+\ell}(\A)} 
\varphi_{\ell-1}\left(\hat{\gamma}u\hat{\gamma}^{-1}n_\ell(\gamma x)\hat{\gamma}\delta_\ell g\right) 
\Psi_{n+\ell}^{-1}(un_\ell(\gamma x)) \, du \, dx.
\end{aligned}
\]
Now perform a change of variables $u\mapsto \hat{\gamma}^{-1}u\hat{\gamma}$ and 
$x\mapsto \gamma^{-1}x$ to obtain that the $\eta$-Fourier coefficient of $\phi_g$ is given by
\[
\begin{aligned}
a_\eta(\phi_g) 
&= \int\limits_{k^n\bs \A^n} \phi_g(n_\ell(x)) \psi^{-1}(^t\eta\cdot x) \, dx \\ 
&= \int\limits_{k^n\bs\A^n} \int\limits_{U^{\ell-1}_{n+\ell}(k)\bs U^{\ell-1}_{n+\ell}(\A)} 
\varphi_{\ell-1}\left(un_\ell(x)\hat{\gamma}\delta_\ell g\right) \Psi_{n+\ell}^{-1}(un_\ell(x)) \, du \, dx.
\end{aligned}
\]
We once again write $\left\{n_\ell(x)\mid x\in \A^n \right\} \cdot U_{n+\ell}^{\ell-1}=U^\ell_{n+\ell}$ and this becomes
\[
a_\eta(\phi_g) = \int\limits_{k^n\bs \A^n} 
\phi_g(n_\ell(x)) \psi^{-1}( ^t\eta\cdot x ) \, dx 
= \int\limits_{U^{\ell}_{n+\ell}(k)\bs U^{\ell}_{n+\ell}(\A)} 
\varphi_{\ell-1}\left( u\hat{\gamma}\delta_\ell g \right) \Psi_{n+\ell}^{-1}(u) \, du \, dx,
\]
where ${^t\eta}=(0,\dots,0,1)\gamma$

If we now combine these expressions, we have
\[
\begin{aligned}
&\left(\varphi^{(N_\ell,\Psi_\ell)}\right)^{(U_G,1)}(g) = \phi_g(0) = \sum_\eta a_\eta(\phi_g) \\ 
&= \left(\varphi_{\ell-1}^{(U_n,1)}\right)^{(N_\ell, \psi_\ell)}(\delta_\ell g) + 
\sum\limits_{\gamma\in P^1_{n-1,1}(k)\bs \GL_{n}(k)} \int\limits_{U^{\ell}_{n+\ell}(k)\bs U^{\ell}_{n+\ell}(\A)} 
\varphi_{\ell-1}\left(u\hat{\gamma}\delta_\ell g\right) \Psi_{n+\ell}^{-1}(u) \, du \, dx.
\end{aligned}
\]

We need to work in the above form for convergence (see the above remark). In order to prove our integrals are Eulerian, 
we need to express this in terms of our original cusp form $\varphi$.  This happens Fourier coefficient by Fourier coefficient. 
For a non-trivial Fourier coefficient we have
\[
\begin{aligned}
a_\eta(\phi_g) 
&= \int\limits_{U^\ell_{n+\ell}(k)\bs U^\ell_{n+\ell}(\A)} 
\varphi_{\ell-1}\left(u\hat{\gamma}\delta_\ell g\right) \Psi^{-1}_{n+\ell}(u) \, du \\
&=\int\limits_{k^n\bs \A^n}\int\limits_{U^{\ell-1}_{n+\ell}(k)\bs U^{\ell-1}_{n+\ell}(\A)} 
\varphi_{\ell-1}\left(un_\ell(x)\hat{\gamma}\delta_\ell g\right) \Psi^{-1}_{n+\ell}(un_\ell(x)) \, du \, dx \\
&=\int\limits_{k^n\bs \A^n}\int\limits_{\X(\A)}\int\limits_{U^{\ell-1}_{n+\ell}(k)\bs U^{\ell-1}_{n+\ell}(\A)} 
\varphi\left(un_\ell(x)\hat{\gamma}\lambda \delta_\ell g\right) \Psi^{-1}_{n+\ell}(un_\ell(x)) \, du \, d\lambda \, dx \\
&=\int\limits_{\X(\A)}\int\limits_{k^n\bs \A^n} \int\limits_{U^{\ell-1}_{n+\ell}(k)\bs U^{\ell-1}_{n+\ell}(\A)} 
\varphi\left(un_\ell(x)\hat{\gamma}\lambda \delta_\ell g\right) \Psi^{-1}_{n+\ell}(un_\ell(x)) \, du \, dx \, d\lambda\\
&=\int\limits_{\X(\A)}\int\limits_{U^{\ell}_{n+\ell}(k)\bs U^{\ell}_{n+\ell}(\A)} 
\varphi\left(u\hat{\gamma}\lambda \delta_\ell g\right) \Psi^{-1}_{n+\ell}(u) \, du \, dx \, d\lambda\\
\end{aligned}
\]
and the authors of \cite{grs} go to great lengths to justify the exchange of the integral over $k^n\bs \A^n$ and that over $\X(\A)$.  
Along the way, in the successive interchanges, they use
\begin{enumerate}
\item[$\bullet$] $\Psi^{-1}_{n+\ell}([x,\lambda])=\Psi_{n+\ell}([x^{-1},\lambda])$ 
for $x\in X(\A)$ where $X=\left\{ \bpm I_n & 0 & x\\& 1 & 0\\ & & 1\epm^\wedge\in N_{n+\ell}\right\}$
\item[$\bullet$] $\hat{\gamma}[x,\lambda]=[x,\lambda]\hat{\gamma}$
\item[$\bullet$] $\hat{\gamma}x\hat{\gamma}^{-1}\in U_{n+\ell}^\ell(\A)$
\item[$\bullet$] $\Psi_{n+\ell}(\hat{\gamma}x\hat{\gamma}^{-1})=1$.
\end{enumerate}
These all involve either unipotent matrices of lifts of elements of $\GL_{n}$ and remain valid in the $\GSpin$ context.

Similarly, for the $0$-Fourier coefficient they show
\[
\left(\varphi_{\ell-1}^{(U_n,1)}\right)^{(N_\ell, \psi_\ell)}(\delta_\ell g) 
= \int\limits_{\X(\A)} \left(\varphi^{(U_n,1)}\right)^{(N_\ell, \psi_\ell)}(\lambda \delta_\ell g) \, d\lambda.
\]
Note that since $\varphi$ is assumed to be cuspidal, the first unipotent period $\varphi^{(U_n,1)}\equiv 0$ 
since $U_n$ is the unipotent radical of the maximal parabolic subgroup $P_n$.  
Hence this contribution to $\left(\varphi^{(N_\ell,\Psi_\ell)}\right)^{(U_G,1)}(g)$ vanishes.  Hence
\[
\begin{aligned}
\left(\varphi^{(N_\ell,\Psi_\ell)}\right)^{(U_G,1)}(g) 
&= \sum\limits_{\gamma\in P^1_{n-1,1}(k)\bs \GL_{n}(k)} \int\limits_{U^{\ell}_{n+\ell}(k)\bs U^{\ell}_{n+\ell}(\A)} 
\varphi_{\ell-1}\left(u\hat{\gamma}\delta_\ell g\right) \Psi_{n+\ell}^{-1}(u) \, du \, dx \\
&= \sum\limits_{\gamma\in P^1_{n-1,1}(k)\bs \GL_{n}(k)} \int\limits_{\X(\A)} 
\int\limits_{U^{\ell}_{n+\ell}(k)\bs U^{\ell}_{n+\ell}(\A)} 
\varphi\left(u\hat{\gamma}\lambda \delta_\ell g\right) \Psi^{-1}_{n+\ell}(u) \, du \, dx \, d\lambda.
\end{aligned}
\]

To proceed, they further Fourier expand $\varphi$ or $\varphi_{\ell-1}$ along the ideas of Piatetski-Shapiro and Shalika. 
\begin{enumerate}
\item[$\bullet$] If $f$ is an automorphic form on $\GL_{n}$, to obtain its Fourier expansion, 
first restrict to the mirabolic $P^1_{n-1,1}$ and expand along its unipotent radical 
$U_{n-1}=U(P^1_{n-1,1})\simeq \A^{n-1}$.
\item[$\bullet$] For each such Fourier coefficient, which can now be viewed as a function on 
\[ P^1_{n-2,1}(k)\bs P^1_{n-2,1}(\A), \] 
expand along the unipotent radical of $P^1_{n-2,1}\simeq \A^{n-2}$. 
\item[$\bullet$] Continue inductively until one obtains the sum of the Whittaker function over $Z_n(k)\bs \GL_{n}(k)$. 
\end{enumerate}

To this end, set 
\[
n_{\ell+1}(x)=\bpm I_{n-1} & x \\ & 1\epm^\wedge\in N_{n+\ell}(\A), \quad x\in \A^{n-1}
\]
and let
\[
\phi_{1,\hat{\gamma},g}(x) = 
\int\limits_{U^\ell_{n+\ell}(k)\bs U^\ell_{n+\ell}(\A)} 
\varphi_{\ell-1}\left(un_{\ell+1}(x)\hat{\gamma}\delta_\ell g\right) \Psi^{-1}_{n+\ell}(u) \, du.
\]
Then, exactly as before, when we write the Fourier expansion of $\phi_{1,\hat{\gamma},g}(x)$ in $x$ 
and evaluate at $x=0$ we obtain
\[
\begin{aligned}
\phi_{1,\hat{\gamma},g}(0) 
= \sum_{\gamma'\in P^1_{n-2, 1}(k)\bs \GL_{n-1}(k)} 
\int\limits_{U^{\ell+1}_{n+\ell}(k)\bs U^{\ell+1}_{n+\ell}(\A)} 
&\varphi_{\ell-1}\left(u\hat{\gamma'}\hat{\gamma}\delta_\ell g\right) \Psi^{-1}_{n+\ell}(u) \, du \\ 
&+\left(\varphi_{\ell-1}^{(U_{n-1},1)}\right)^{(N_{\ell+1},\Psi_{\ell+1})}\left(\hat{\gamma}\delta_\ell g\right).
\end{aligned}
\]
If, following \cite{grs}, we denote
\[
P_{n-i, 1,\dots,1}^{1,\dots,1}(k)=\left\{\bpm g & x \\ & z\epm \in \GL_{n}(k) \, : \, z\in Z_i(k)\right\}
\]
then, as above, this gives
\[
\begin{aligned}
\left(\varphi^{(N_\ell,\Psi_\ell)}\right)^{(U_G,1)}(g) 
&= \sum_{\gamma\in P^{1,1}_{n-2,1,1}(k)\bs \GL_{n}(k)} 
\int\limits_{U^{\ell+1}_{n+\ell}(k)\bs U^{\ell+1}_{n+\ell}(\A)} 
\varphi_{\ell-1}\left(u\hat{\gamma}\delta_\ell g\right) 
\Psi_{n+\ell}^{-1}(u) \, du \\  
&+ \sum_{\gamma\in P^1_{n-1}(k)\bs \GL_{n}(k)} 
\left(\varphi_{\ell-1}^{(U_{n-1},1)}\right)^{(N_{\ell+1},\psi_{\ell+1})}(\hat{\gamma} \delta_\ell g).\\
\end{aligned}
\]
 If we, as above, write this in terms of our original $\varphi$, so again interchanging a number of integrals, this becomes
 \[
\begin{aligned}
\left(\varphi^{(N_\ell,\Psi_\ell)}\right)^{(U_G,1)}(g) 
&= \sum_{\gamma\in P^{1,1}_{n-2,1,1}(k)\bs \GL_{n}(k)} 
\int\limits_{\X(\A)} \int\limits_{U^{\ell+1}_{n+\ell}(k)\bs U^{\ell+1}_{n+\ell}(\A)} 
\varphi\left(u\hat{\gamma}\lambda\delta_\ell g\right) \Psi_{n+\ell}^{-1}(u) \, du \, d\lambda \\  
&+ \sum\limits_{\gamma\in P^1_{n-1}(k)\bs \GL_{n}(k)} \int\limits_{\X(\A)} 
\left(\varphi^{(U_{n-1},1)}\right)^{(N_{\ell+1},\psi_{\ell+1})}\left(\hat{\gamma} \lambda \delta_\ell g\right) \, d\lambda.\\
\end{aligned}
\] 
Once again, by the cuspidality of $\varphi$, we see $\varphi^{(U_{n-1},1)}\equiv 0$ and so we have
\[
\begin{aligned}
\left(\varphi^{(N_\ell,\Psi_\ell)}\right)^{(U_G,1)}(g) 
&= \sum\limits_{\gamma\in P^{1,1}_{n-2,1,1}(k)\bs \GL_{n}(k)} 
\int\limits_{\X(\A)} \int\limits_{U^{\ell+1}_{n+\ell}(k)\bs U^{\ell+1}_{n+\ell}(\A)} 
\varphi\left(u\hat{\gamma}\lambda\delta_\ell g\right) \Psi_{n+\ell}^{-1}(u) \, du \, d\lambda .
\end{aligned}
\] 

Continuing in this fashion inductively, we arrive at the statement of Proposition \ref{sum-int-exp}.

%%%%%%%%%%%%%% NEW SECTION %%%%%%%%%%%%%%%%%%%% 
%%%%%%%%%%%%%%%%%%%%%%%%%%%%%%%%%%%%%%%%%%%%%%% 
\section{Euler Product Expansion}\label{sec-euler-exp}
As a consequence of Theorems \ref{thm-BasicIdentityI} and \ref{thm-BasicIdentityII} we can now obtain an Euler product expansion 
for the global integrals therein, which will in turn allow us to relate the integrals to the 
generic Rankin-Selberg $L$-functions for $\GSpin \times \GL$ (with arbitrary rank and including the quasi-split forms).   

Recall that $\A = \A_k$ and we have the embedding $G = G_{n'} \hookrightarrow H = H_{m'}$ in all cases.  
We can summarize the various cases from Sections \ref{sec-integralB} and \ref{sec-integralA} in the table below.  
In each case, we also indicate the Rankin-Selberg $L$-function the integral will produce (see below).

\begin{center} 
\begin{tabular}{|c|c|c|} 
\hline
&&\\
& \mbox{\bf case A: Sec.\ref{sec-integralA} integrals} & {\bf case B: Sec.\ref{sec-integralB} integrals } \\
&&\\
\hline
&&\\
\multirow{4}{*} {\bf odd} & 
$G=\GSpin_{2n} \xhookrightarrow{} H=\GSpin_{2m+1}$ & 
$G=\GSpin_{2n+1} \xhookrightarrow{} H=\GSpin_{2m}$ 
\\  
& $m \ge n, (\ell = m-n)$
& $m > n, (\ell = m-n-1)$ 
\\
& 
$L(s, \pi \times \tau)$ for $H \times \GL_{n}$ &
$L(s, \pi \times \tau)$ for $G \times \GL_{m}$
\\ 
& 
($G$ split, $H$ split)&
($G$ split, $H$ split) 
\\
&&\\
\hline
&&\\
\multirow{4}{*} {\bf even} & 
$G=\GSpin_{2n+1} \xhookrightarrow{} H=\GSpin_{2m}$ & 
$G=\GSpin_{2n} \xhookrightarrow{} H=\GSpin_{2m+1}$ 
\\
& $m > n, (\ell = m-n-1)$ & 
$m \ge n, (\ell = m-n)$ 
\\
& 
$L(s, \pi \times \tau)$ for $H \times \GL_{n}$ &
$L(s, \pi \times \tau)$ for $G \times \GL_{m}$
\\ 
& 
($G$ split, $H$ quasi-split)&
($G$ quasi-split, $H$ split) 
\\
&&\\
\hline
\end{tabular} 
\label{table}
\end{center} 
% 

% % % % % % % % % % % % % % % 
\begin{rem}
When we do the so-called ``unramified computation'' later, following the ideas of Soudry \cite[\S 12]{soudry-mem} 
there will be a certain duality, 
to be made more precise later, between the diagonal entries of 
this table.  The case of $\GSpin_{2a+1} \times \GL_{b}$ will be related to the case of 
$\GSpin_{2b} \times \GL_{a}.$  In the local setting, both will give a Rankin-Selberg $L$-function 
of degree $2 a b$ (in $q^{-s}$).  We will discuss this  in Section \ref{sec:uc}. 
\end{rem}
% % % % % % % % % % % % % % % 

% % % % % % % % % % % % % % % 
\begin{rem}
We note, as it is pointed out in \cite[page 58]{gpsr}, that when $H=\SO_{2m+1},$  Method A with 
$m=n$ gives the $\SO_{2m+1}\times\GL_{m}$ case while Method B when $G=\SO_{2n+1}$ with $m=n+1$ 
gives $\SO_{2n+1} \times \GL_{n+1}.$  
On the other hand, when $H=\SO_{2m},$ Method A with $n=m-1$ gives the $\SO_{2m} \times \GL_{m-1}$ case 
while Method B when $G=\SO_{2n}$ with $m=n$ gives the $\SO_{2n} \times \GL_{n}$ case. These are consistent with our choice of 
names for (case A) and (case B). In \cite[Part B]{gpsr}, Gelbart and Piatetski-Shapiro write down the details of 
the first and the last of these four cases. 
(They also include their Method C that treats the symplectic groups and remark that they could also be applied to 
$\Spin$ groups.)  
\end{rem}
% % % % % % % % % % % % % % % 

Our goal in this section is to factor the right hand sides of the basic identities in Theorems \ref{thm-BasicIdentityI} and \ref{thm-BasicIdentityII} 
as products of local zeta integrals (over all the places of $k$).  Given that the adelic domains of integration factor over the places of $k$, 
we need to show that we may choose the data in the integrands that also factor and that the resulting local zeta integrals converge absolutely 
in some right half plane, not depending on the place of $k$, and that the Euler product converges in that half plane.  Fortunately, the same 
proofs as in the case of special orthogonal groups in \cite{soudry-mem} and \cite{kaplan-thesis} work without the need for much modification; 
the existence of the nontrivial center in the general spin groups makes little difference for the purposes of the results in this section. Therefore, 
we only briefly review the steps below.

Let $\pi = \otimes_v \pi_v$ be an irreducible, cuspidal, globally $\psi$-generic, automorphic representation of 
$G(\A)$ in (case B), resp. of $H(\A)$ in (case A).   
Choose a decomposable $\varphi \in V_\pi$ so that 
the corresponding Whittaker function $W_\varphi = W^\psi_\varphi$ in the Whittaker model $\W(\pi, \psi)$ 
is a product of local Whittaker functions 
\begin{equation*} 
W_\varphi(x) 
= 
\prod\limits_{v} 
W_v(x_v),  
\end{equation*} 
with each $W_v$ a nonzero function in the local Whittaker model $\W(\pi_v,\psi_v)$.  In fact, since 
$\W(\pi, \psi) \not= \{0\}$ by assumption, the map $\varphi \mapsto W_\varphi$ is an isomorphism 
onto $\bigotimes'_v \W(\pi_v,\psi_v)$.

Similarly, let $\tau = \otimes \tau_v$ be an irreducible, cuspidal, automorphic representation of $\GL_{m}(\A)$ in (case B), 
resp. of $\GL_{n}(\A)$ in (case A), and assume that $f_s$ in (\ref{fs-I}), resp. (\ref{fs-II}), 
is a decomposable section.  
When we apply the $\GL$-Whittaker coefficient on $\tau$ with respect to the character $\psi^{-1}$ to it, we can write  
\begin{equation*} 
f_s^{(Z_m,\psi)} (y) 
= 
\prod\limits_{v} 
f_{v,s}(y_v; I_m), \quad \mbox{(case B),}
\end{equation*} 
resp. 
\begin{equation*} 
f_s^{(Z_n,\psi_n)} (y) 
= 
\prod\limits_{v} 
f_{v,s}(y_v; I_n), \quad \mbox{(case A).} 
\end{equation*} 
Here $f_{s,v}$ is a $K$-finite holomorphic section in $\rho_{\tau_v,s}$ taking values in the local Whittaker model 
$\W(\tau_v,\psi^{-1}_v)$ of $\tau_v$.   
For a fixed $y_v$ we denote the corresponding Whittaker function in the Whittaker model of $\tau_v$ by 
$x \mapsto f_{\tau_v,s}(y_v; x).$

We choose a finite set $S$ of places of $k,$ containing all the places at infinity, outside of which 
all data are unramified.  For $v \not\in S,$ we take $W_v = W_v^0$ 
to be the unique Whittaker function such that its value at the identity element is $1$ and take the function 
$x \mapsto f_{\tau_v,s}(I;x)$ be the unique spherical and normalized Whittaker function in the corresponding Whittaker model of $\tau_v.$   
Again, denote this unique section by $f^0_{v,s}.$

%%%%%%%%%%%
\begin{thm}\label{euler} 
Let $\L(\varphi, f_s)$ be as in Theorem \ref{thm-BasicIdentityI}, resp. Theorem \ref{thm-BasicIdentityII}.  
\begin{itemize}
\item[(i)] 
With a choice of data as above, 
for $\Re(s) \gg 0$ we have  
\begin{equation*} 
\L(\varphi, f_s) = \prod_{v} 
\L_v(W_v, f_{v,s}), 
\end{equation*} 
where the local factor is given by 
\begin{flalign}
\label{xi}
&\L_v(W_v, f_{v,s}) 
 =  
 \\ 
& 
\begin{cases} 
\int\limits_{N_G(k_v)Z_G(k_v) \backslash G(k_v)} W_v(g) 
\int\limits_{N_\ell(k_v) \cap \beta^{-1} P_m(k_v) \beta \backslash N_\ell(k_v)} 
f_{v,s} \left( \beta u g; I_m \right) \Psi_{\ell,v}(u)^{-1}  \, du \, dg, 
	& 
\mathrm{(case B),} 
	\\
	\\
\int\limits_{N_G(k_v) Z_G(k_v) \backslash G(k_v)} f_{v,s}(g; I_n) 
\int\limits_{\X(k_v)} 
W_v (\lambda \delta_\ell g) 
 \, d\lambda \, dg, 
	& 
\mathrm{(case A).} 
\end{cases} 
\end{flalign} 
Here, each local factor $\L_v$ converges absolutely in a fixed right half plane independent of the $v$ 
and continues to a meromorphic function on the complex plane. 
\item[(ii)] 
When $v$ is finite, we can choose the data $W_v$ and $f_{v,s}$  such that 
$\L_v(W_v, f_{v,s})$ is identically 1 as a function of $s.$  
\item[(iii)] 
When $v$ is infinite, given $s_0 \in \C$ we can choose the data $W_v$ and $f_{v,s}$ 
such that $\L_v(W_v^{\psi_v}, f_{v,s})$ is holomorphic and nonzero in a neighborhood of $s_0.$  
\end{itemize}
\end{thm}
%%%%%%%%%%%

%%%%%%%%%%%
\begin{proof} 
As we mentioned already the proof follows exactly as in the case of special orthogonal groups and the 
existence of the large center in the case of the general spin groups does not make much of a difference 
for this proof. 

For (i) one starts out by estimating the Whittaker functions $W_v(g)$ in (case B). For the non-archimedean 
$v$ this is done exactly as in \cite[\S 2]{soudry-mem} by estimating the Whittaker function by a ``gauge", 
originally introduced by Jacquet, Piatetski-Shapiro, and Shalika in the case of the general linear groups.  
These estimates result in the convergence of the local integrals for $\Re(s) \ge s_0$, with $s_0 \in \R$ only 
depending on the groups and the embeddings, and not the place $v$ as in \cite[\S 4]{soudry-mem} for 
the non-archimedean $v$.

For the archimedean places $v$ again we follow Sourdy as in \cite[\S 3]{soudry-mem} for the estimates, 
where one also appeals to results of Dixmier and Malliavin \cite{DM}.  The convergence then follows 
as in \cite[\S 5]{soudry-mem} for the archimedean $v$.  
While Soudry focuses on the odd case in \cite{soudry-mem} 
the even case is also covered in \cite{kaplan-thesis}.  The analogous results for (case A) are already 
covered by Ginzburg \cite{ginzburg} and also reviewed by \cite{kaplan-thesis} in the even case.  Since in the 
case of the general spin groups we already divide by the center in the domains of the integrals, no modifications 
in the above proofs are necessary and we can conclude (i).

We should note here that the results we cite above show that region of convergence of the local integrals 
depend only on the representations, and not on the data $W_v$ and $f_{v,s}$.  The dependence on 
the local representations is through their exponents which can be uniformly bounded.  Therefore, 
when $\Re(s)$ is sufficiently large we have convergence that is valid for all $v$. 
Similarly part (ii) proceeds as in \cite[\S 6]{soudry-mem}.

Also, part (iii) follows as in \cite[\S7]{soudry-mem}.  We need to choose $K$-finite data and show that 
the integral admits meromorphic continuation which is continuous in the input data. Here, the argument 
of \cite{soudry-annsci} applies to (case A) and (case B) where we can follow \cite[p. 402]{kaplan4} in the spit case. 
The non-split, quasi-split case for the special orthogonal groups (and for $\GSpin$ groups) has not been yet 
appeared in a published paper and, as in the $\SO$ case in \cite{kaplan4}, we assume it for the $\GSpin$ 
groups.   
\end{proof}
%%%%%%%%%%%

%%%%%%%%%%%%%% NEW SECTION %%%%%%%%%%%%%%%%%%%% 
%%%%%%%%%%%%%%%%%%%%%%%%%%%%%%%%%%%%%%%%%%%%%%% 
\section{The Unramified Computations}\label{sec:uc}
In this section we compute the local zeta integrals for ``unramified data''.  This local 
analysis allows us to relate $\L_v(s, W_v, f_{v,s})$ in Theorem \ref{euler} 
to the local $L$-functions when $v \not\in S,$ $W_v=W^0_v$ and $f_{v,s}=f^0_{v,s}$ and 
is usually referred to as the ``unramfied computation''.

In the split case with $n$ equal, or nearly equal to $m$, the unramfied computation was already 
done for the odd and even special orthogonal groups in \cite[Appendix]{gpsr}.  More precisely, Gelbart, 
Piatetski-Shapiro, and Rallis worked out the cases of $\SO_{2n+1} \times \GL_n$ and $\SO_{2n} \times \GL_{n-1}$ 
in (case A) and the cases of $\SO_{2n+1} \times \GL_{n+1}$ and $\SO_{2n} \times \GL_n$ in (case B), 
even though the emphasis in \cite[Appendix]{gpsr} is on the first and last of these four cases. Their method 
combines the Casselman-Shalika formula with a decomposition of the symmetric algebra of polynomials 
defined on complex matrices of appropriate size on which the (connected components) of the Langlands 
duals of the special orthogonal group and the general linear group act. They then show that the 
local zeta integral with unramified data is a quotient of $L$-functions, by explicitly calculating the integral 
and the local $L$-functions as series in $q^{-s}$, where $q$, as usual, 
denotes the cardinality of the residue field of the local non-Archimedean field.  
For the symmetric algebra decomposition they invoke 
some results of Ton-That \cite{ton-that1, ton-that2}.

The unramified computation in the case of $\SO_{2n+1} \times \GL_m$ with $m \le n$ was then carried 
out by Ginzburg in \cite{ginzburg} as well as the case of $\SO_{2n} \times \GL_m$ with $m \le n-1$, 
where he uses a certain inductive argument to reduced the proof for the more general $m$ to those of 
$m=n$ or $m=n-1$ in the odd and even cases, respectively.

Instead of extending Ginzburg's inductive argument to the general spin groups, we have followed the original 
approach of Gelbart, Piatetski-Shapiro, and Rallis mentioned above in (case A).  This is possible because Ton-That's results 
on the decomposition of the symmetric algebra are fortunately available for $m \le n$, resp. $m \le n-1$, in the 
odd, resp. even, cases (and not just $m=n$ and $m=n-1$ respectively, which is what was used in \cite[Appendix]{gpsr}).

However, in (case B) when the rank of the general linear group is larger than the rank of the $\SO$ groups 
the decomposition of the symmetric algebra becomes too complicated to be helpful. Soudry \cite[\S 12]{soudry-mem} 
then showed how to deal with the case of $\SO_{2n+1} \times \GL_m$ with $m > n$ by relating 
the local Rankin-Selberg $L$-functions in this case to that of $\SO_{2m} \times \GL_n$, where we already 
have the (case A) results.  E. Kaplan has also extended Soudry's method to the case of $\SO_{2n} \times \GL_m$ with 
$m \ge n$ as well as considering the quasi-split forms of $\SO_{2n}$ \cite{kaplan-thesis}. Kaplan also 
suggests modifications of the method in \cite{kaplan2}.  Both Soudry and Kaplan use a certain uniqueness 
result that is fortunately now available for the general spin groups as well thanks to \cite{kll}.  This allows us to apply 
Soudry's ideas in (case B) for the general spin groups.

As mentioned, we consider the unramified computation for the general spin groups, both when 
the rank of the general spin group is larger and when it is smaller than the rank of the general 
linear group, following the above works. We also consider the quasi-split case in the even case, 
following \cite[\S 3.2.1]{kaplan-thesis} where a similar argument is given for the quasi-split even 
special orthogonal groups. 
Given that the method of proof is somewhat different in (case A) and (case B) as we just explained,  
we state the results in two separate theorems below (cf. Theorem \ref{unram-A} and Theorem \ref{unram-B}) 
even though the statements end up being similar.

For this section only, we let $F$ denote a non-Archimedean local field of characteristic zero with 
ring of integers $O_F$ and the cardinality of the residue field $q=p^f.$  
We also fix a uniformizer $\varpi \in O_F.$  The usual $p$-adic absolute value on 
$F$ is denoted by $|\cdot | = |\cdot|_F,$ with $|\varpi| = q^{-1}$.   
Also, let $\psi$ denote an additive character of $F$ which is unramified, i.e., it is trivial on $O_F,$ 
but non-trivial on $\varpi^{-1} O$.

Recall that we have $G=\GSpin_{n'} \hookrightarrow H=\GSpin_{m'}.$ 
In (case A) we have either 
$n'=2n$, $m'=2m+1$,  and $n \le m,$ (odd case) or 
$n'=2n+1$, $m'=2m,$  and $n < m$ (even case). 
In the latter case, $H$ may be quasi-split.   
Let $(\pi,V_\pi)$ denote an unramified globally $\psi$-generic representation 
of $H(F)$ and let $(\tau,V_\tau)$ be an unramified representation of $\GL_{n}(F).$  
(By unramified we mean a representation which has fixed vectors under the action of 
the maximal compact subgroup $H(O_F)$, resp. $\GL_n(O_F)$.  
Such a representation is then induced from an unramified quasi-character of the torus.)  
The representation $\tau$ is determined by its Frobenius-Hecke or Satake parameter $t_\tau$, 
a semi-simple conjugacy class in $\GL_{n}(\C).$ 
Similarly, let $t_\pi$ denote the parameter of $\pi$, a semi-simple conjugacy 
class in the Langlands dual ${}^LH$. 
We may take ${}^LH \cong \GSp_{2m}(\C)$ when $m' = 2m+1$. When $m'=2n$, we 
have ${}^LH \cong \GSO_{2m}(\C) \times \operatorname{Gal}(\overline{F} / F)$ when $H$ 
is split and ${}^LH \cong \GSO_{2m}(\C) \rtimes \operatorname{Gal}(\overline{F} / F)$ when 
$H$ is quasi-split non-split; cf. \S\ref{duals}(B) and (C).

Similarly, in (case B) we let $(\pi,V_\pi)$ denote an unramified globally $\psi$-generic representation of $G(F)$ 
and let $(\tau,V_\tau)$ be an unramified representation of $\GL_{m}(F),$  
with $n'=2n+1,$ $m'=2m,$ and $m>n$ (odd case), or $n'=2n,$ $m'=2m+1,$ and $m \ge n$ (even case).   
Now, $t_\tau$ is a semi-simple 
conjugacy class in $\GL_{m}(\C)$ and $t_\pi$ is a semi-simple conjugacy 
class in ${}^LG$.  
Again, we may take ${}^LG \cong \GSp_{2n}(\C)$ when $n'=2n+1$ and 
we have ${}^LG = \GSO_{2n}(\C) \times \operatorname{Gal}(\overline{F} / F)$ if $n' = 2n$ with $G$ split and 
${}^LG = \GSO_{2n}(\C) \rtimes \operatorname{Gal}(\overline{F} / F)$ if $n' = 2n$ with 
$G$ quasi-split non-split; cf. \S\ref{duals}(B) and (C).

We fix a Haar measure on the additive group $F$ with the volume of 
$O_F$ equal to $1.$  We can then use this measure to fix a left Haar measure on $G(F)$ and $H(F)$ 
in such a way that 
\begin{equation}\label{measure-normalize}
\operatorname{Vol}(u_\alpha(O_F)) = \operatorname{Vol}\left(\left\{u_\alpha(x) : x \in F, |x| \le 1 \right\}\right) = 1,
\end{equation}
for all roots $\alpha$ in $G$ and $H.$ Here, the image of $u_\alpha$ is the root group associated with $\alpha$ 
in the ambient reductive group.  
In particular, we will have $\operatorname{Vol}(K_G) = \operatorname{Vol}(T_G \cap K_G) = 1$ and similarly 
for $H$. Here, $K_G$, resp., $K_H$, denotes the (fixed choice of a) maximal compact in $G(F)$, resp., $H(F)$, 
and $T_G$, reps. $T_H$, denotes the (fixed) maximal torus in $G$, resp., $H$.

We now review some preliminary facts that help us relate the local integrals of Theorem \ref{euler} to 
the local $L$-functions.  In the local setting, the local field $F$ will always be the completion of the 
number field $k$ at a non-Archimedean place $v$ of $k$ and the local representations $\pi$ 
and $\tau$ above will be the component at $v$ of the corresponding global representations bearing 
the same names in Sections \ref{sec-integralB} and \ref{sec-integralA}.

%%%%%%%%%%%%%%%%%%%%%%%%%%%%%%%%%%%%%%%%%%%%%%% 
\subsection{Symmetric Algebra Decompositions}\label{symm-alg-subsubsec} 
We next recall some preliminary facts about the decomposition of a symmetric algebra that is needed 
for our results.  Note that this analysis is only feasible in (case A) with both parities.  Indeed, it is precisely 
the complications with the analysis of the symmetric algebra decompositions in (case B) that  
will force us to use an alternative method in that case, as we will see below.

We refer to \cite{ton-that1,ton-that2} for the details.  While the results of Ton-that are 
proved for the case when the rank of the classical group is larger than or equal to the 
rank of $\GL,$ Gelbart, Piatetski-Shapiro and Rallis only used them in the case of 
equal (or almost equal) ranks in \cite[Appendix]{gpsr}.  
We should add that Ton-that's results are for the non-similitude situation.  
We will use \cite{ton-that1} with $m \ge n$ when $m'=2m+1$ and \cite{ton-that2} with $m > n$ when $m' = 2m$ 
and make the necessary adjustments in the calculations for the similitude situation.

Assume that $m \ge n.$  
Let $E = \C^{n \times 2m}$ denote the space of $n \times 2m$ complex matrices. Then 
$\widehat{\GL}_{n}=\GL_{n}(\C)$ acts on $E$ on the left and 
$\widehat{H} = \GSp_{2m}(\C)$, resp., $\GSO_{2m}(\C)$, acts on $E$ on the right. 
Let $S(E^*)$ denote the symmetric algebra of complex-valued polynomial functions 
on $E.$ It becomes a $\left(\widehat{H} \times \widehat{\GL}_{n}\right)$-module via 
\begin{equation} 
\left((g_1,g_2) \cdot P\right)(X) = P\left( {}^t\!g_2 X g_1 \right), \quad 
X \in E, P \in S(E^*), g_1\in\widehat{H}, g_2\in\GL_{n}(\C).
\end{equation} 

Our goal is to decompose this module in a way that is useful for our setting.  
Recall that the group $\widehat{H} = \GSp_{2m}(\C),$ resp., $\GSO_{2m}(\C)$, is defined as in \cite[\S 2.3]{duke}, i.e., 
it is the connected component of the group 
\begin{equation}\label{gspgso}  
\left\{ g\in\GL_{2m}(\C) : {}^t g J g = \mu(g) J \right\}, 
\end{equation} 
where the $2m\times 2m$ matrix $J$ is defined via  
\begin{equation}\label{jmatrix} 
J =  \left(    \begin{array}{cccccc}  
                                   &   & & & & 1 \\ 
                                      &&&& \adots & \\  
                                   &&& 1&& \\ 
                                   &&-1&&& \\ 
                                   &\adots&&&& \\ 
                                   -1 &&&&&  
                      \end{array}   \right), 
\quad\mbox{resp., }  
J = \left(    \begin{array}{cccccc}  
                                   &   & & & & 1 \\ 
                                      &&&& \adots & \\  
                                   &&& 1&& \\ 
                                   &&1&&& \\ 
                                   &\adots&&&& \\ 
                                   1 &&&&&  
                      \end{array}   \right),    
\end{equation}  
and $\mu$ denotes the similitude character. (Notice that the algebraic group defined in (\ref{gspgso}) 
is connected for the former $J,$ but it has two connected components for the latter $J,$ cf. \cite[\S 2.3]{duke}.)

Next, we introduce two subalgebras $I(E^*)$ and $H(E^*)$ of $S(E^*)$ such that 
\begin{equation} 
S(E^*) \cong I(E^*) \otimes H(E^*) 
\end{equation} 
as $\left(\widehat{H} \times \widehat{\GL}_{n}\right)$-modules. 
We let $I(E^*)$ be the subalgebra of all $\Sp_{2n}(\C)$-, resp., $\SO_{2n}(\C)$-invariant polynomials in $S(E^*).$  
Equivalently, $I(E^*)$ is the algebra of polynomials on the space 
\begin{equation*} 
\left\{ Y = X \, J \,\, {}^tX:  X \in E \right\}, 
\end{equation*} 
where $J$ is as in (\ref{jmatrix}). The action of $\GL_{n}(\C)$ on the space of polynomials of degree $i$ in 
$I(E^*)$ is given by 
\begin{eqnarray}
\nonumber
\sym^i \left( \wedge^2(g_2) \right),& \quad & g_2 \in \GL_{n}(\C), \mbox{ or}\\ 
\sym^i \left( \sym^2(g_2) \right), &\quad & g_2 \in \GL_{n}(\C),  
\end{eqnarray} 
respectively, while the action of $\widehat{H}=\GSp_{2m}(\C),$ resp., $\GSO_{2m}(\C),$ is given simply by  
\begin{equation} 
\mu(g_1)^{i}, \quad g_1 \in \widehat{H}, 
\end{equation}
in both cases.

Similarly, $H(E^*)$ denotes the subspace of $S(E^*)$ 
consisting of $\Sp_{2m}(\C)$-, resp., $\SO_{2m}(\C)$-harmonic polynomials. 
Equivalently, $H(E^*)$ is isomorphic, as $\Sp_{2m}(\C)$-, resp., $\SO_{2m}(\C)$-module, to the symmetric algebra of 
polynomials on the space 
\begin{equation*} 
\left\{ X \in E : X \, J \,\, {}^tX = 0 \right\}.  
\end{equation*} 
Let 
\begin{equation}\label{delta} 
\delta = \left(k_1,k_2,\dots,k_n \right) \in \Z^n 
\end{equation} 
with $k_1 \ge k_2 \ge \dots k_n \ge 0.$ For such a 
``dominant'' $\delta$ we define  
\begin{equation}\label{deltabar}
\bar{\delta} = \left(k_1,\dots,k_n, 0, \dots, 0\right) \in \Z^m. 
\end{equation} 

Let $H(E^*,\delta)$ be the subspace of $H(E^*)$ consisting of polynomials transforming 
under $\GL_{n}(\C)$ according to the irreducible finite-dimensional representation 
\begin{equation*}
\rho_\delta^{\GL_{n}(\C)}(g_2), \quad g_2 \in \GL_{n}(\C), 
\end{equation*}
of highest weight $\delta.$ Then, as an $\widehat{H}$-module, $H(E^*,\delta)$ is equivalent 
to the representation 
\begin{equation*}
\rho_{(\bar{\delta};\tr\delta)}^{\widehat{H}}(g_1), \quad g_1 \in \widehat{H}, 
\end{equation*}
where 
\begin{equation} \label{gbar}
\rho_{(\bar{\delta};\tr\delta)}^{\widehat{H}}(g_1) 
= 
\begin{cases} 
\mu(g_1)^{\tr \delta} \cdot \rho_{\bar{\delta}}^{\Sp_{2m}(\C)}\left(\overline{g}_1\right),& 
\mbox{ if } g_1 \in \GSp_{2m}(\C), \\ 
\\ 
\mu(g_1)^{\tr \delta} \cdot \rho_{\bar{\delta}}^{\SO_{2m}(\C)}\left(\overline{g}_1\right),& 
\mbox{ if } g_1 \in \GSO_{2m}(\C),  
\end{cases} 
\end{equation} 
with 
\begin{equation} \label{g1bar}
\overline{g}_1 = \mu(g_1)^{-1/2} g_1 \in \Sp_{2m}(\C), \mbox{ resp., } \SO_{2m}(\C).
\end{equation} 
We recall here that an irreducible finite-dimensional 
representation of $\widehat{H} = \GSp_{2m}(\C),$ resp., $\GSO_{2m}(\C),$ is given as  
\begin{equation*}
\rho_{\left((k_1,k_2,\dots,k_m), k_0\right)} (g) 
= 
\mu(g)^{k_0} \cdot \rho_{(k_1,k_2,\dots,k_m)} \left(\bar{g}\right), 
\quad g \in \widehat{H},
\end{equation*}
where $(k_1,k_2,\dots,k_m) \in \Z^m$ satisfies $k_1 \ge \dots k_n \ge 0,$ $k_0 \in \Z,$ and  
$\rho_{(k_1,k_2,\dots,k_m)}$ denotes the irreducible finite-dimensional representation of 
$\Sp_{2m}(\C),$ resp., $\SO_{2m}(\C),$ of highest weight $(k_1,\dots,k_n).$ 
It follows that 
\begin{flalign} 
&\tr \sym^r \left(g_1 \otimes g_2 \right) 
 =  \label{sym-alg-decomp} 
 \\ 
& 
\begin{cases} 
\sum\limits_{2i+j=r} 
\tr\sym^i(\wedge^2 g_2) 
\mu(g_1)^{i} 
\sum\limits_{\substack{\tr \delta = j \\ \delta \text{ dominant}} } 
\mu(g_1)^j 
\chi_{\bar{\delta}}^{\Sp_{2m}(\C)}\left(\overline{g}_1\right) 
\chi_\delta^{\widehat{\GL}_{n}}(g_2), &
\mbox{ if } \widehat{H}=\GSp_{2m}(\C), 
	\\
	\\
\sum\limits_{2i+j=r} 
\tr\sym^i(\sym^2 g_2) 
\mu(g_1)^{i} 
\sum\limits_{\substack{\tr \delta = j \\ \delta \text{ dominant}} } 
\mu(g_1)^j 
\chi_{\bar{\delta}}^{\SO_{2m}(\C)}\left(\overline{g}_1\right) 
\chi_\delta^{\widehat{\GL}_{n}}(g_2), & 
\mbox{ if } \widehat{H}=\GSO_{2m}(\C). 
\end{cases} \nonumber 
\end{flalign} 
%

%%%%%%%%%%%%%%%%%%%%%%%%%%%%%%%%%%%%%%%%%%%%%%% 
\subsection{The Casselman-Shalika Formula}\label{cs-subsec}
The Casselman-Shalika formula \cite{cs} evaluates the normalized spherical Whittaker function 
of an unramified representation of a connected reductive group.  
Here, normalized means that we choose the Whittaker function to have the value $1$ at the identity 
as we explained in Section \ref{sec-euler-exp}.  
For a split group one can combine the formula with the Weyl character formula (for the dual group) 
to arrive at the form we state below, as can be found, for example, in \cite[\S 3]{bump}.  When the 
group is a quasi-split general spin group, we modify the formula accordingly.

Consider an irreducible, admissible, unramified representation $\pi$ of the $F$-points of 
a split, connected, reductive group $G$ (over $F$), with a fixed Borel $B = T U_G$, where $T$ is 
a maximal torus and $U_G$ is the unipotent radical of $B$.  Let $\chi$ be an unramified character of $U_G$ 
and let $W^0 \in \cW(\pi, \chi)$ denote a normalized spherical 
Whittaker function. Then, $W^0$ is right $K_G$-invariant and satisfies $W(1) =1$.  Hence, it is completely 
determined by its values on $T(F)=T_G(F)$ and in fact, as we will recall below, by the dominant elements 
$t$ in $T(F)/T(O_F).$  Here, dominant means that $|\alpha(t)| \le 1$ for all simple roots $\alpha.$

Moreover, there is a parametrization of the irreducible representations of the 
complex dual group $\widehat{G}$ by the dominant elements $t$ in $T(F)/T(O_F)$.   
Let $t_\lambda$ be such a dominant element and assume that it corresponds to 
the representation whose character $\chi_\lambda$ has highest weight vector $\lambda.$ 
Then the Casselman-Shalika formula can be stated as 
\begin{equation}\label{cs-formula}
W^0(t_\lambda) = \delta_G(t_\lambda)^{1/2} \chi_\lambda^{\widehat{G}}(t_\pi), 
\end{equation} 
where $\delta_G$ denotes the modulus character of the standard Borel subgroup in $G$  
and $t_\pi$ denotes the semi-simple conjugacy class in the complex dual group parametrizing $\pi$. 
See \cite[Prop. 1]{tamir} for another example of the formulation \eqref{cs-formula}.

We recall the explicit form of the $t_\pi$ for the groups of interest to us.  Let $\pi = \operatorname{Ind}_{B_G}^G (\mu)$, 
where $\mu$ is a character of $T_G(F)$.  When $G = \GL_m$, we have $\mu = \chi_1 \otimes \cdots \otimes \chi_m$, with 
$\chi_i$ unramified quasi-characters of $F^\times$.  Similarly, when $G = \GSpin_{2n+1}$ or $\GSpin_{2n}$ (split), we have 
$\mu = \chi_0 \otimes \chi_1 \otimes \cdots \chi_n$.  Here, $\chi_0$ is 
the pullback of the central character of $\pi$ to $e^*_0(\GL_1(F))$.  Then, a representative 
$t_\pi$ for the semi-simple conjugacy class in the dual group corresponding to $\pi$ is given by 
\begin{equation}\label{tpi}
t_\pi = \begin{cases} 
\operatorname{diag}\left( \chi_1(\varpi), \dots, \chi_m(\varpi)  \right) \in \GL_m(\C) & \mbox{ if } G = \GL_m, \\ 
\operatorname{diag}\left( \chi_1(\varpi), \dots, \chi_n(\varpi), \chi^{-1}_n \chi_0(\varpi),\dots, 
\chi^{-1}_1 \chi_0(\varpi)  \right) \in \GSp_{2n}(\C) & \mbox{ if } G = \GSpin_{2n+1}, \\ 
\operatorname{diag}\left( \chi_1(\varpi), \dots, \chi_n(\varpi), \chi^{-1}_n \chi_0(\varpi),\dots, 
\chi^{-1}_1 \chi_0(\varpi)  \right) \in \GSO_{2n}(\C) & \mbox{ if } G = \GSpin_{2n}.  \\ 
\end{cases} 
\end{equation} 

Next, assume that $G = \GSpin_{2n}^a$ with $a$ a non-square in $F^\times$ is quasi-split, but not split over $F$. 
Let $E = F(\sqrt{a})$ be a quadratic extension of $F$ over which $G$ splits, cf. \S\ref{duals}(C).  Then 
\[ T_G(F) = F^\times \times \left( F^\times \right)^{n-1} \times \GSpin_2^a(F) \] 
with $\GSpin_2^a(F) = \left(\operatorname{Res}_{E/F} \GL_1\right)(F) = E^\times$. 
The unramified character $\mu$ of $T_G(F)$ can be written as 
$\mu =  \chi_0 \otimes \chi_1 \otimes \cdots \chi_{n-1} \otimes \chi\circ\operatorname{Norm}_{E/F}$, and we may 
identify $t_\pi \in {}^LG \cong \GSO_{2n}(\C) \rtimes \operatorname{Gal}(E/F)$ with 
\[ 
t_\pi =  \operatorname{diag}\left( \chi_1(\varpi), \dots, \chi_{n-1}(\varpi), 
\begin{pmatrix} 
\alpha & \beta a \\ 
\beta & \alpha
\end{pmatrix}, 
\chi^{-1}_{n-1} \chi_0(\varpi),\dots, \chi^{-1}_1 \chi_0(\varpi)  \right) \in \GL_{2n}(\C), 
\]  
with $\alpha^2 - a \beta^2 = \chi_0(\varpi)$.

We also let 
\begin{equation} \label{t-prime}  
t'_\pi = \operatorname{diag}\left( \chi_1(\varpi), \dots, \chi_{n-1}(\varpi), \chi^{-1}_{n-1} \chi_0(\varpi), 
\dots, \chi^{-1}_1 \chi_0(\varpi)  \right) \in \GSp_{2(n-1)}(\C).  
\end{equation}  
Then, the analog of the Casselman-Shalika formula for $G$ becomes 
\begin{equation}\label{cs-formula-qsplit}
W^0(t_\lambda) = \delta_G(t_\lambda)^{1/2} \chi_\lambda^{\GSp_{2(n-1)}}(t'_\pi).  
\end{equation} 
%

%%%%%%%%%%%%%%%%%%%%%%%%%%%%%%%%%%%%%%%%%%%%%%% 
%%%%%%%%%%%%%%%%%%%%%%%%%%%%%%%%%%%%%%%%%%%%%%% 
\subsection{Local Identity} \label{loc-id-sec}
We now state and prove the main results of this section: the computation of the integrals with 
unramified data.

Let $\omega$ denote the character of $F^\times$ such that 
$\pi(e_0^*(\lambda)) = \omega(\lambda) \operatorname{Id}_{V_\pi}$ 
with $e_0^*$ is as in Section \ref{sec-pre}.  
For 
$s \in \C$ consider the representation 
$\tau_s = \tau  |\det|^{s} \otimes \omega^{-1}$ 
of $M(F) \cong \GL_{n}(F) \times \GL_{1}(F)$ in (case A), resp., of $M(F) \cong \GL_{m}(F) \times \GL_{1}(F)$ 
in (case B), and define $\rho_{s},$ a representation of $G(F),$ resp., $H(F),$ similar to the global 
setting in (\ref{fs-II}), resp., (\ref{fs-I}), as follows: 
\begin{equation}\label{induced-local} 
\rho_{s} = 
\begin{cases} 
\operatorname{Ind}_{P_{n}(F)}^{G(F)} \left( \tau_s \right),  
& \mbox{(case A),}
\\ 
\\
\operatorname{Ind}_{P_{m}(F)}^{H(F)} \left( \tau_s \right),  
& \mbox{(case B).}
\end{cases}
\end{equation}

Let $W^0_\pi$ be the unique spherical and normalized Whittaker function in the $\psi$-Whittaker model of $\pi$.  
Also, let $f^0_s$ be unramified and normalized so that the vector-valued function $b \to f^0_s(e, b)$, taking 
values in the $\psi^{-1}$-Whittaker model of $\tau$, is the normalized Whittaker function of the general linear group. 
Here, $e$ denotes the identity element of $G(F)$ or $H(F)$ as appropriate.

Similar to the global situation, for $\Re(s) \gg 0,$ we set  
\begin{flalign} 
\nonumber 
\label{loc-zeta} 
&\xi(W^0_\pi, f^0_{s}) 
 =  
  \\ 
& 
\begin{cases} 
\int\limits_{N_G(F) Z_G(F) \backslash G(F)} f^0_s(g; I_n) 
\int\limits_{\X(F)} 
W^0_\pi (\lambda \delta_\ell g) 
 \, d\lambda \, dg, 
	& 
\mathrm{(case A),} 
	\\
	\\
\int\limits_{N_G(F)Z_G(F) \backslash G(F)} W^0_\pi(g) 
\int\limits_{N_\ell(F) \cap \beta^{-1} P_m(F) \beta \backslash N_\ell(F)} 
 f^0_s \left( \beta u g; I_m \right) \Psi_{\ell}(u)^{-1}  \, du \, dg, 
	& 
\mathrm{(case B),} 
\end{cases} 
\end{flalign} 
where 
$\beta$ is as in Theorem \ref{thm-BasicIdentityI} and $\delta_\ell$ is as in Theorem \ref{thm-BasicIdentityII}.

As we mentioned before, there is a significant difference in the way unramified computation 
is carried out in (case A) versus (case B), due to the feasibility of the symmetric algebra 
decompositions. Therefore, we consider these two cases separately even though the statements of the results 
are similar. We first consider (case A).

%%%%%%%%%%%%%%%%%%%
\begin{thm} \label{unram-A}
Let $\pi$ be an irreducible, admissible, unramified, globally $\psi$-generic representation of $H(F)$ where 
$H = \GSpin_{2m+1}$ referred to as the (odd case), or $H = \GSpin_{2m}$, possibly 
non-split quasi-split, referred to as the (even case).  Assume that $n \le m$ in the (odd case),  
or $n < m$ in the (even case). Let $\tau$ be an irreducible, admissible, unramified, globally $\psi^{-1}$-generic 
representation of $\GL_n(F)$ as above.  Choose $W^0_\pi$ and $f^0_{s}$ as before. 
With the Haar measures normalized as in (\ref{measure-normalize}), we have 
\begin{equation}
\xi(W^0_\pi, f^0_{s}) 
 =  
\begin{cases} 
\displaystyle
\frac{L(s,\pi \times \tau)}{L\left(2s, \tau, \wedge^2 \otimes \omega\right)}, & 
\mbox{ (case A), odd,}
\\
&\\
\displaystyle
\frac{L(s,\pi \times \tau)}{L\left(2s, \tau, \sym^2\otimes\omega\right)}, & 
\mbox{ (case A), even, split,}
\\
&\\
\displaystyle
\frac{L(s,\pi \times \tau)}{L\left(2s, \tau, \wedge^2\otimes\omega\right)}, & 
\mbox{ (case A), even, quasi-split.}
\end{cases} 
\end{equation}
(Refer to the table in Section \ref{sec-euler-exp} for the details of the cases.)
\end{thm}
%%%%%%%%%%%%%%%%%%%

%%%%%%%%%%%%%%%%%%%
\begin{proof} 
Recall that we have the embeddings 
\begin{eqnarray*} 
G=\GSpin_{2n}  & \hookrightarrow & H = \GSpin_{2m+1}, \quad n \le m,   \mbox{(case A--odd)}, \\ 
G=\GSpin_{2n+1}  & \hookrightarrow &  H = \begin{cases} 
\GSpin_{2m}, & \mbox{ split},  \\ \GSpin^{a}_{2m}, & \mbox{quasi-split non-split},  
\end{cases} 
\quad n < m,   \mbox{(case A--even)}.  
\end{eqnarray*} 
The above embeddings induce embeddings at the level of $F$-points.  By the Iwasawa decomposition, we have 
\begin{equation}\label{iwasawa} 
G(F) = N_G(F) T_G(F) K_G = N_G(F) Z_G(F) T_1(F) K_G(F),  
\end{equation} 
where $K_G$ is the maximal compact subgroup of $G(F)$, $Z_G(F)$ denotes the connected component of the center of $G(F)$, and   
\begin{equation} \label{t-element}
T_1 = \left\{ t = e^*_1(t_1) \cdots e^*_n(t_n) : t_i \in F^\times \right\}.
\end{equation} 
Below we will also employ $t$ to denote the image of $t$ under the embedding of $G(F)$ into $H(F)$ 
as well as the element 
$\begin{pmatrix} 
t_1 &&& \\ & t_2 && \\ && \ddots & \\
&&& t_n \end{pmatrix} \in \GL_{n}(F),$ 
i.e., we are using $t$ to denote the related elements in $\GL_{n}(F)$, $G(F)$ and $H(F)$.  
Given these identifications, we may write 
\begin{equation} 
f^0_s(t;I_n) = \delta_{P_n}^{1/2}(t)  \mid \det t \mid^{s} W^0_\tau(t)  
\end{equation} 
where 
$W^0_\tau(t)$ is the normalized spherical Whittaker 
function in the $\psi^{-1}$-Whittaker model of $\tau$ (a representation of $\GL_n(F)$).

The integrand on the right hand side of \eqref{loc-zeta} is invariant under multiplying $g$ on the right by an element 
in $K_G$ and a central element on the left.  
Given our normalization of the Haar measures, this means that the integral reduces to 
\begin{equation}
\int\limits_{T_1(F)} 
W^0_{\tau}(t) 
\mid t_1 t_2 \dots t_n \mid^{s+u} \, \delta_G^{-1}(t)
\int\limits_{\X(F)} 
W^0_\pi (\lambda \delta_\ell t) 
 \, d\lambda \, dt, 
\end{equation}
where 
\begin{equation}\label{shift-u}
u = 
\begin{cases} 
\frac{n-2}{2},
& 
\mbox{odd case,}
\\ 
\\
\frac{n-1}{2},
& 
\mbox{even case,}
\end{cases}
\end{equation}
and 
$\delta_G$ 
denotes the modulus function of the Borel subgroup of $G(F)$ (restricted to $T_1$).

Next, we dispose of the integration over $\X(F)$.  
Here, we can argue as in 
\cite[p.~176]{ginzburg} or \cite[p.~98]{soudry-mem}.  For each $t$ we have $W^0_\pi(\lambda \delta_\ell t) = 0$ 
unless $\lambda \in \X(O_F) \subset K_H$.  This follows from the facts that $\delta_\ell \in K_H$, 
that $W^0_\pi$ is right $K_H$ invariant, and that $\delta_\ell t$ normalizes $\X(F)$, leading, by \eqref{X}, to a change of 
variables 
$d\lambda \to |\det t|^{-\ell} d\lambda$.  
where, as in Section \ref{sec-integralA}, 
\begin{equation}\label{shift-ell} 
\ell = 
\begin{cases} 
m-n,
& 
\mbox{odd case,}
\\ 
\\
m-n-1,
& 
\mbox{even case.}
\end{cases}
\end{equation}

Therefore, our integral reduces to 
\begin{equation}\label{zetaWW} 
\zeta(s, W^0_\pi,W^0_\tau)  
= \int\limits_{T_1(F)}  W^0_\pi(t) 
W^0_{\tau}
(t) 
|t_1 t_2 \dots t_n|^{s + u - \ell}
\, \delta_G^{-1}(t)
 \, dt.    
\end{equation}
We are reduced to proving that 
\begin{equation} \label{l-zeta-odd}
L(2s,\tau, \wedge^{2}\otimes\omega) 
\, \cdot 
\zeta(s, W^0_\pi,W^0_\tau) 
= L(s,\pi \times \tau) 
\end{equation}
in the odd case, 
\begin{equation}\label{l-zeta-even-split}
L(2s,\tau, \sym^{2}\otimes\omega) 
\cdot 
\zeta(s, W^0_\pi,W^0_\tau) 
= L(s,\pi \times \tau)   
\end{equation}
in the even split case, and 
\begin{equation}\label{l-zeta-even-qsplit}
L(2s,\tau, \wedge^{2}\otimes\omega) 
\cdot 
\zeta(s, W^0_\pi,W^0_\tau) 
= L(s,\pi \times \tau)   
\end{equation}
in the even quasi-split case.

We prove the two sides of \eqref{l-zeta-odd}, resp., \eqref{l-zeta-even-split} and \eqref{l-zeta-even-qsplit}, 
have equal coefficients when expanded as power series in $q^{-s}.$  
We start by expanding $\zeta(s,W^0_\pi,W^0_\tau)$.

By \cite[Lemma 5.1]{cs} we have $W^0_\pi(t) \equiv 0$ unless 
\begin{equation}\label{h-roots}
\left| \alpha(t)\right| \le 1
\end{equation} 
for all simple roots $\alpha$ of $H$, and 
$W^0_{\tau}(t) \equiv 0$ 
unless 
\begin{equation}\label{g-roots}
\left| \alpha(t) \right| \le 1 
\end{equation} 
for all simple 
roots $\alpha$ of $G$ appearing in the Levi $M_n$. Therefore, \eqref{h-roots} implies that 
\begin{equation}\label{h-ineq} 
\operatorname{ord}(t_1) \ge \cdots \ge \operatorname{ord}(t_n) \ge 0 
\end{equation}
while \eqref{g-roots} implies that 
\begin{equation} \label{g-ineq}
\operatorname{ord}(t_1) \ge \cdots \ge \operatorname{ord}(t_n). 
\end{equation}
Notice the crucial fact that the last inequality in \eqref{h-ineq} holds because 
of the structure of the root system of $H$ (both odd and even case),  
and the fact that $n \le m$ in the odd case and $n < m$ in the even case 
(and that $t_{n+1}=\cdots=t_m=1$, i.e., they do not appear).

For $\delta$ as in \eqref{delta} and $\bar{\delta}$ as in \eqref{deltabar}, set 
\begin{equation} 
\varpi^\delta = e_1^*(\varpi^{k_1}) \cdots e_n^*(\varpi^{k_n}) \in T_1(F) \subset G(F) 
\end{equation}  
and write $\varpi^{\bar{\delta}}$ for its image in the maximal torus of $H(F)$ under the embedding. 
We will also use $\varpi^\delta$ to denote the diagonal element in $\GL_{n}(F).$ 
Also, write 
\begin{equation}
\tr \delta = \tr\bar{\delta} = k_1 + k_2 + \cdots + k_m. 
\end{equation} 
Our integral then reduces to  
\begin{equation} \label{zeta-q} 
\zeta(s, W^0_\pi,W^0_\tau) =  
\sum_{\overset{\delta=(k_1,\cdots,k_n)}{\rm{dominant}}}
W^0_\pi(\varpi^{\bar{\delta}}) 
W^0_\tau(\varpi^\delta)
\delta_G^{-1}(\varpi^\delta) 
q^{- \left(s+u-\ell\right) 
\tr \delta},   
\end{equation} 
with $u$ as in \eqref{shift-u} and $\ell$ as in \eqref{shift-ell}. 
Using the Casselman-Shalika formulas \eqref{cs-formula} in the split case 
and \eqref{cs-formula-qsplit} in the quasi-split case, 
we have 
\begin{equation} 
W^0_\tau
(\varpi^\delta)= \delta_{\GL_{n}}^{1/2}
(\varpi^\delta)\chi_\delta^{\widehat{\GL}_{n}}(t_\tau) 
\end{equation}
and 
\begin{equation} \label{Wpi}
W^0_\pi(\varpi^{\bar{\delta}}) = \delta_{H}^{1/2}(\varpi^{\bar{\delta}}) \cdot 
\begin{cases} 
\chi_{(\bar{\delta};\tr\delta)}^{ \GSp_{2m} }(t_\pi), & \mbox{ odd, split} \\ \\ 
\chi_{(\bar{\delta};\tr\delta)}^{ \GSO_{2m} } (t_\pi), & \mbox{ even, split} \\ \\ 
\chi_{(\bar{\delta};\tr\delta)}^{ \GSp_{2(m-1)} }(t'_\pi), & \mbox{ even, quasi-split} 
\end{cases}
\end{equation}
with $t'_\pi$ as in \eqref{t-prime}. 
Also, by \eqref{g1bar}, the right hand sides of \eqref{Wpi} are equal to  
\begin{equation} 
\delta_{H}^{1/2}(\varpi^{\bar{\delta}}) \cdot 
\begin{cases}
\mu(t_\pi)^{\tr\delta} 
\chi_{\bar{\delta}}^{\Sp_{2m}}(\bar{t}_\pi),  
& 
\mbox{ odd, split,}
\\ 
\\  
\mu(t_\pi)^{\tr\delta} 
\chi_{\bar{\delta}}^{\SO_{2m}}(\bar{t}_\pi),  
& 
\mbox{ even, split,}
\\ 
\\
\mu(t'_\pi)^{\tr\delta} 
\chi_{\bar{\delta}}^{\Sp_{2(m-1)}}(\bar{t}'_\pi),  
& 
\mbox{ even, quasi-split.}
\end{cases}
\end{equation} 
(Note that by \eqref{t-prime} we know that $\mu(t_\pi) = \mu(t'_\pi)$ in the even, quasi-split case.) 

For $t \in T_1(F)$ as in \eqref{t-element} (where $t_{n+1}=\cdots=t_m=1$) 
we see from the root data that 
\begin{equation*} 
\delta_G(t) = 
\begin{cases}
\left| t_1^{2n-2} t_2^{2n-4} \cdots t_n^0 \right|,  
& 
\mbox{ odd case,} 
\\ 
\\
\left| t_1^{2n-1} t_2^{2n-3} \cdots t_n^{1} \right|,  
& 
\mbox{ even case,} 
\end{cases}
\end{equation*}
\begin{equation*} 
\delta_H(t) = 
\begin{cases}
\left| t_1^{2m-1} t_2^{2m-3} \cdots t_n^{2m-2n+1} \right|, 
& 
\mbox{ odd case,} 
\\ 
\\
\left| t_1^{2m-2} t_2^{2m-4} \cdots t_n^{2m-2n} \right|, 
& 
\mbox{ even case,} 
\end{cases}
\end{equation*}
and 
\begin{equation*} 
\delta_{\GL_{n}}
(t) = \left| t_1^{n-1} t_2^{n-3} \cdots t_n^{1-n} \right|. 
\end{equation*}
Therefore,  for $t$ as in \eqref{t-element} 
\begin{equation}
\delta_G(t) = \delta_H^{1/2}(t) 
\cdot 
\delta_{\GL_{n}}^{1/2}
(t) \cdot 
\begin{cases} 
\left| t_1 \cdots t_n \right|^{\frac{3n-2m-2}{2}}, 
& 
\mbox{ odd case,} 
\\ 
\\ 
\left| t_1 \cdots t_n \right|^{\frac{3n-2m+1}{2}}, 
& 
\mbox{even case.} 
\end{cases}
\end{equation}

Substituting in (\ref{zeta-q}) we get 
\begin{equation}\label{zeta-exp} 
\zeta(s, W^0_\pi, W^0_\tau) =  
\begin{cases} 
\sum\limits_{\overset{\delta=(k_1,\cdots,k_n)}{\rm{dominant}}}
\mu(t_\pi)^{\tr\delta}
\chi_{\bar{\delta}}^{\Sp_{2m}(\C)}(\bar{t}_\pi) 
\chi_\delta^{\widehat{\GL}_{n}}(t_\tau)
q^{-s \tr\delta}, 
& 
\mbox{odd, split,} 
\\ 
\\ 
\sum\limits_{\overset{\delta=(k_1,\cdots,k_n)}{\rm{dominant}}}
\mu(t_\pi)^{\tr\delta}
\chi_{\bar{\delta}}^{\SO_{2m}(\C)}(\bar{t}_\pi) 
\chi_\delta^{\widehat{\GL}_{n}}(t_\tau)
q^{-s \tr\delta}, 
& 
\mbox{even, split,} 
\\ 
\\ 
\sum\limits_{\overset{\delta=(k_1,\cdots,k_n)}{\rm{dominant}}}
\mu(t'_\pi)^{\tr\delta}
\chi_{\bar{\delta}}^{\Sp_{2(m-1)}(\C)}(\bar{t}'_\pi) 
\chi_\delta^{\widehat{\GL}_{n}}(t_\tau)
q^{-s \tr\delta}, 
& 
\mbox{even, quasi-split.} 
\end{cases}
\end{equation}

Next, recall the well-known identity 
\begin{equation}
\det(I-AX)^{-1} = \sum\limits_{r=0}^{\infty} \tr \left( \sym^r(A) \right) X^r, 
\end{equation}
where $A$ is an arbitrary square complex matrix and $X$ is a sufficiently small complex variable. 
Applying this identity to $A = \left( \omega(\varpi)\wedge^{2}(t_\tau)  \right)$, 
resp., $A = \left( \omega(\varpi)\sym^{2}(t_\tau)  \right)$, 
we obtain 
\begin{eqnarray}\label{L-exp-odd}
L(2s,\tau,R\otimes\omega) 
&=& 
\det \left( I - \omega(\varpi) (R t_\tau)  q^{-2s} \right)^{-1}
\\
\nonumber
&=& 
\sum\limits_{i=0}^\infty \tr \left( \sym^i \left(\omega(\varpi) R t_\tau \right) \right) q^{-2 i s} 
\\
\nonumber
&=& 
\sum\limits_{i=0}^\infty \omega(\varpi)^{i} \tr \sym^i \left( R t_\tau \right)
q^{-2 i s}  
\end{eqnarray}
where 
\begin{equation} \label{R} 
R= \begin{cases} 
\wedge^2& \mbox{odd case,} \\
\sym^2 & \mbox{even, split case,} \\  
\wedge^2& \mbox{even, quasi-split case.} \\
\end{cases} 
\end{equation} 

Multiplying (\ref{zeta-exp}) by \eqref{L-exp-odd}, we see that the 
left hand side of \eqref{l-zeta-odd}, resp., \eqref{l-zeta-even-split}, \eqref{l-zeta-even-qsplit}, is equal to 
\begin{equation}\label{LHS-odd}
\sum\limits_{r=0}^\infty 
\left[ 
\sum\limits_{2i+j=r} \tr \sym^i \left( \wedge^2 t_\tau \right) \omega(\varpi)^{i}  
\sum\limits_{\overset{\tr\delta = j}{\delta \, \rm{dominant}}} 
\mu(t_\pi)^j 
\chi_{\bar{\delta}}^{\Sp_{2m}(\C)}(\bar{t}_\pi) 
\chi_\delta^{\widehat{\GL}_{n}}(t_\tau) 
\right]
q^{-s r}  
\end{equation}
in the odd case, to 
\begin{equation}\label{LHS-even-split}
\sum\limits_{r=0}^\infty 
\left[ 
\sum\limits_{2i+j=r} \tr \sym^i \left( \sym^2 t_\tau \right) \omega(\varpi)^{i}  
\sum\limits_{\overset{\tr\delta = j}{\delta \, \rm{dominant}}} 
\mu(t_\pi)^j 
\chi_{\bar{\delta}}^{\SO_{2m}(\C)}(\bar{t}_\pi) 
\chi_\delta^{\widehat{\GL}_{n}}(t_\tau) 
\right]
q^{-s r}  
\end{equation}
in the even, split case, and to 
\begin{equation}\label{LHS-even-qsplit}
\sum\limits_{r=0}^\infty 
\left[ 
\sum\limits_{2i+j=r} \tr \sym^i \left( \wedge^2 t_\tau \right) \omega(\varpi)^{i}  
\sum\limits_{\overset{\tr\delta = j}{\delta \, \rm{dominant}}} 
\mu(t_\pi)^j 
\chi_{\bar{\delta}}^{\Sp_{2(m-1)}(\C)}(\bar{t}'_\pi) 
\chi_\delta^{\widehat{\GL}_{n}}(t_\tau) 
\right]
q^{-s r}    
\end{equation}
in the even, quasi-split case. 

On the other hand, the right hand sides of \eqref{l-zeta-odd}, 
resp., \eqref{l-zeta-even-split}, \eqref{l-zeta-even-split}, are 
\begin{equation}
L(s, \pi \times \tau) = \det\left( I - (t_\pi \otimes t_\tau) q^{-s}\right)^{-1}
= 
\sum\limits_{r=0}^\infty \tr \sym^r \left(t_\pi\otimes t_\tau \right) q^{-sr}
\end{equation}
and it is enough to verify that $\tr \sym^r\left(t_\pi\otimes t_\tau\right)$ 
is equal to the expression in brackets in \eqref{LHS-odd} in the odd case, 
in \eqref{LHS-even-split} in the even split case, and in 
\eqref{LHS-even-qsplit} in the even, quasi-split case. 
We do this with the help of our earlier discussion of the symmetric algebra decompositions in 
Section \ref{symm-alg-subsubsec}. 
Noting that $\mu(t_\pi) = \omega_\pi(\varpi) = \omega(\varpi)$, 
the equation \eqref{sym-alg-decomp} finishes the proof. 
\end{proof}

Next we consider (case B).

%%%%%%%%%%%%%%%%%%%
\begin{thm} \label{unram-B}
Let $\pi$ be an irreducible, admissible, unramified, globally $\psi$-generic representation of $G(F)$ where 
$G = \GSpin_{2n+1}$ referred to as the (odd case), or $G = \GSpin_{2n}$, split or  
non-split quasi-split, referred to as the (even case).  
Assume that $m > n$ in the odd case or $m \ge n$ in the even case. 
Let $\tau$ be an irreducible, admissible, unramified, globally $\psi^{-1}$-generic representation of $GL_{m}(F)$ as above.  
Also, choose $W^0$ and $f^0_{s}$ as before. 
With the Haar measures normalized as in (\ref{measure-normalize}), we have 
\begin{equation}
\xi(W^0_\pi, f^0_{s}) 
 =  
\begin{cases} 
\displaystyle
\frac{L(s,\pi \times \tau)} { L\left(2s, \tau, \wedge^2 \otimes \omega\right)}, & 
\mbox{ (case B), odd,}
\\
&\\
\displaystyle
\frac{L(s,\pi \times \tau)}{L\left(2s, \tau, \sym^2\otimes\omega\right)}, & 
\mbox{ (case B), even, split.}
\\
&\\
\displaystyle
\frac{L(s,\pi \times \tau)}{L\left(2(s), \tau, \wedge^2\otimes\omega\right)}, & 
\mbox{ (case B), even, quasi-split.}
\end{cases} 
\end{equation}
Here, $\xi(W^0_\pi, f^0_{s})$ is as in (case B) of \eqref{loc-zeta}.  
(Again refer to the table in Section \ref{sec-euler-exp} for the details of the cases.)
\end{thm}

\begin{proof}
The proof in the non-split quasi-split case is a little more involved than the split case so 
we present them separately, indicating the extra issues in that case. 
In the split cases, we follow Soudry's technique for the odd special orthogonal groups \cite[\S 12]{soudry-mem} 
and its adaptation to the even special orthogonal case by Kaplan \cite[\S 3.2.3]{kaplan-thesis}.

%%%
{\it The split cases.}  We first deal with the odd case. The even split case will be similar, as we explain below. 
Since $\pi$ is in general a quotient of a full parabolically induced representation from 
an unramified representation of the Siegel Levi, we may assume that 
\[ \pi = \operatorname{Ind}_{\bar{P}_n(F)}^{G(F)}\left( \sigma \otimes \omega \right) \] 
with $\sigma$ an unramified representation of $\GL_n(F)$ and $\omega = \omega_\pi$ the central character of $\pi.$ 
Here, $\bar{P}_n = M_n \bar{N}_n$ denotes the opposite of the Siegel parabolic $P_n$ (cf. \ref{subsec:siegel}).

Let $\phi$ be a function in $V_\pi$ and assume that for $g \in G(F)$, $\phi(g)$ takes values in $\W(\sigma\otimes\omega,\psi)$.  
As a function of $g$ it is smooth and 
\[ \phi(a_0 \bar{y} g) = | \det a_0 |^{-n/2} \phi(g) \]  
for $\bar{y} \in \bar{N}_n(F)$ and $a_0 \in \GL_n(F)$ considered as the factor in the Siegel Levi of $G(F)$.  Consider 
\begin{equation}\label{Wg}
W_\phi(g) = \int\limits_{N_n(F)} \phi(yg) \psi_n^{-1}(y) \, dy,  
\end{equation} 
with $\psi_n$ obtained from $\psi$ as in Section \ref{sec-integralB} (with everything local now). 
We formally have 
\[ 
W_\phi(ug) = \psi_n(u) W_\phi(g), \quad u \in N_n(F), g \in G(F). 
\] 
However, the integral \eqref{Wg} may not converge absolutely. To remedy this problem, 
we replace 
$\sigma$ by 
\[ 
\sigma_\zeta = \sigma \otimes |\det\cdot|^{-\zeta} 
\] 
for $\Re(\zeta) \gg 0$.  
Replacing $\sigma$ by $\sigma_\zeta$ and taking a holomorphic section $\phi_\zeta$ instead of $\phi$, we 
see just as in \cite{soudry-mem} that the integral defining $W_{\phi_\zeta}$ converges absolutely 
for $\Re(\zeta)$ large enough and has a continuation to a holomorphic function on the whole plane. This is seen 
exactly as in the case of special orthogonal groups by noting that the integral \eqref{Wg} always has a principal 
value and if $\phi_\zeta$ is a standard section, this principal value is a polynomial in $q^{-\zeta}$. Thus, for 
$\Re(\zeta)$ large enough, $W_{\phi_\zeta}(g)$ is a polynomial in $q^{-\zeta}$ which provides the holomorphic 
continuation (cf. \cite[\S 12]{soudry-mem}).

Choose the vector $\phi^0_\zeta = \phi^0_{\sigma,\zeta}$ that gives the normalized unramified Whittaker 
function $W^0_{\sigma} \in \W(\sigma,\psi)$.  
It follows from the Casselman-Shalika formula \cite{cs} (see also \cite[Remark 3.5.14]{shahidi-book}) that 
\[ W^0_\pi = W_{\phi^0_{\sigma,\zeta}}(I)^{-1} W_{\phi^0_{\sigma,\zeta}} = 
L(1+2\zeta, \hat{\sigma}, \sym^2\otimes\omega) W_{\phi^0_{\sigma,\zeta}}. \]

The proof now proceeds the same way as in \cite[\S 12]{soudry-mem} and we  
briefly review the steps indicating the points of difference with our case, which include 
 the appearance of the twisted versions of the symmetric and exterior square $L$-functions.  
 We should also mention that the local $\gamma$- and $\epsilon$-factors 
enter in the steps, which we did not define earlier. However, these factors, which are defined via 
applying intertwining operators to $\rho_s$, are defined for the $\GSpin$ groups 
in \cite{kll}.

Using the above expression for $W_\pi^0$ and proceeding as in \cite[\S 12]{soudry-mem} we see that 
$\xi(W^0, f^0_{s})$ is equal to 
\[ 
\frac{L(1+2\zeta, \hat{\sigma}, \sym^2\otimes\omega)}
{\gamma{\left(s-\zeta, \sigma \times \tau, \psi^{-1}\right)}  L{\left(2s, \tau, \wedge^2\otimes\omega\right)} }
\] 
times a double integral.  This double integral can be manipulated just as in \cite[pp. 97--98]{soudry-mem} 
as the presence of the center in the $\GSpin$ case does not disturb those arguments. As a result, 
the double integral is seen to represent the Rankin-Selberg convolution for $\GSpin(2m) \times \GL(n)$ (as in (case A) 
we considered earlier), which is equal to 
\[ 
\frac{L(\zeta_1+(s - 1/2), \tau\times\hat{\sigma}\omega) L(\zeta_1-(s - 1/2), \hat{\tau}\times\hat{\sigma})} 
{L(2\zeta_1, \hat{\sigma}, \sym^2\otimes\omega)} \cdot 
\] 
with $\zeta_1 = \zeta + \frac{1}{2}$.  As we pointed out above, the only new phenomenon is the appearance 
of the central character $\omega$ and the twisted symmetric square $L$-function in the above expression. 
Therefore, 
\[ 
\xi(W^0, f^0_{s}) = 
\frac{L(s-\zeta, \sigma\times\tau) L(s+\zeta, \hat{\sigma}\omega \times\tau)} 
{L{\left(2s, \tau, \wedge^2\otimes\omega\right)}} \, \cdot  
\] 
Finally, note that the above equalities are all in the sense of equality of rational functions in $q^{-s}$ and $q^{-\zeta}$ 
in their region of convergence and have analytic continuation to $\zeta = 0$. As such, since they are defined for $\zeta = 0$, 
we may indeed substitute $\zeta = 0$. Consequently we obtain 
\[ \xi(W^0, f^0_{s}) = 
\frac{L(s, \sigma\times\tau) L(s, \hat{\sigma}\omega \times\tau)} 
{L{\left(2s, \tau, \wedge^2\otimes\omega\right)}}   
=\frac 
{L(s, \pi\times\tau)} 
{L{\left(2s, \tau, \wedge^2\otimes\omega\right)}} \cdot 
\] 
This finishes the proof in the odd case. 
The even case proceeds in a similar way.  The steps, for the even special orthogonal groups, are detailed in 
\cite[\S 3.2.3]{kaplan-thesis}.  As above, the same steps go through in the split general spin groups, except for 
the fact that again $\omega$ appears and now twisted symmetric square $L$-function shows up. Again, similar 
to the above, we will have 
\[ 
\xi(W^0, f^0_{s})  
=\frac 
{L(s, \pi\times\tau)} 
{L{\left(2s, \tau, \sym^2\otimes\omega\right)}} \cdot 
\]

%%%
{\it The non-split quasi-split case.}  It only remains to consider the case where $G$ is non-split, quasi-split, 
i.e., when $G = \GSpin_{2n}^* = \GSpin_{2n}^a$, as in Section \ref{sec:qsplit}, 
and $H = \GSpin_{2m+1}$ with $m \ge n$. (Recall that a quasi-split $\GSpin_{2n+1}$ is already split, so covered 
above.) 
In this case, the technique described above does not quite extend in a straightforward way due to ``the presence of 
the compact modulo center $\GSpin_{2}^*$ in the middle'' so we need to make some modifications in the proof.  
Here, we can follow the similar proof for the case of non-split quasi-split even special orthogonal groups 
due to E. Kaplan in his thesis \cite[\S 3.2.3 and \S 7.2.3]{kaplan-thesis} and in \cite{kaplan3}.  
Essentially, we may assume the representation $\pi$ of $\GSpin_{2n}^*$ is (a quotient of) an induced 
representation from $\GL_{n-1} \times \GSpin_2^*$ and Kaplan's arguments go through because the 
presence of the center in the general spin groups does not impact those arguments. 
Along the way one also needs a 
local multiplicity at most one result \cite{AGRS,MW2}, which fortunately for the case of 
$\GSpin_{2n} \times \GL_m$ with $n \le m$ is now available 
due to the work of E. Kaplan, J.-F. Lau and B. Liu \cite[Theorem 3.3]{kll}.  Similar results hold for the case of 
$\GSpin_{2n+1} \times \GL_m$ with $n < m$ as in \cite[\S8.3]{soudry-mem}, where it is worked out for the 
odd special orthogonal groups, and similar to \cite{kll}. (There is typo in \cite[\S8.3]{soudry-mem} where (1.2.4) 
is for $\ell < n$ in the notation of that paper.)  
To complete the proof in the case of special orthogonal groups Kaplan uses a formal identity of power series \cite[(4)]{kaplan3} 
which follows from \cite{bff}, particularly its Appendix, and one would use an analog of that for $\GSpin$ groups. 
\end{proof}

%%%%%%%%%%%%%% NEW SECTION %%%%%%%%%%%%%%%%%%%% 
%%%%%%%%%%%%%%%%%%%%%%%%%%%%%%%%%%%%%%%%%%%%%%% 
\section{Global $L$-functions}\label{sec:L}

We now state the major consequence of the above discussions in the global setting, which we need. 
We use the notation of the earlier sections.

Let $G$ and $H$ be as in (case A), resp. (case B), of the Table in Section \ref{sec-euler-exp}.  
In (case A) consider $H \times \GL_n$ as the Levi of a maximal parabolic inside a larger $\GSpin$ 
group of the same type and similarly in (case B) consider $G \times \GL_m$ as the Levi of a maximal 
parabolic inside the larger $\GSpin$ group.  This is the setup for the Langlands-Shahidi method and 
in the cases we are considering the adjoint action of the dual of the Levi on the unipotent radical of the dual parabolic 
decomposes into two irreducible components.  The first component is the tensor product representation, 
leading to the Rankin-Selberg $L$-functions. The second component, $\varrho$, will be either the twisted symmetric 
square or the twisted exterior square representation of the complex $\GL_n \times \GL_1$, resp. $\GL_m \times \GL_1$. 
To be more precise, we have 
% % % % 
\begin{eqnarray} \label{varrho} 
\varrho = 
\begin{cases} 
\wedge^2 \otimes \omega, & 
\mbox{ (case A) or (case B) , odd,}
\\
%&\\
\sym^2\otimes\omega, & 
\mbox{ (case A) or (case B), even, split,} 
\\ 
%&\\
\wedge^2\otimes\omega, & 
\mbox{ (case A) or (case B), even, quasi-split,} 
\end{cases} 
\end{eqnarray}
% % % % 
where $\omega$ denotes the character on $\GL_1$.

% % % % % % % % % % % % % % % 
\begin{thm} \label{rs-Lqt}
Let $\pi$ be a unitary, cuspidal, globally generic, automorphic representation of $H(\A)$ in (case A), 
resp. of $G(\A)$ in (case B), and 
let $\tau$ be a unitary, cuspidal representation of $\GL_n(\A)$, resp. of $\GL_m(\A)$. 
Let $\omega$ be as in \eqref{omega} (essentially the central character  $\omega_\pi$ of $\pi$).  
Let $S$ be a sufficiently large finite set of places, 
including all the archimedean places, such that for $v \not\in S$ all data are unramified. 
Then we have 
\begin{equation} \label{Lqt} 
\L(\varphi, f_s) = 
\displaystyle
\frac{L^S(s,\pi \times \tau)}{L\left(2s, \tau, \varrho \right)} \cdot R(s),  
\end{equation} 
where, $R(s)$ is a meromorphic function, which can be made holomorphic and nonzero in 
a neighborhood of any given $s=s_{0}$ for an appropriate choice of $\varphi$ and $f_s$ as in 
\eqref{fs-I} or \eqref{fs-II}. 
Here,  
$\L(\varphi, f_s)$ is as in Theorem \ref{thm-BasicIdentityI} or Theorem \ref{thm-BasicIdentityII}, 
as appropriate, and $\varrho$ is as in \eqref{varrho}
\end{thm}
% % % % % % % % % % % % % % % 

%
\begin{proof}
The theorem follows from Theorems \ref{unram-A} and \ref{unram-B} if we take $R(s)$ to be equal 
to the product of the local zeta integrals (\ref{xi}) over $v \in S$. The fact that 
$R(s)$ is meromorphic is clear. To show that it can be made holomorphic in 
the neighborhood of any point $s=s_{0}$ the argument in \cite[\S\S 6-7]{soudry-mem} 
applies.  The presence of the nontrivial center in the $\GSpin$ case does not have an impact 
on those arguments.  
\end{proof}

As a corollary we obtain the following result, which is the precise statement we already used 
in an earlier work as \cite[Prop. 4.9]{manuscripta}.

%%%%%%%%%% 
\begin{prop} \label{L-pole}
Let $\pi$ be globally generic unitary cuspidal automorphic representations of $\GSpin_{2n+1}(\A)$ or $\GSpin_{2n}(\A)$ and 
let $\tau$ be a unitary cuspidal automorphic representation of $\GL_n(\A)$, as above. 
If $L^S(s,\pi\times\tau)$ has a pole at $s_0$ with $\Re(s_0) \ge 1$, then for a choice of $f_s$ the Eisenstein series 
$E(h, f_s)$ has a pole at $s = s_0.$ 
\end{prop}
%%%%%%%%%% 

%%%%%%%%%% 
\begin{proof} 
This statement follows immediately from Theorem \ref{rs-Lqt}.  Recall that any pole of  
$\L(\varphi, f_s)$ must come from a pole of the Eisenstein series that is used to define it. 
Moreover, the twisted symmetric and twisted exterior square $L$-functions appearing on the 
right hand sides in Theorem \ref{rs-Lqt} are holomorphic in $\Re(s) \ge 1$ so they 
can not cancel a possible pole of $L(s, \pi \times \tau)$.  

We note that it follows from \cite[Appendix]{ckpss2} that the above argument holds for generic 
representations with respect to an arbitrary $\psi$ and any particular ``standard'' one that we 
may have fixed.   This issue is relevant when the some of our groups are not of adjoint type, as here. 
See \cite[Appendix]{ckpss2} for more details. 
\end{proof} 
%%%%%%%%%% 

%%%%%%%%%%%%%% NEW SECTION %%%%%%%%%%%%%%%%%%%% 
%%%%%%%%%%%%%%%%%%%%%%%%%%%%%%%%%%%%%%%%%%%%%%% 

\end{document}